\documentclass[DIV=15, bibliography=totoc]{scrartcl}  

\usepackage[a4paper,margin=2.25cm,bottom=2.5cm]{geometry}

\usepackage{xcolor} 
\definecolor{lightblue}{rgb}{0,0.5,1.0}
\definecolor{linkblue}{rgb}{0,0.1,0.6}
\definecolor{citegreen}{rgb}{0,0.4,0.0}
\definecolor{linkred}{rgb}{0.8,0,0.005}
\definecolor{mailviolet}{rgb}{0.3,0,0.35}
\definecolor{tumblue}{rgb}{0,0.396,0.741}
\definecolor{darkgreen}{rgb}{0,0.4,0} 
\definecolor{darkbrown}{rgb}{0.5, 0.396, 0.09}

\usepackage[hyperindex,colorlinks=true,linkcolor=linkblue,citecolor=citegreen,urlcolor=mailviolet,filecolor=linkred]{hyperref}

\usepackage{caption, subcaption}
\usepackage{fancyhdr}

\usepackage{siunitx}
\usepackage[square,comma,nonamebreak,compress,numbers]{natbib}

\usepackage{soulpos}
\usepackage{ulem}
\usepackage{authblk} 

\usepackage[inkscapearea=page]{svg}
\usepackage{amsmath}
\usepackage[noabbrev]{cleveref}
\usepackage{upgreek}
\usepackage{adjustbox}
\usepackage{color}
\usepackage{enumitem}

\usepackage{pgfplots}
\usepackage{pgfplotstable}
\usepackage{filecontents}
\usepackage{tikz}
\usepackage{tikzscale}
\usepackage{tikz-3dplot}

\usetikzlibrary{spy}
\usetikzlibrary{intersections}
\usetikzlibrary{arrows,shapes}
\usetikzlibrary{spy}
\usetikzlibrary{backgrounds}  
\usetikzlibrary{decorations}
\usetikzlibrary{decorations.markings}
\usetikzlibrary{positioning}
\usetikzlibrary{patterns}
\usetikzlibrary{calc} 
\usetikzlibrary{quotes} 
\usetikzlibrary{external} 
\usetikzlibrary{matrix}

\pgfplotscreateplotcyclelist{customCycleList}
{%
  {blue, mark=o},
  {red, mark=square},
  {darkgreen, mark=diamond},
  {cyan, mark=diamond},
  {darkbrown, mark=triangle*},
  {black, mark=pentagon},
}

\pgfplotsset{every axis/.append style= {
    cycle list name=customCycleList,
}}

\title{On the Use of Neural Networks for Full Waveform Inversion}

\author[1]{Leon Herrmann\thanks{\href{mailto:leon.herrmann@tum.de}{\texttt{leon.herrmann@tum.de}},
    Corresponding author}}
\author[1]{Tim Bürchner}
\author[2]{Felix Dietrich}
\author[1]{Stefan Kollmannsberger}

\affil[1]{Chair of Computational Modeling and Simulation, Technical University of Munich, School of Engineering and Design, Arcisstraße 21, Munich, 80\,333, Germany}

\affil[2]{Chair of Scientific Computing in Computer Science (SCCS), Technical University of Munich, School of Computation, Information and Technology, Boltzmannstraße 3, Garching, 85\,748, Germany}

\newcommand{\publicationDate}{\today}

\date{}
\fancypagestyle{plain}{%
\fancyhf{}%
\fancyfoot[L]{
\vspace{-0.0cm} 
\footnotesize{Preprint submitted to \textit{Engineering with Computers}.  \hfill
\publicationDate
\\
}} }

\begin{document}     

\normalem \maketitle  
\normalfont\fontsize{11}{13}\selectfont

\vspace{-1.5cm} \hrule 

\section*{Abstract}

Neural networks have recently gained attention in solving inverse problems. One prominent methodology are Physics-Informed Neural Networks (PINNs) which can solve both forward and inverse problems. In the paper at hand, full waveform inversion is the considered inverse problem. The performance of PINNs is compared against classical adjoint optimization, focusing on three key aspects: the forward-solver, the neural network Ansatz for the inverse field, and the sensitivity computation for the gradient-based minimization. Starting from PINNs, each of these key aspects is adapted individually until the classical adjoint optimization emerges. It is shown that it is beneficial to use the neural network only for the discretization of the unknown material field, where the neural network produces reconstructions without oscillatory artifacts as typically encountered in classical full waveform inversion approaches. Due to this finding, a hybrid approach is proposed. It exploits both the efficient gradient computation with the continuous adjoint method as well as the neural network Ansatz for the unknown material field. This new hybrid approach outperforms Physics-Informed Neural Networks and the classical adjoint optimization in settings of two and three-dimensional examples.

\vspace{0.25cm}
\noindent \textit{Keywords:} 
physics-informed neural networks, full waveform inversion, deep learning, adjoint state method, inverse problems, finite difference method
\vspace{0.25cm}


\newpage 
\textbf{List of Symbols}
\begin{table}[h!]

	\centering

	\begin{tabular}{llll}

		$\boldsymbol{x}$ & coordinates 	& $\phi_0$ & source amplitude \\

		$t$ & time & $\boldsymbol{x}_s$ & source position\\

		$\gamma(\boldsymbol{x}, t)$ & indicator function & $\boldsymbol{m}$ & indicator parametrization\\

		$u(\boldsymbol{x}, t)$ & wavefield & $u_{\mathcal{M}}(\boldsymbol{x}_i, t)$ & sensor measurements\\

		$\hat{\gamma}(\boldsymbol{x})$ & indicator prediction & $\mathcal{N}[u, \gamma]$ & PDE parametrization\\

		$\hat{u}(\boldsymbol{x}, t)$ & wavefield prediction & $\mathcal{B}[u]$ & boundary condition operator \\

		$\boldsymbol{\gamma}$ & indicator at $N$ grid points & $\mathcal{I}[u]$ & initial condition operator\\

		$\boldsymbol{u}$ & wavefield at $N$ grid points & 	$N_{\mathcal{M}}$ & number of sensors\\

		$A_{\gamma}(\boldsymbol{\theta};\boldsymbol{x})$ & indicator neural network & $N_\mathcal{N}$ & number of internal grid points\\

		$A_{u}(\boldsymbol{\theta};\boldsymbol{x})$ & wavefield neural network & $N_\mathcal{B}$ & number of boundary grid points\\

		$\boldsymbol{\theta}$ & neural network parameters & $N_\mathcal{I}$ & number of inital grid points \\

		$C_{\gamma}(\boldsymbol{c}, \boldsymbol{x})$ & constant Ansatz & $\mathcal{L}$ & cost function\\

		$c_i, \boldsymbol{c}$ & constant Ansatz coefficients & $\mathcal{L}_{\mathcal{M}}$ & measurement loss \\

		$N_i(x)$ & constant Ansatz shape functions & $\mathcal{L}_\mathcal{N}$ & PDE loss\\

		$F[\boldsymbol{\gamma}, \boldsymbol{x}, t]$ & non-trainable forward operator & $\mathcal{L}_\mathcal{B}$ & boundary condition loss\\

		$\rho(\boldsymbol{x}), \rho_0$ & density & $\mathcal{L}_\mathcal{I}$ & initial condition loss\\

		$c(\boldsymbol{x}), c_0$ & wave velocity & $\boldsymbol{\lambda}, \lambda$ & learnable penalty weights\\

		$f(\boldsymbol{x}, t)$ & volume force & $\alpha$ & step size, learning rate\\

		$\phi(t)$ & sine burst source & $\epsilon$ & small positive value $<<1$ \\

		$n_c$ & number of cycles & $K_\gamma$ & Fréchet kernel\\

		$\omega, f_{\phi}$ & source frequency & $u^{\dagger}$ & adjoint wavefield \\	

	\end{tabular}

\end{table}

\section{Introduction}\label{sec:introduction}

\subsection{Motivation}\label{sec:motivation}

Since the emergence of Physics-Informed Neural Networks (PINNs)~\cite{lagaris_artificial_1998, psichogios_hybrid_1992, raissi_physics-informed_2019}, they have been of great interest to the scientific community of computational mechanics. Despite their novelty, multiple review articles have already emerged~\cite{cuomo_scientific_2022, karniadakis_physics-informed_2021,Kollmannsberger2021a}. The essence of PINNs is a variational approach, where a solution field is approximated with a neural network whose weights and biases are learned by minimizing the residual of a differential equation evaluated at collocation points. This general formulation allows for solving a large variety of problems that can at least partially be described by differential equations. At the same time, current machine learning frameworks like PyTorch~\cite{paszke_pytorch_2019} and TensorFlow \cite{tensorflow2015-whitepaper} provide generic tools to construct PINNs with very little additional and rather simple code compared to the complexity of classical differential equation solvers. \\

The method shows potential in settings where the amount of data available is insufficient to build supervised surrogate models, as e.g. presented in ~\cite{thuerey_deep_2020, zhu_physics-constrained_2019, lino_simulating_2021, sanchez-gonzalez_learning_2020, pfaff_learning_2021, bhatnagar_prediction_2019, lu_learning_2021}. If used as a forward solver, PINNs are often computationally more expensive for many common differential equations which are efficiently solvable with established forward solvers~\cite{markidis_old_2021, leiteritz_how_2021}, as e.g. finite element~\cite{hughes_finite_2000} or finite difference methods~\cite{langtangen_finite_2017}. However, in settings where classical methods show difficulties, PINNs are promising. Possible applications are high-dimensional differential equations~\cite{han_solving_2018, sirignano_dgm_2018}, surrogate modeling~\cite{goswami_physics-informed_2022, oldenburg_geometry_2022, wong_improved_2021}, accelerating traditional solvers~\cite{markidis_old_2021}, and inverse problems~\cite{shukla_physics-informed_2020, rashtbehesht_physicsinformed_2022, cai_flow_2021, wang_deep_2021, mao_physics-informed_2020, jagtap_conservative_2020, jagtap_physics-informed_2022, chen_physics-informed_2020}. \\

The novelty in this paper is twofold. Firstly, it pinpoints the strengths and weaknesses of PINNs for inverse problems in comparison to classical full waveform inversion where the continuous adjoint state method~\cite{givoli_tutorial_2021, plessix_review_2006} is considered as the state-of-the-art reference method. To this end, a thorough investigation of the strengths and weaknesses of the individual components of each method is carried out by going from PINNs to the adjoint optimization step by step. The investigation is conducted on the scalar wave equation, popular in non-destructive testing~\cite{fichtner_full_2011, sayag_shape_2022}, and recently explored with PINNs in~\cite{shukla_physics-informed_2020, rashtbehesht_physicsinformed_2022}. In particular, the following three aspects will be examined extensively

\begin{enumerate}[label=(\alph*)]

	\item the use of PINNs as a forward-solver,

	\item neural networks as an Ansatz for the coefficient field of the PDE,

	\item gradient computation.

\end{enumerate}

This dissection gives rise to the methods summarized in~\cref{tab:methods} which are described and studied in detail in~\cref{sec:effect} apart from the more general setting of PINNs described in the introductory~\cref{sec:pinns}.


\begin{table}[htb]

	\caption{Investigated methods.} \label{tab:methods}

	\centering

	\begin{tabular}{lll}

		\hline

		\textbf{Method}                                                                             & \textbf{Ansatz}                                                                                            & \textbf{Gradient Computation}                                                                                                                                                                                     \\ \hline

		\begin{tabular}[c]{@{}l@{}} PINNs \\   (\cref{sec:pinns}) \end{tabular}      & \begin{tabular}[c]{@{}l@{}}

			$\hat{\gamma}(\boldsymbol{x})=A_{\gamma}(\boldsymbol{\theta_{\gamma}};\boldsymbol{x})$ \\

			$\hat{u}(\boldsymbol{x}, t)=A_{u}(\boldsymbol{\theta_{u}};\boldsymbol{x}, t)$

		\end{tabular} & $\nabla_{\boldsymbol{\theta_{\gamma}}}\mathcal{L}$ with automatic differentiation                                                                                                                        \\ \hline

		\begin{tabular}[c]{@{}l@{}}PINNs with non-trainable\\ forward operator \\ (\cref{sec:PINNSnonLearnableFO}) \end{tabular}   &

		\begin{tabular}[c]{@{}l@{}}

			$\hat{\gamma}(\boldsymbol{x})=A_{\gamma}(\boldsymbol{\theta_{\gamma}}; \boldsymbol{x})$ \\ $\hat{u}(\hat{\gamma},\boldsymbol{x}, t)=F[\boldsymbol{\hat{\gamma}},\boldsymbol{x}, t]$ 

		\end{tabular}    & $\nabla_{\boldsymbol{\theta_{\gamma}}}\mathcal{L}$ with automatic differentiation                                                                                                                        \\ \hline

		\begin{tabular}[c]{@{}l@{}} Adjoint Method \\   (\cref{sec:AdjointMethod}) \end{tabular}                                                                     & 		\begin{tabular}[c]{@{}l@{}}

			$\hat{\gamma}(\boldsymbol{x})=C_{\gamma}(\boldsymbol{c},\boldsymbol{x})$\\ 

			$\hat{u}(\hat{\gamma},\boldsymbol{x}, t)=F[\boldsymbol{\hat{\gamma}},\boldsymbol{x}, t]$

		\end{tabular}      & \begin{tabular}[c]{@{}l@{}}$\nabla_{\boldsymbol{c}}\mathcal{L}$ with automatic differentiation\\ or $\nabla_{\boldsymbol{c}}\mathcal{L}$ with the adjoint method\end{tabular}                            \\ \hline \hline

		\begin{tabular}[c]{@{}l@{}} Hybrid Approach \\   (\cref{sec:hybridApproach}) \end{tabular}                                                                    & \begin{tabular}[c]{@{}l@{}}

			$\hat{\gamma}(\boldsymbol{x})=A_{\gamma}(\boldsymbol{\theta_{\gamma}},\boldsymbol{x})$ \\ $\hat{u}(\hat{\gamma},\boldsymbol{x}, t)=F[\boldsymbol{\hat{\gamma}},\boldsymbol{x}, t]$ 

		\end{tabular}         & 

		\begin{tabular}[c]{@{}l@{}} $\nabla_{\boldsymbol{\hat{\gamma}}} \mathcal{L}$ with the adjoint method \\ and $\nabla_{\boldsymbol{\theta_{\gamma}}} \hat{\gamma}$ with automatic differentiation  \end{tabular} \\ \hline

	\end{tabular}

\end{table}

Most of these methods naturally emerge as individual steps moving from a fully PINN-based approach to a classical adjoint optimization. By contrast, and as a second novelty, the paper at hand presents a hybrid method, which combines the strengths of the investigated approaches and delivers more accurate results.


\subsection{Full Waveform Inversion}\label{sec:fwi}

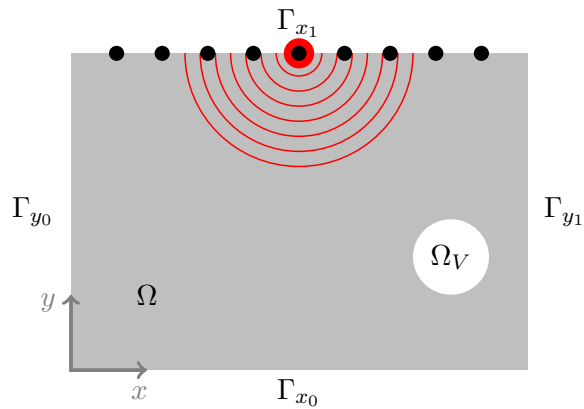
\begin{figure}[htb]

	\centering

	\begin{tikzpicture}

	\fill [lightgray, line width=0.3mm] (0,0) rectangle (6,4.2);

	\fill [line width=0.3mm, white] (5,1.5) circle (0.5cm);

	\node at (5,1.5) {$\Omega_V$};

	\begin{scope}

	\clip (0,0) rectangle (6,4.2);

	\foreach {\r} in {0.3,0.5,...,1.5} {

		\draw [line width=0.2mm, red] (3,4.2) circle (\r cm);

	}

	\end{scope}

	\draw [gray, line width=0.5mm,<->] (0,1) -- (0,0) -- (1,0);

	\node [gray] at (-0.3,0.9) {$y$};

	\node [gray] at (0.9, -0.3) {$x$};

	\fill [red] (3,4.2) circle (0.2cm);


	\foreach {\x} in {0.6,1.2,...,5.4} {


		\foreach {\y} in {4.2} {

			\fill [black] (\x,\y) circle (0.1cm);

		}

	}






	\node at (-0.5,2.1) {$\Gamma_{y_0}$};

	\node at (6.5,2.1) {$\Gamma_{y_1}$};

	\node at (3,-0.3) {$\Gamma_{x_0}$};

	\node at (3,4.6) {$\Gamma_{x_1}$};

	\node at (1,1) {$\Omega$};

	\end{tikzpicture}

	\caption{Ultrasonic Nondestructive Testing on a two-dimensional rectangular domain $\Omega$ with a void $\Omega_V$. Horizontal boundaries are defined as $\Gamma_{x}=\Gamma_{x_0}\cup\Gamma_{x_1}$, while vertical boundaries are defined as $\Gamma_{y}=\Gamma_{y_0}\cup\Gamma_{y_1}$.}

	\label{fig:ndt}

\end{figure}

Full waveform inversion aims to non-destructively estimate the material distribution in a specimen. A conceptual illustration is provided in~\cref{fig:ndt} with a rectangular domain $\Omega$ and a circular void $\Omega_V$, that is to be detected. Voids are identified using a fixed number of sources emitting a signal (red circles) and sensors (black circles) recording the image of the signal.

The physics of this problem can be modeled with the scalar wave equation in an isotropic heterogeneous medium, which for the solution field $u(\boldsymbol{x},t)$, the density $\rho(\boldsymbol{x})$, the wave speed $c(\boldsymbol{x})$ and the volume force $f(\boldsymbol{x},t)$ is given as
\begin{align}
\rho(\boldsymbol{x}) u_{tt}(\boldsymbol{x}, t)-\nabla\cdot(\rho(\boldsymbol{x}) c(\boldsymbol{x})^2 \nabla u(\boldsymbol{x}, t)) - f(\boldsymbol{x},t) = 0 \qquad \text{on } \Omega \times\mathcal{T}. \label{eq:waveequationOG}
\end{align}
Additionally, homogeneous Neumann boundary conditions and homogeneous initial conditions are prescribed:
\begin{alignat}{2}
u_x(\boldsymbol{x}, t)&=0  && \qquad \text{on }  \Gamma_{y}\times\mathcal{T} \label{eq:neumanninx} \\ 
u_y(\boldsymbol{x}, t)&=0  && \qquad \text{on } \Gamma_{x}\times\mathcal{T} \label{eq:neumanniny} \\
u(\boldsymbol{x}, 0)&=u_t(\boldsymbol{x}, 0)=0 && \qquad \text{on } \Omega \label{eq:initialconditions} 
\end{alignat}
Parametrizing the inverse problem with an indicator function $\gamma$ by scaling the density leads to the most accurate reconstruction of voids and unknown homogeneous Neumann boundaries. The constant intact material has the density $\rho_0$ and $\gamma$ is constrained between $[\epsilon, 1]$, where $\epsilon \ll 1$ is a non-zero lower bound:
\begin{align}
\rho(\boldsymbol{x})&=\gamma(\boldsymbol{x}) \rho_0 \label{eq:rhoscaling}
\end{align}
In contrast, the wave speed is assumed to be constant $c(\boldsymbol{x})=c_0$. With this parametrization, the scalar wave equation takes the form:
\begin{align}
\gamma(\boldsymbol{x}) \rho_0 u_{tt}(\boldsymbol{x}, t)-\nabla\cdot(\gamma(\boldsymbol{x}) \rho_0 c_0^2 \nabla u(\boldsymbol{x}, t)) - f(\boldsymbol{x},t) = 0\qquad \text{on } \Omega \times\mathcal{T}. \label{eq:waveequation}
\end{align}
The volume force $f(\boldsymbol{x},t)$ emits the signal at the position $\boldsymbol{x}_s$. This is usually modeled by a spatial Dirac delta using a source term $\psi(t)$.
\begin{equation}
f(\boldsymbol{x},t)=\psi(t)\delta(\boldsymbol{x}-\boldsymbol{x}_s) \label{eq:source}
\end{equation}
A typical source in the context of non-destructive testing is the sine burst with $n_c$ cycles, a frequency $\omega = 2 \pi f_{\psi}$, and an amplitude $\psi_0$, so that
\begin{equation}
\psi(t)=\begin{cases}
\psi_0 \sin(\omega t)\sin(\frac{\omega t}{2 n_c}), & \text{for } 0\leq t\leq \frac{2\pi n_c}{\omega} \\
0, & \text{for } \frac{2\pi n_c}{\omega}<t. \label{eq:sineburst}
\end{cases}
\end{equation}
This burst type is used in all examples in this paper.

The goal of the inversion is to reconstruct a material distribution, i.e. the indicator $\hat{\gamma}(\boldsymbol{x}; \boldsymbol{m})$ parametrized by $\boldsymbol{m}$, that best reproduces the $N_{\mathcal{M}}$ sensor measurements
$u_{\mathcal{M}}(\boldsymbol{x}_i, t)$ during a forward simulation which yields $\hat{u}(\boldsymbol{x}, t; \boldsymbol{m})$. The inversion is posed as an optimization problem, where the cost function is the mean squared error between the prediction $\hat{u}$ and the measurements $u_{\mathcal{M}}$:
\begin{equation}
\mathcal{L}_{\mathcal{M}}(\hat{\boldsymbol{u}}, \boldsymbol{u}_{\mathcal{M}})=\frac{1}{2}\int_\mathcal{T}\int_\Omega \sum_{i=1}^{N_{\mathcal{M}}} \biggl(\hat{u}(\boldsymbol{x}, t; \boldsymbol{m}) - u_{\mathcal{M}}(\boldsymbol{x}_i, t)\biggr)^2 \delta(\boldsymbol{x} - \boldsymbol{x}_i) d\Omega d\mathcal{T}. \label{eq:measurementloss}
\end{equation}
The gradient of $\mathcal{L}_{\mathcal{M}}$ w.r.t. the parameters $\boldsymbol{m}$ is expressed by $\nabla_m \mathcal{L}_{\mathcal{M}}(\hat{\boldsymbol{u}}(\boldsymbol{m}), \boldsymbol{u}_{\mathcal{M}})$ and computed with the adjoint state method. The gradient is then used in a gradient-based optimization procedure to minimize 
$\mathcal{L}_{\mathcal{M}}$ by iteratively updating the parameters $\boldsymbol{m}$:  
\begin{equation}
\boldsymbol{m}^{(j+1)} = \boldsymbol{m}^{(j)} - \alpha \nabla_{m^{(j)}} \mathcal{L}_{\mathcal{M}} \label{eq:gradientDecent}
\end{equation}
in which $\alpha$ is the step size and the superindex $(j)=1,2,\dots$ is the iteration counter. This iterative process then yields an optimized estimate of the indicator $\hat{\gamma}(\boldsymbol{x}; \boldsymbol{m})$. In case of multiple sources, an average of the mean squared errors $\mathcal{L}_{\mathcal{M}}$ and their gradients $\nabla_m \mathcal{L}_{\mathcal{M}}$ must be considered in the iteration update from~\cref{eq:gradientDecent}. More details on the computation of the gradient using the adjoint method are found in~\cite{buerchner_2022}. For a more in depth explanation of the adjoint method for inverse problems, see Givoli~\cite{givoli_tutorial_2021}.




\subsection{Physics-Informed Neural Network}\label{sec:pinns}

The idea of PINNs is to use an artificial neural network $A$ to approximate either the solution $\hat{u}$ or and/or the coefficients of the differential equation $\hat{\gamma}$ using a collocation type approach. In the paper at hand,~\cref{eq:waveequation} is considered where the solution $\hat{u}$ is the wave field and $\rho_0$ is scaled by $\gamma(\boldsymbol{x})$ as given in~\cref{eq:rhoscaling}. In the most general case, both $u$ and $\gamma$ are approximated by a neural network $A=(A_u, A_{\gamma})$:
\begin{align}
\hat{u}&=A_{u}(\boldsymbol{x},t; \boldsymbol{\theta}_u) \label{eq:solution} \\
\hat{\gamma}&=A_{\gamma}(\boldsymbol{x}, t; \boldsymbol{\theta}_{\gamma}) \label{eq:coefficient}
\end{align}
where $\boldsymbol{\theta}=\{\boldsymbol{\theta}_u, \boldsymbol{\theta}_{\gamma}\}$ are the model parameters of the neural network and the hat symbol indicates that the values of $u$ and $\gamma$ are predictions from a neural network. In the remaining paper, in models using only one of the two neural networks, i.e. either $A_{u}$ or $A_{\gamma}$, the notation $\boldsymbol{\theta}$ without subindex is used as the corresponding neural network parameters.

The initial boundary value problem defined by~\cref{eq:waveequation,eq:neumanninx,eq:neumanniny,eq:initialconditions} is then solved by an optimization procedure which is similar to the one used in the classical full waveform inversion described in the previous section. To this end, the initial boundary problem is first cast into its residual form: 
\begin{align}
\mathcal{N}[u; \gamma; \boldsymbol{x}, t] &= \gamma(\boldsymbol{x}) \rho_0 u_{tt}(\boldsymbol{x}, t)-\nabla\cdot(\gamma(\boldsymbol{x}) \rho_0 c_0^2 \nabla u(\boldsymbol{x}, t)) - f(\boldsymbol{x},t) \label{eq:N} \\
\mathcal{B}_x[u; \boldsymbol{x}, t] &= u_x(\boldsymbol{x}, t) \label{eq:Bx}\\
\mathcal{B}_y[u; \boldsymbol{x}, t] &= u_y(\boldsymbol{x}, t) \label{eq:By}\\
\mathcal{I}_0[u; \boldsymbol{x}, t] &= u(\boldsymbol{x}, t) \label{eq:I0}\\
\mathcal{I}_1[u; \boldsymbol{x}, t] &= u_t(\boldsymbol{x}, t). \label{eq:I1}
\end{align}
The operators $\mathcal{N}, \mathcal{B}_x, \mathcal{B}_y, \mathcal{I}_0, \mathcal{I}_1$ are the residuals of the corresponding~\cref{eq:waveequation,eq:neumanninx,eq:neumanniny,eq:initialconditions}. These residuals are evaluated at a number of collocation points $N=N_{\mathcal{N}}+ N_{\mathcal{B}_x}+N_{\mathcal{B}_y}+N_{\mathcal{I}_0}+N_{\mathcal{I}_1}$ using the current approximation of the solution coefficient fields given by~\cref{eq:solution,eq:coefficient}. This requires the computation of the first and second derivatives of~\cref{eq:solution,eq:coefficient} w.r.t $\boldsymbol{x}$ and $t$ which is carried out by either using reverse-mode automatic differentiation~\cite{baydin_automatic_2018} or numerical differentiation~\cite{zhu_physics-constrained_2019,geneva_modeling_2020}.

The loss functional to be minimized by the overall optimization procedure, i.e. the cost function $\mathcal{L}$, is then defined as the weighted sum of squares of the residual evaluations (the differential equation residual $\mathcal{L}_{\mathcal{N}}$, the boundary condition residual $\mathcal{L}_{\mathcal{B}}$, and the initial condition residual $\mathcal{L}_{\mathcal{I}}$) at the collocation points:
\begin{align}
\mathcal{L}_{\mathcal{N}}(\boldsymbol{u},\boldsymbol{\gamma})&=\frac{1}{2N_{\mathcal{N}}}\sum_{i=1}^{N_{\mathcal{N}}}\lambda_{\mathcal{N},i}(\mathcal{N}[u_i; \gamma_i])^2 \label{eq:losspde} \\
\mathcal{L}_{\mathcal{B}}(\boldsymbol{u})&=\frac{1}{2N_{\mathcal{B}_x}}\sum_{i=1}^{N_{\mathcal{B}_x}}\lambda_{\mathcal{B}_x,i}(\mathcal{B}_x[u_i])^2+\frac{1}{2N_{\mathcal{B}_y}}\sum_{i=1}^{N_{\mathcal{B}_y}}\lambda_{\mathcal{B}_y,i}(\mathcal{B}_y[u_i])^2 \label{eq:lossbc}\\
\mathcal{L}_{\mathcal{I}}(\boldsymbol{u})&=\frac{1}{2N_{\mathcal{I}_0}}\sum_{i=1}^{N_{\mathcal{I}_0}}\lambda_{\mathcal{I}_0,i}(\mathcal{I}_0[u_i])^2+\frac{1}{2N_{\mathcal{I}_1}}\sum_{i=1}^{N_{\mathcal{I}_1}}\lambda_{\mathcal{I}_1,i}(\mathcal{I}_1[u_i])^2 \label{eq:lossic}\\
\mathcal{L}(\boldsymbol{u},\boldsymbol{\gamma})&=\lambda_{\mathcal{N}}\mathcal{L}_{\mathcal{N}}(\boldsymbol{u},\boldsymbol{\gamma})+\lambda_{\mathcal{B}}\mathcal{L}_{\mathcal{B}}(\boldsymbol{u})+\lambda_{\mathcal{I}}\mathcal{L}_{\mathcal{I}}(\boldsymbol{u}) \label{eq:cost}	
\end{align}
To counteract imbalances in the magnitude of the loss terms $\mathcal{L}_{\mathcal{N}}, \mathcal{L}_{\mathcal{B}}, \mathcal{L}_{\mathcal{I}}$ usually each loss term is multiplied with a weight $\lambda$ individual to that loss term. The penalty weights can be computed automatically using attention mechanisms from deep vision~\cite{wang_residual_2017,zhang_occluded_2018}, where the penalty weights are updated via a maximization w.r.t. the cost function $\mathcal{L}$ resulting in a minimax optimization procedure. Expanding on this idea, as discussed in~\cite{nandwani_primal-dual_2019, mcclenny_self-adaptive_2022}, each evaluation at a collocation point is considered an individual equality constraint, thereby allocating one penalty weight per collocation point, as shown in the loss equations~\eqref{eq:losspde},~\eqref{eq:lossbc},~\eqref{eq:lossic}. This extension makes it possible to assign a greater emphasis on more important collocation points, i.e. points which lead to larger residuals. This yields an overall minimax optimization of the form 
\begin{equation}
\min_{\boldsymbol{\theta}} \max_{\boldsymbol{\lambda}} \mathcal{L}(\boldsymbol{\theta}, \boldsymbol{\lambda}), \label{eq:minimax}
\end{equation} 
This minimax optimization introduces the learning rate of the maximization $\alpha_{\max}$ as an additional hyperparameter to the learning rate for the minimization $\alpha_{\min}$. 

The optimization problem posed in~\cref{eq:minimax} is solved using the gradients \textbf{$\nabla_{\boldsymbol{\theta}} \mathcal{L}$} and \textbf{$\nabla_{\boldsymbol{\lambda}} \mathcal{L}$} obtained via reverse-mode automatic differentiation. The Adam optimizer~\cite{kingma_adam_2017} is used in the paper at hand as an optimization procedure in all examples. However, more advanced optimizers such as combinations of Adam and L-BFGS~\cite{liu_limited_1989} often yield a faster convergence (see, e.g. ~\cite{samaniego_energy_2020}).  

Important differences exist between the application of PINNs to either forward or inverse problems. These differences are summarized in~\cref{tab:approaches} using three problem categories: (a) Forward problems, (b) Inverse problems with full domain knowledge and (c) Inverse problems with partial domain knowledge.  
\begin{table}[htb]
	\caption{Overview of differences in the PINN formulation between forward and inverse problems.} \label{tab:approaches}
	\centering
	\begin{tabular}{llll}
		\hline
		\textbf{Problem}   & \textbf{Cost Function} &  \textbf{Neural Networks}                                                                     \\ \hline 
		\begin{tabular}[c]{@{}l@{}}\textbf{(a)} Forward Problem\\
		\end{tabular}	
		& \begin{tabular}[c]{@{}l@{}}$\mathcal{L}(\boldsymbol{\hat{u}}, \boldsymbol{\gamma})$\\$=\mathcal{L}_{\mathcal{N}}(\hat{\boldsymbol{u}},\boldsymbol{\gamma})+\mathcal{L}_{\mathcal{B}}(\hat{\boldsymbol{u}})+\mathcal{L}_{\mathcal{I}}(\hat{\boldsymbol{u}})$\end{tabular}  &$\hat{u}_i=A_u(\boldsymbol{x}_i, t_i; \boldsymbol{\theta})$                                \\ \hline
		\begin{tabular}[c]{@{}l@{}}\textbf{(b)} Inverse Problem\\ (Full Domain Knowledge)\\ \end{tabular}    & \begin{tabular}[c]{@{}l@{}}$\mathcal{L}(\boldsymbol{u_\mathcal{M}}, \boldsymbol{\hat{\gamma}})$\\$=\mathcal{L}_{\mathcal{N}}(\boldsymbol{u}_{\mathcal{M}},\hat{\boldsymbol{\gamma}})$\end{tabular}                                                                 &  $\hat{\gamma}_i=A_{\gamma}(\boldsymbol{x}_i, t_i, \boldsymbol{\theta})$                                                                                                                 \\ \hline
		\begin{tabular}[c]{@{}l@{}}\textbf{(c)} Inverse Problem\\(Partial Domain Knowledge)\\
		\end{tabular} & \begin{tabular}[c]{@{}l@{}}$\mathcal{L}(\boldsymbol{\hat{u}},\boldsymbol{u_\mathcal{M}}, \boldsymbol{\hat{\gamma}})$\\$=\mathcal{L}_{\mathcal{M}}(\boldsymbol{\hat{u}}, \boldsymbol{u}_{\mathcal{M}})+\mathcal{L}_{\mathcal{N}}(\boldsymbol{\hat{u}}, \hat{\boldsymbol{\gamma}})$\\$+\mathcal{L}_{\mathcal{B}}(\boldsymbol{\hat{u}})+\mathcal{L}_{\mathcal{I}}(\boldsymbol{\hat{u}})$\end{tabular} & \begin{tabular}[c]{@{}l@{}}$\hat{u}_i=A_u(\boldsymbol{x}_i, t_i; \boldsymbol{\theta}_u)$\\ and $\hat{\gamma}_i=A_{\gamma}(\boldsymbol{x}_i, t_i; \boldsymbol{\theta}_\gamma)$\end{tabular} \\ \hline
	\end{tabular}
\end{table}

\textbf{Problem type (a) - Forward Problem}: In the forward problem, only the neural network $A_u$ from~\cref{eq:solution} is used as an approximator of the wave field $u(\boldsymbol{x}, t)$, as $\gamma(\boldsymbol{x},t)$ is known. The gradients of the approximated field $\hat{u}$ are computed w.r.t. the inputs $\boldsymbol{x}$ and $t$, which enables the computation of the cost function $\mathcal{L}(\boldsymbol{\hat{u}},\boldsymbol{\gamma})$ from~\cref{eq:cost}. The wavefield is then solved by performing the minimax optimization from~\cref{eq:minimax}. Variations of PINNs for the forward problem of the wave equation can be found in literature~\cite{moseley_solving_2020, song_solving_2021, song_solving_2021-1, karimpouli_physics_2020}. 

With a minor modification, the framework of PINNs can also be applied to the inverse problem, i.e. the identification of the spatially varying indicator function $\gamma(\boldsymbol{x})$. 

\textbf{Problem type (b) - Inverse problems with full domain knowledge of the solution field} $u(\boldsymbol{x}, t)$:  Since $u(\boldsymbol{x}, t)$ is known in the form of measurements $\boldsymbol{u}_\mathcal{M}$ everywhere, a neural network to predict $A_u$ is not needed. Instead, $\boldsymbol{u}_\mathcal{M}$ can directly be used to evaluate the cost function $\mathcal{L}$ from~\cref{eq:cost}. Since the initial and boundary conditions are already approximately satisfied by the wave field, the cost function $\mathcal{L}(\boldsymbol{u}_\mathcal{M},\boldsymbol{\hat{\gamma}})$ is now solely composed of the residual term of the differential equation $\mathcal{L}_{\mathcal{N}}(\boldsymbol{u}_\mathcal{M},\boldsymbol{\hat{\gamma}})$. Therein, a neural network $A_\gamma$ is used to predict $\boldsymbol{\hat{\gamma}}$ and $u(\boldsymbol{x}, t)$ is directly plugged into $\mathcal{L}_{\mathcal{N}}(\boldsymbol{u}_\mathcal{M},\boldsymbol{\hat{\gamma}})$ as $\boldsymbol{u}_\mathcal{M}$. A detailed description of solving this inverse problem with PINNs can be found in~\cite{shukla_physics-informed_2020} which documents the state-of-the-art.

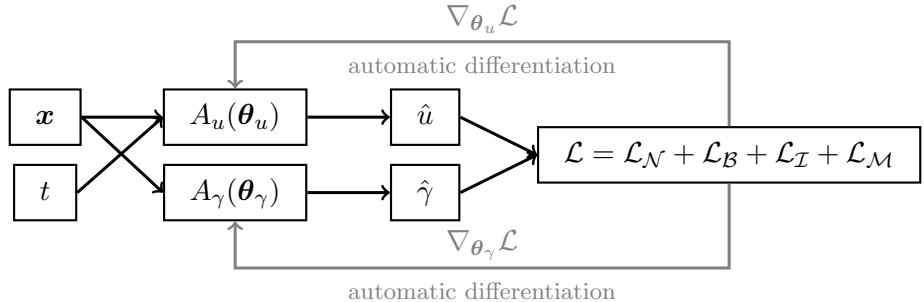
\begin{figure}[htb]
	\centering
	\begin{tikzpicture}
	\node (I1) [draw, thick] at (-1,0) {\begin{tabular}{c} $\boldsymbol{x}$ \end{tabular}};
	\node (I2) [draw, thick] at (-1,-1) {\begin{tabular}{c} $t$ \end{tabular}};
	\node (A1) [draw, thick] at (1.5,0) {\begin{tabular}{c} $A_u(\boldsymbol{\theta}_u)$ \end{tabular}};
	\node (A2) [draw, thick] at (1.5,-1) {\begin{tabular}{c} $A_{\gamma}(\boldsymbol{\theta}_\gamma)$ \end{tabular}};
	\node (O1) [draw, thick] at (4,0) {\begin{tabular}{c} $\hat{u}$ \end{tabular}};
	\node (O2) [draw, thick] at (4,-1) {\begin{tabular}{c} $\hat{\gamma}$ \end{tabular}};
	\node (L) [draw, thick] at (8,-0.5) {\begin{tabular}{c} $\mathcal{L}=\mathcal{L}_{\mathcal{N}}+\mathcal{L}_{\mathcal{B}}+\mathcal{L}_{\mathcal{I}}+\mathcal{L}_{\mathcal{M}}$ \end{tabular}};
	\draw [line width=0.4mm,,->] (I1.east) -- (A1.west);
	\draw [line width=0.4mm,,->] (I2.east) -- (A1.west);
	\draw [line width=0.4mm,,->] (I1.east) -- (A2.west);
	\draw [line width=0.4mm,,->] (A1.east) -- (O1.west);
	\draw [line width=0.4mm,,->] (A2.east) -- (O2.west);
	\draw [line width=0.4mm,,->] (O1.east) -- (L.west);
	\draw [line width=0.4mm,,->] (O2.east) -- (L.west);
	
	\draw [line width=0.4mm,,->,gray] (L.north) -- (8, 1) -- (1.5, 1) -- (A1.north);
	\node [gray] at (4.75,1.3) {$\nabla_{\boldsymbol{\theta}_u}\mathcal{L}$};
	\node [gray] at (4.75,0.7) {\footnotesize automatic differentiation};
	\draw [line width=0.4mm,->,gray] (L.south) -- (8, -2) -- (1.5, -2) -- (A2.south);
	\node [gray] at (4.75,-1.7) {$\nabla_{\boldsymbol{\theta}_\gamma}\mathcal{L}$};
	\node [gray] at (4.75,-2.3) {\footnotesize automatic differentiation};
	\end{tikzpicture}
	\caption{Prediction and training of PINNs for inversion with the scalar wave equation. Predictions of the wavefield $\hat{u}$ and indicator $\hat{\gamma}$ are obtained through the neural networks $A_u(\boldsymbol{\theta}_u)$ and $A_\gamma(\boldsymbol{\theta}_\gamma)$. These are updated during training with the gradients of the cost function $\mathcal{L}$ obtained by automatic differentiation.}
	\label{fig:pinnvisualization}
\end{figure}

\textbf{Problem type (c) - Inverse problems with partial knowledge of the solution field} $u(\boldsymbol{x}, t)$: Full knowledge of the wave field, i.e. everywhere in $\Omega \times \mathcal{T}$ is rare in the context of non-destructive testing. Typically receivers can only be placed at parts of the boundary as depicted in~\cref{fig:ndt} such that only partial knowledge is given. Yet, to compute the differential equation loss $\mathcal{L}_{\mathcal{N}}$ of~\cref{eq:losspde}, the wave field $u(\boldsymbol{x}, t)$ and the scaling function $\gamma(\boldsymbol{x})$ must be known everywhere in $\Omega \times \mathcal{T}$. Thus, in the inverse problem with partial domain knowledge, the evaluation of the differential equation loss is computed by approximating both quantities with two separate neural networks, i.e. using both~\cref{eq:solution} and~\cref{eq:coefficient}, as visualized in~\cref{fig:pinnvisualization}. The cost $\mathcal{L}(\boldsymbol{\hat{u}}, \boldsymbol{\hat{\gamma}})$ is then composed of the boundary condition loss,~\cref{eq:lossbc}, initial condition loss,~\cref{eq:lossic}, the differential equation loss $\mathcal{L}_{\mathcal{N}}(\boldsymbol{\hat{u}},\boldsymbol{\hat{\gamma}})$,~\cref{eq:losspde}, and the measurement loss $\mathcal{L}_{\mathcal{M}}(\boldsymbol{\hat{u}},\boldsymbol{u}_\mathcal{M})$,~\cref{eq:measurementloss}. Note, that the measurement data $\boldsymbol{u}_\mathcal{M}$ is now incorporated through the measurement loss, evaluated with the prediction of the forward solution $\hat{\boldsymbol{u}}$. The current state-of-the-art of this inverse problem with PINNs is given e.g.  in~\cite{rashtbehesht_physicsinformed_2022} with a representative Python code provided in~\cite{rasht-behesht_physics-informed_2021}.


\section{Effect of Neural Networks in the Inversion}\label{sec:effect}
To investigate the strengths and weaknesses of PINNs for full waveform inversion, the PINN and its training procedure as illustrated in~\cref{fig:pinnvisualization} is incrementally modified until it is indistinguishable from the classical adjoint optimization. 

The work in~\cite{shukla_physics-informed_2020, rashtbehesht_physicsinformed_2022} discusses the use of PINNs for full waveform inversion, where fully-connected neural networks with inputs $\boldsymbol{x}$ and $t$ were employed. By constrast, in the paper at hand convolutional neural networks are employed, since their spatial equivariance leads to a significant speed-up in training and faster convergence. This has already been reported in a more general context in computational mechanics in e.g.~\cite{zhu_physics-constrained_2019, geneva_modeling_2020} and is not the main focus of the paper at hand. This modification additionally enables a better transition and comparison to the presented methods and is not seen as a detriment to the PINN methodology. Despite this minor difference, we will still consider~\cite{rashtbehesht_physicsinformed_2022} as the reference case for what is possible within the current state-of-the-art using pure PINNs and only partial domain knowledge. 

\subsection{Physics Informed Neural Networks with non-trainable Forward Operator} \label{sec:PINNSnonLearnableFO}
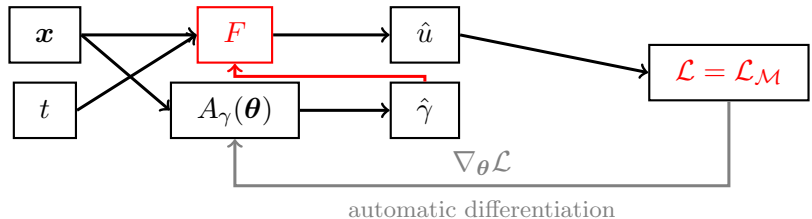
\begin{figure}[htb]
	\centering
	\begin{tikzpicture}
	\node (I1) [draw, thick] at (-1,0) {\begin{tabular}{c} $\boldsymbol{x}$ \end{tabular}};
	\node (I2) [draw, thick] at (-1,-1) {\begin{tabular}{c} $t$ \end{tabular}};
	\node (A1) [draw, thick, red] at (1.5,0) {\begin{tabular}{c} $F$ \end{tabular}};
	
	\node (A2) [draw, thick] at (1.5,-1) {\begin{tabular}{c} $A_{\gamma}(\boldsymbol{\theta})$ \end{tabular}};
	\node (O1) [draw, thick] at (4,0) {\begin{tabular}{c} $\hat{u}$ \end{tabular}};
	\node (O2) [draw, thick] at (4,-1) {\begin{tabular}{c} $\hat{\gamma}$ \end{tabular}};
	\node (L) [draw, thick] at (8,-0.5) {\begin{tabular}{c} $\textcolor{red}{\mathcal{L}=\mathcal{L}_{\mathcal{M}}}$ \end{tabular}};
	\draw [line width=0.4mm,,->] (I1.east) -- (A1.west);
	\draw [line width=0.4mm,,->] (I2.east) -- (A1.west);
	\draw [line width=0.4mm,,->] (I1.east) -- (A2.west);
	\draw [line width=0.4mm,,->] (A1.east) -- (O1.west);
	\draw [line width=0.4mm,,->] (A2.east) -- (O2.west);
	\draw [line width=0.4mm,,->] (O1.east) -- (L.west);
	
	\draw [line width=0.4mm,,->,gray] (L.south) -- (8, -2) -- (1.5, -2) -- (A2.south);
	\node [gray] at (4.75,-1.7) {$\nabla_{\boldsymbol{\theta}}\mathcal{L}$};
	\node [gray] at (4.75,-2.3) {\footnotesize automatic differentiation};
	
	\draw [line width=0.4mm,red,->] (O2.north) -- (4,-0.55) -- (1.5,-0.55) -- (A1.south);
	
	\end{tikzpicture}
	\caption{PINN with non-trainable forward operator. The modification of the original PINN from~\cref{fig:pinnvisualization} is marked in red. Instead of learning the forward problem, a non-trainable forward operator $F$ with $\hat{\gamma}$ as additional input is used. Simultaneously, this leads to a simplification of the cost function $\mathcal{L}$, as $\mathcal{L}_{\mathcal{N}}, \mathcal{L}_{\mathcal{B}}, \mathcal{L}_{\mathcal{I}}$ are approximately satisfied due to $F$.}
	\label{fig:pinnForwardSolver}
\end{figure}

The first modification of PINNs visualized in~\cref{fig:pinnvisualization} is to replace the neural network $A_u(\boldsymbol{\theta})$, whose task is to learn the forward solution, with a non-trainable forward operator $F$ with $\hat{\gamma}$ as additional input. The new structure of a Physics-Informed Neural Network with a non-trainable Forward Operator is illustrated in~\cref{fig:pinnForwardSolver}. The forward operator could be any type of surrogate model such as pre-trained neural networks, or classical solution methods such as the finite element method~\cite{hughes_finite_2000} or the finite difference method~\cite{langtangen_finite_2017}. An important computational requirement on the forward operator $F$ is the viability of a GPU-based implementation, which is an advantage of pre-trained neural networks over CPU-based solvers. However, finite difference methods are commonly also accelerated with GPUs, as e.g. in~\cite{michea_accelerating_2010}. Therefore, finite differences are used in the sequel as the forward operator $F$. A basic but relatively efficient implementation in PyTorch is used by exploiting its convolutional network infrastructure as stencils. This results in simulation runs which take about $500\upmu$s wall clock time for all 1200 timesteps on a $252\times 124$ grid\footnote{tested on an Nvidia A100-SXM4-40GB}. Details are laid out in Appendix~\ref{sec:appendixFD}.

The introduction of a non-trainable forward operator should, in the ideal case, lead to very small values of the loss components related to the initial-boundary value problem, i.e. $\mathcal{L}_{\mathcal{N}},\mathcal{L}_{\mathcal{B}},\mathcal{L}_{\mathcal{I}}$. Typically, these are only approximately zero due to modelling and discretization errors. Their contribution is, however, still negligible in the presence of the measurement loss $\mathcal{L}_{\mathcal{M}}$. For this reason, they can be eliminated from the cost function $\mathcal{L}$, such that it is now solely composed of the measurement loss $\mathcal{L}_{\mathcal{M}}$.

The introduction of a non-trainable forward operator leads to a significantly less complex optimization task. In the case of a time-stepping scheme as forward operator $F$, the disadvantage of the design depicted in~\cref{fig:pinnForwardSolver} is a more expensive sensitivity computation $\nabla_{\boldsymbol{\theta}}\mathcal{L}(F(\boldsymbol{\hat{\gamma}}))~=~\nabla_{\boldsymbol{\hat{\gamma}}}\mathcal{L}(F(\boldsymbol{\hat{\gamma}}))~\cdot~ \nabla_{\boldsymbol{\theta}}\boldsymbol{\hat{\gamma}}(\boldsymbol{\theta})$ because the backpropagation needs to be piped through each timestep of the forward operator $F$ by means of automatic differentiation to compute $\nabla_{\boldsymbol{\hat{\gamma}}}\mathcal{L}(F(\boldsymbol{\hat{\gamma}}))$. The cost of computing the gradient of the inverse field w.r.t. the network parameters $\nabla_{\boldsymbol{\theta}}\boldsymbol{\hat{\gamma}}(\boldsymbol{\theta})$ is unchanged.

\subsection{Adjoint Method} \label{sec:AdjointMethod}
\begin{figure}[htb]
	\centering
	\begin{tikzpicture}
	\node (I1) [draw, thick] at (-1,0) {\begin{tabular}{c} $\boldsymbol{x}$ \end{tabular}};
	\node (I2) [draw, thick] at (-1,-1) {\begin{tabular}{c} $t$ \end{tabular}};
	\node (A1) [draw, thick] at (1.5,0) {\begin{tabular}{c} $F$ \end{tabular}};
	
	\node (A2) [draw, thick, red] at (1.5,-1) {\begin{tabular}{c} $C_{\gamma}(\boldsymbol{c})$ \end{tabular}};
	\node (O1) [draw, thick] at (4,0) {\begin{tabular}{c} $\hat{u}$ \end{tabular}};
	\node (O2) [draw, thick] at (4,-1) {\begin{tabular}{c} $\hat{\gamma}$ \end{tabular}};
	\node (L) [draw, thick] at (8,-0.5) {\begin{tabular}{c} $\mathcal{L}=\mathcal{L}_{\mathcal{M}}$ \end{tabular}};
	\draw [line width=0.4mm,,->] (I1.east) -- (A1.west);
	\draw [line width=0.4mm,,->] (I2.east) -- (A1.west);
	\draw [line width=0.4mm,,->] (I1.east) -- (A2.west);
	\draw [line width=0.4mm,,->] (A1.east) -- (O1.west);
	\draw [line width=0.4mm,,->] (A2.east) -- (O2.west);
	\draw [line width=0.4mm,,->] (O1.east) -- (L.west);
	
	\draw [line width=0.4mm,,->] (O2.north) -- (4,-0.55) -- (1.5,-0.55) -- (A1.south);
	
	\draw [line width=0.4mm,,->,gray] (L.south) -- (8, -2) -- (1.5, -2) -- (A2.south);
	\node [gray] at (4.75,-1.7) {$\textcolor{red}{\nabla_{\boldsymbol{c}}}\mathcal{L}$};
	\node [gray] at (4.75,-2.45) {\footnotesize \begin{tabular}{c}
		automatic differentiation\\
		or adjoint method
		\end{tabular}};
	\end{tikzpicture}
	\caption{Adjoint method. The modification of the previously modified PINN from~\cref{fig:pinnForwardSolver} is indicated in red. The neural network Ansatz is replaced by a constant Ansatz parametrized by the coefficients $\boldsymbol{c}$.}
	\label{fig:PINNconstantAnsatz}
\end{figure}
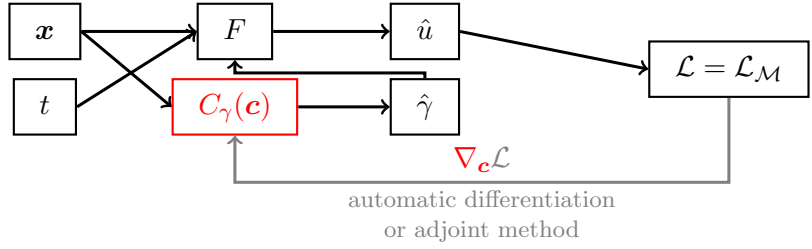

The next step is to replace the neural network Ansatz of the indicator function $\gamma$ with a piece-wise polynomial Ansatz. Here, a piece-wise constant Ansatz is selected and defined in the one-dimensional case as $\hat{\gamma}(x)=\sum_{i=1}^{N_c} c_i N_i(x)$ with
\begin{equation}
N_i(x)=\begin{cases}
1 \text{ if } x_i\leq x<x_{i+1} \\
0 \text{ else. }
\end{cases} \label{eq:PWConstAnsatz}
\end{equation}
For the two-dimensional and three-dimensional cases the discretization is handled in the same way for the additional dimensions. This modification is depicted in~\cref{fig:PINNconstantAnsatz}. The coefficients $\boldsymbol{c}=\{c_i\}_{i=1}^{N_c}$ are still updated with the Adam optimizer using the gradients obtained through automatic differentiation. The reason for choosing a piecewise constant Ansatz over a piecewise higher-order polynomial is that a piecewise constant Ansatz is already capable of representing the jump of going from a constant material to a void, as is the case in all examples presented in~\cref{sec:configuration}. 
Without the neural network, the Sigmoid activation function is missing and the range of $\hat{\gamma}$ has to be constrained in a different manner. To this end, an additional weight, i.e. coefficient clipping is used, where the absolute value is clipped to the range of a small value $\epsilon$ close to zero and an upper bound one or slightly exceeding one \cite{buerchner_2022}.

Given~\cref{eq:PWConstAnsatz}, the sensitivity w.r.t. the coefficients
\begin{align}
\nabla_{\boldsymbol{c}}\mathcal{L}(F(\hat{\boldsymbol{\gamma}}))~=~\nabla_{\hat{\boldsymbol{\gamma}}} \mathcal{L}(F(\hat{\boldsymbol{\gamma}}))\cdot \nabla_{\boldsymbol{c}} \hat{\boldsymbol{\gamma}}(\boldsymbol{c})
\label{eq:chainRuleForC}
\end{align}
can either be computed via automatic differentiation or the continuous adjoint method, as indicated in~\cref{fig:PINNconstantAnsatz}.  In the continuous adjoint method, the gradient is derived analytically before its discretization, as e.g. derived in~\cite{buerchner_2022}. This third modification has transformed the approach into the classical full waveform inversion based on finite differences. No major implications on the accuracy are expected due to this modification, as the gradients are the same, apart from numerical errors introduced by the finite difference discretization in the adjoint method. These originate in the computation of the Fréchet kernel
\begin{equation}
K_{\gamma}=\int_{\mathcal{T}} -\rho_0 u_t^{\dagger} u_t+\rho_0 c_0^2 \nabla u^{\dagger} \cdot \nabla u dt \label{eq:frechet}
\end{equation}
used to compute the gradient $\nabla_{\hat{\gamma}} \mathcal{L}_{\mathcal{M}}$ via the integration of $K_{\gamma}$ over the spatial and time domains (c.f.~\cite{buerchner_2022} for more details). To compute~\cref{eq:frechet}, the temporal and spatial gradients of the solution field $u$ and the adjoint field $u^{\dagger}$ are computed with finite differences. Note, that the sensitivity w.r.t. the coefficients $c_i$ is computed in a second step, by multiplication with the corresponding shape function $N_i(x)$, obtained by deriving~\cref{eq:PWConstAnsatz} w.r.t. $c_i$, i.e. $\nabla_{\boldsymbol{c}} \hat{\boldsymbol{\gamma}}(\boldsymbol{c})$. 

In fact, the key difference between the adjoint method and automatic differentiation is the order in which the gradient computation and discretization of the state equation is performed. In the continuous adjoint method, the gradient is first derived analytically, and thereupon the state equation is discretized as described above. By contrast, in the case of automatic differentiation, the differential equation is discretized first, and the gradient is propagated through each timestep. Thus, automatic differentiation bears a great similarity to the discrete adjoint method and is, as stated in~\cite{norgaard_applications_2017}, fundamentally the same approach. Manually deriving the sensitivities in the discrete adjoint or computing them through reverse mode automatic differentiation as in~\cite{norgaard_applications_2017, dilgen_topology_2018} will not result in a difference in accuracy. Thus, regarding the accuracy, the difference between automatic differentiation and the continuous adjoint method is only affected by the difference between the discrete and continuous adjoint method.

Unlike the accuracy, the computational cost of the automatic differentiation approach compared to the continuous adjoint method is different. Automatic differentiation adds a significant memory cost due to the backpropagation. In addition, the backpropagation along all timesteps is costly with automatic differentiation, as also mentioned in the work on the ordinary differential equation networks~\cite{chen_neural_2019}. The theoretical complexity of the backpropagation and continuous adjoint method w.r.t. the number of timesteps is in both cases $O(n)$~\cite{williams_gradient-based_1995} and in practice only a factor of roughly two is observed between the computation times.  
To save memory, a truncated backpropagation through time~\cite{williams_gradient-based_1995, sutskever_training_2013} can be used, as typically carried out for recurrent neural networks. However, a truncation was not used in the examples investigated in~\cref{sec:configuration} since it was found to be detrimental to the learning behavior of the indicator function. 

\subsection{Hybrid Approach} \label{sec:hybridApproach}
\begin{figure}[htb]
	\centering
	\begin{tikzpicture}
	\node (I1) [draw, thick] at (-1,0) {\begin{tabular}{c} $\boldsymbol{x}$ \end{tabular}};
	\node (I2) [draw, thick] at (-1,-1) {\begin{tabular}{c} $t$ \end{tabular}};
	\node (A1) [draw, thick] at (1.5,0) {\begin{tabular}{c} $F$ \end{tabular}};
	
	\node (A2) [draw, thick] at (1.5,-1) {\begin{tabular}{c} $A_{\gamma}(\boldsymbol{\theta_\gamma})$ \end{tabular}};
	\node (O1) [draw, thick] at (4,0) {\begin{tabular}{c} $\hat{u}$ \end{tabular}};
	\node (O2) [draw, thick] at (4,-1) {\begin{tabular}{c} $\hat{\gamma}$ \end{tabular}};
	\node (L) [draw, thick] at (8,-0.5) {\begin{tabular}{c} $\mathcal{L}=\mathcal{L}_{\mathcal{M}}$ \end{tabular}};
	\draw [line width=0.4mm,,->] (I1.east) -- (A1.west);
	\draw [line width=0.4mm,,->] (I2.east) -- (A1.west);
	\draw [line width=0.4mm,,->] (I2.east) -- (A2.west);
	\draw [line width=0.4mm,,->] (A1.east) -- (O1.west);
	\draw [line width=0.4mm,,->] (A2.east) -- (O2.west);
	\draw [line width=0.4mm,,->] (O1.east) -- (L.west);
	
	\draw [line width=0.4mm,,->] (O2.north) -- (4,-0.55) -- (1.5,-0.55) -- (A1.south);
	
	\draw [line width=0.4mm,red] (L.south) -- (8, -2.2) -- (4, -2.2) -- (O2.south);
	
	\draw [line width=0.6mm,,->,red] (O2.south) -- (4, -2.2) -- (1.5, -2.2) -- (A2.south);
	
	\node [red] at (8.7,-1.8) {$(\nabla_{\boldsymbol{\hat{\gamma}}}\mathcal{L})^T$};
	
	\node [red] at (6.7,-1.8) {$\boldsymbol{\cdot}$};
	
	\node [red] at (4.7,-1.8) {$(\nabla_{\boldsymbol{\theta}}\boldsymbol{\hat{\gamma}})^T$};
	
	\node [red] at (2.75,-1.8) {$(\nabla_{\boldsymbol{\theta}}\mathcal{L})^T=$};

	\node [red] at (6,-2.5) {\footnotesize adjoint method};
	\node [red] at (2.75,-2.5) {\footnotesize automatic differentiation};

	\end{tikzpicture}
	\caption{Hybrid approach. Modification of the previously modified PINN from~\cref{fig:pinnForwardSolver} is marked in red. The gradient computation is now split into two parts. First, the sensitivity of the cost function $\mathcal{L}$ w.r.t. the predicted indicator field $\hat{\gamma}$ is computed with the adjoint method using~\cref{eq:frechet}, whereupon the sensitivity of the indicator field $\hat{\gamma}$ w.r.t. the neural network parameters $\boldsymbol{\theta}$ is found with the automatic differentiation  used to train the neural network. This yields the desired gradient $\nabla_{\boldsymbol{\theta}}\mathcal{L}$ by using the chain rule~\cref{eq:chainrule}.}
	\label{fig:HybridApproach}
\end{figure}
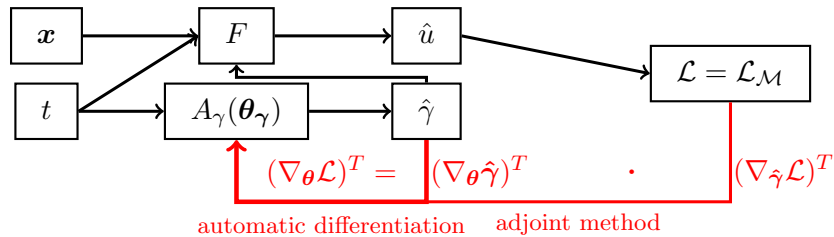

This section introduces a hybrid approach which combines the advantages of a neural network Ansatz for the scaling function $\gamma$ and the more efficient gradient computation of the gradient of the loss functional~\cref{eq:measurementloss} using the continuous adjoint method~\cref{eq:frechet}. To this end, the neural network is reintroduced as an Ansatz for $\gamma$. In analogy to~\cref{eq:chainRuleForC}, the gradient of the cost function can be computed with the chain rule but must now be computed with respect to the neural network parameters $\boldsymbol{\theta}$:
\begin{equation}
\nabla_{\boldsymbol{\theta}}\mathcal{L}(\boldsymbol{\hat{u}}(\boldsymbol{\hat{\gamma}}(\boldsymbol{\theta})))= \nabla_{\boldsymbol{\hat{\gamma}}}\mathcal{L}(\boldsymbol{\hat{u}}(\boldsymbol{\hat{\gamma}})) \cdot \nabla_{\boldsymbol{\theta}}\boldsymbol{\hat{\gamma}}(\boldsymbol{\theta}) \label{eq:chainrule}
\end{equation}
The sensitivity of the cost function w.r.t. the neural network outputs, i.e. the indicator field, $\nabla_{\boldsymbol{\hat{\gamma}}}\mathcal{L}$ is efficiently computed with the continuous adjoint method~\cref{eq:frechet}, while the sensitivity of the indicator field w.r.t. the neural network parameters $\nabla_{\boldsymbol{\theta}}\boldsymbol{\hat{\gamma}}$ is obtained via backpropagation through the neural network $A_{\gamma}(\boldsymbol{\theta})$. The hybrid approach results in similar computational times as the adjoint method for a sufficient number of timesteps. The scheme is illustrated in~\cref{fig:HybridApproach}. A similar scheme was recently proposed in the context of topology optimization~\cite{chandrasekhar_tounn_2021}.\\
In addition to saving memory compared to the scheme discussed in the previous section, the hybrid scheme bears the advantage of separating neural networks from forward solvers such that an efficient GPU implementation of the forward solver is no longer strictly necessary. Instead, the neural network can provide the indicator prediction $\boldsymbol{\hat{\gamma}}$ and its sensitivity w.r.t. the parameters $\nabla_{\boldsymbol{\theta}}\boldsymbol{\hat{\gamma}}$ on a CPU while, the forward solver can also be executed independently on a CPU. The increase in computational effort of evaluating and updating the neural network on a CPU instead of a GPU is very small compared to the cost of the forward solver. Furthermore, the forward solver can even be a completely independent black-box solver which provides the wave-field $\hat{u}$ and the sensitivity of the loss functional $\nabla_{\boldsymbol{\hat{\gamma}}}\mathcal{L}$. 

\section{Configuration} \label{sec:configuration}
Two model problems with increasing complexity are considered. First, a two-dimensional plate is used to determine the strengths and weaknesses of neural networks in the inversion process. Afterwards, a more complex three-dimensional case is used to compare the proposed hybrid approach to the adjoint method. The forward operator is always a finite difference solver. To avoid an inverse crime~\cite{wirgin_inverse_2004}, the reference data was generated with a boundary-conforming finite element solver for the partial domain knowledge cases (Voxel-FEM for the three-dimensional case). However, in the case of the two-dimensional full domain knowledge case the finite cell method \cite{parvizian_finite_2007} is employed for the reference data generation. In all studies, the sine burst, presented in~\cref{eq:sineburst} is applied with $n_c=2$ cycles, a frequency of $f_{\psi} = 500$ kHz, and an amplitude of $\psi_0 = 10^{12} \text{ kg}(\text{ms})^{-2}$. Special care must be taken when using finite differences during the source application in \cref{eq:source}. The Dirac delta has to be scaled by a factor of $1/ \vert \vert d\boldsymbol{x}\vert \vert  _2^2$ to account for the grid spacings~\cite{hicks_arbitrary_2002}.
Additionally, both cases share the same material with a density of $\rho_0=2700\text{ kg/m}^3$ and a wave speed of $6000 \text{ ms}^{-1}$.

\subsection{Two-Dimensional Case} \label{sec:2dcase}
\begin{figure}[htb]
	\centering
	\begin{tikzpicture}
	\draw [line width=0.5mm] (-9.42857,0) -- (-0.7,0);
	
	\fill [red] (-6,0) circle (0.15cm);
	
	\fill [red] (-2.4,0) circle (0.15cm);
	
	\foreach {\x} in {-6,-5.8,...,-0.8} {
		\fill [black] (\x,0) circle (0.1cm);
	}
	
	\draw [line width=0.5mm] (-9.42857,0) -- (-9.42857,-0.5);
	\draw [dashed, line width=0.5mm] (-9.42857,-0.5) -- (-9.42857,-1);
	
	\draw [dash pattern={on 7pt off 2pt on 1pt off 3pt}, gray, line width=0.3mm] (-0.7,1) -- (-0.7,-1);
	
	\node [gray] at (-7.75,0.7) {\footnotesize 0.0398mm};
	\draw [gray, |-|, line width=0.3mm] (-9.42857, 0.5) -- (-2, 0.5);
	
	\node [gray] at (-7.75, 1.2) {\footnotesize 0.0183mm};
	\draw [gray, |-|, line width=0.3mm] (-9.42857, 1) -- (-6, 1);
	
	\node [gray] at (-5.9, -0.55) {\footnotesize 0.00120mm};
	\draw [gray, |-|, line width=0.3mm] (-5.8, -0.3) -- (-6, -0.3);
	
	
	\node [gray] at (-8.5, -2.35) {$x$};
	\node [gray] at (-9.7,-1.3) {$y$};
	\draw [gray, <->, line width=0.3mm] (-9.42857, -1.1) -- (-9.42857, -2.1) -- (-8.42857, -2.1);
	
	\node at (-8.7,0.25) {$\Gamma_{x_{1}}$};
	\node at (-9.75,-0.4) {$\Gamma_{y_{0}}$};
	
	\end{tikzpicture}
	\caption{Symmetrical sensor (black circles) and source (red circles) configuration.}\label{fig:sensorconfig}
\end{figure}
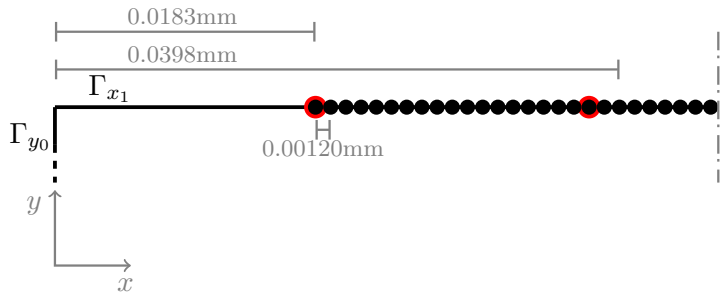

\begin{figure} 
	\centering
	\begin{subfigure}[b]{0.32\textwidth}
		\includegraphics[width=\textwidth]{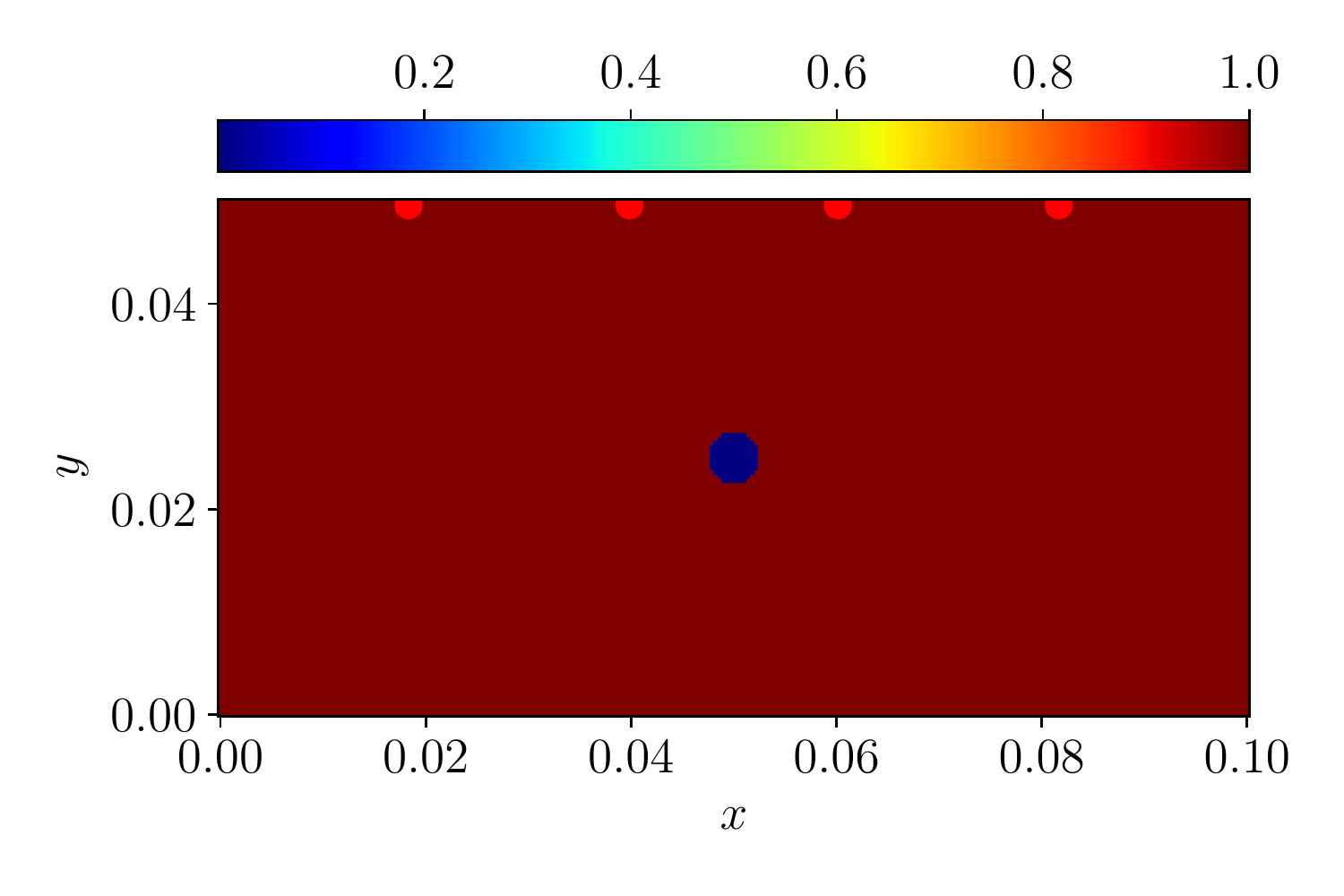}
		\caption{Material defined by the indicator function $\gamma(x,y)$} \label{fig:forwardsimulationa}
	\end{subfigure}
	\hfill
	\begin{subfigure}[b]{0.32\textwidth}
		\includegraphics[width=\textwidth]{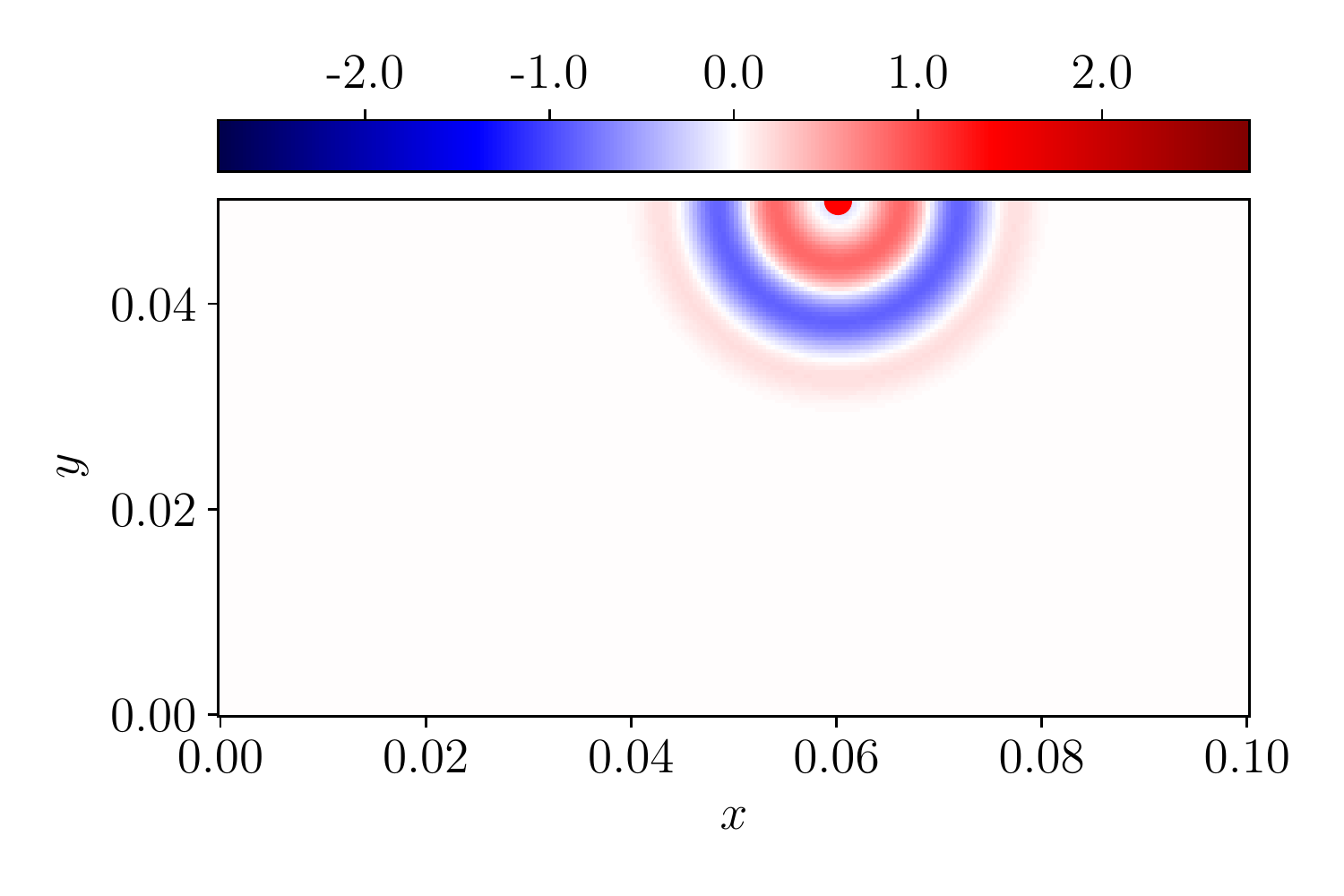}
		\caption{Wavefield $u(x,y,t)$ at $t=0.00375\upmu$s.} \label{fig:forwardsimulationb}
	\end{subfigure}
	\hfill
	\begin{subfigure}[b]{0.32\textwidth}
		\includegraphics[width=\textwidth]{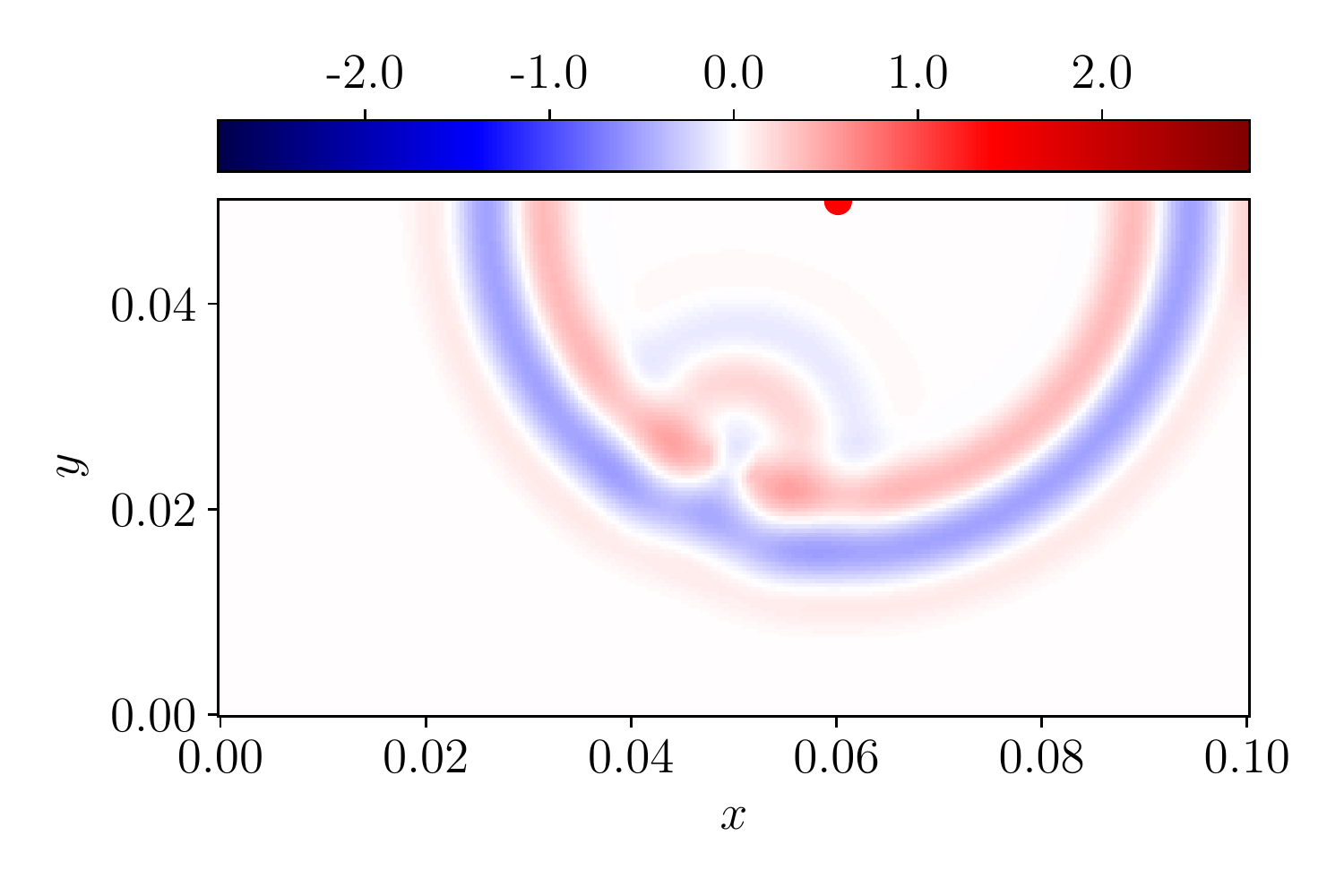}
		\caption{Wavefield $u(x,y,t)$ at $t=0.0075\upmu$s.} \label{fig:forwardsimulationc}
	\end{subfigure}
	\caption{The wavefield is computed with the underlying material distribution, depicted in the leftmost figure, with the third source. When the waves reach the defect, a disturbance in the wavefield can be seen.} \label{fig:forwardsimulation}
\end{figure}

The two-dimensional case is strongly inspired by~\cite{buerchner_2022}, where a circular hole with a diameter of $5$ mm is to be detected on a $100$ mm $\times$ $50$ mm rectangular plate, illustrated in~\cref{fig:forwardsimulationa}. We do not provide specific information about the hole to the machine learning algorithms, which implies the entire material distribution $\gamma(\boldsymbol{x})$ has to be estimated. While the geometry, material, and source type are identical, the sensor and source locations are slightly different from~\cite{buerchner_2022}. In total, 54 sensors and four sources are used on the top edge in four separate measurements. The exact configuration is provided in~\cref{fig:sensorconfig}, showing the upper left edge of the domain using the domain definitions from~\cref{fig:ndt}. The setup is symmetric along the vertical midline indicated by the dashed grey line. The slight peculiarity of the sensor and source locations is due to them being aligned to the finite difference nodes and neural network grid points. This is, however, not a limitation of the presented methods, as it can easily be circumvented by, e.g. using Kaiser windowed sinc functions as interpolation of the sources~\cite{hicks_arbitrary_2002} and linear interpolation for the sensors.

A forward simulation is illustrated in~\cref{fig:forwardsimulation}, where the third source emits a pulse at $t=0$ s, which then leads to a disturbance of the wavefield in the center due to the circular void. The finite difference simulation for the inversion uses $252\times 124$ grid points and two ghost cells per dimension, a lower indicator function bound of $\epsilon=10^{-5}$, and a timestep size of $\Delta t=3 \cdot 10^{-8}\text{ s}^{-1}$ for $1200$ timesteps. 

\subsection{Three-Dimensional Case}
\begin{figure}
	\centering
	\begin{subfigure}[b]{0.8\textwidth}
		\vspace{-0.5cm}
		\begin{tikzpicture}
		
		\node at (0,-3.9) (red){
			\includegraphics[width=\textwidth, trim={50cm 0 0 35cm},clip]{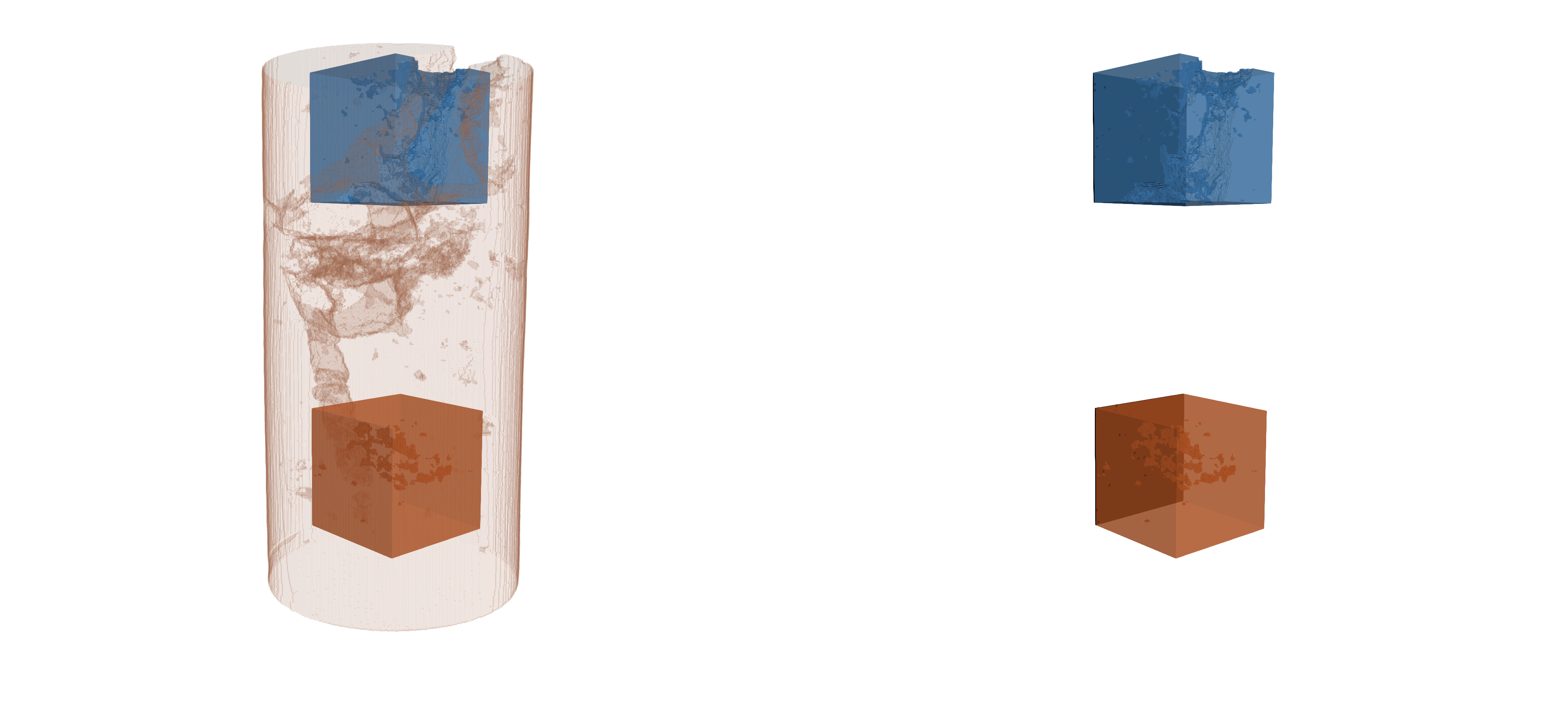}
		};
		
		\node at (0,0.4) (blue){
			\includegraphics[width=\textwidth, trim={50cm 45cm 0 0},clip]{code/Ctscan.png}
		};
		
		\node at (-4,-1.7) (cylinder){
			\includegraphics[width=0.5\textwidth, trim={0 0 80cm 0},clip]{code/Ctscan.png}
		};				
		
		\node at (1.6,1.3) {case 1};
		\node at (1.6,-1.7) {case 2};
		\draw [line width=0.5mm,->] (-2.5,0.4) -- (0,0);
		\draw [line width=0.5mm,->] (-2.5,-2.9) -- (0,-2.8);


		\end{tikzpicture}
	\end{subfigure}
	\vspace{-1cm}
	\caption{The two cases extracted from a CT-scanned drill core~\cite{hug_three-field_2022}. The sources and sensors are applied on the backside of the cubes indicated by the darker surface on the right.} \label{fig:3dcase}
\end{figure}

The purpose of the three-dimensional case is to showcase the possible strengths of neural networks when the adjoint optimization is used in a setting where it does not perform well anymore. For this, two cuboid slices from a CT-scanned drill core~\cite{hug_three-field_2022} are considered and illustrated in~\cref{fig:3dcase}. The edge length is $20$ mm and discretized by $92$ grid points in each dimension during the inversion, where 1200 timesteps with a timestep size of $\Delta t=1.8 \cdot 10^{-8}\text{ s}^{-1}$ are considered. Again a lower indicator function bound of $\epsilon=10^{-5}$ is used, while the upper bound is increased to $1.2$ for the adjoint method, as this was found to be beneficial. The sources and sensors are applied to the backside of the cubes, as indicated by the darker shading on the right of~\cref{fig:3dcase}. The four sources are placed at a distance of $5.27$ mm from each edge, while the sensors coincide with the grid points of the plane. \\
The case was selected such that a single forward and adjoint simulation can be performed within the limits of the $40$ GB GPU memory. To compute larger three-dimensional, efficient parallelizations of the simulation have to be used, which can either be implemented on GPU clusters or alternatively on CPUs using efficient communication between the GPU and CPU for the neural network.

\subsection{Neural Network Architecture}
\begin{figure} 
	\centering
	\includegraphics[width=\textwidth]{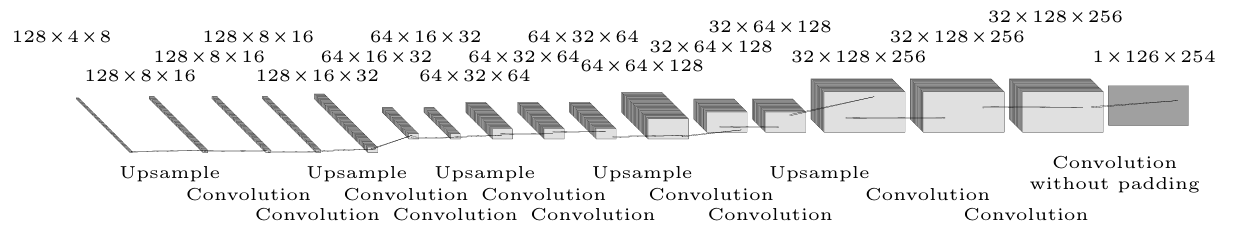}
	\caption{The convolutional neural network for the two-dimensional inversion case. It performs upsampling and convolutions in a repeating pattern on a random input of size $128\times 4\times 8$, yielding a prediction $\hat{\gamma}$ of size $1\times 126 \times 254$. For more details, see Appendix~\ref{sec:appendixNN}.} \label{fig:neuralnetwork}
\end{figure}

The neural network considered as Ansatz for the indicator function $\hat{\gamma}$ has a generator structure with pixel-wise normalization inspired by generative adversarial networks~\cite{karras_progressive_2018}. It upsamples a random three-dimensional matrix of size $128 \times 8\times 4$ with subsequent convolutional filters. The upsampling and convolutions are repeated five times to get the desired resolution for the indicator function $\hat{\gamma}$, see~\cref{fig:neuralnetwork} for the neural network used for the two-dimensional case. A similar network is employed in the three-dimensional case using three-dimensional convolutions and a four-dimensional random input matrix of size $256 \times 3 \times 3 \times 3$. To increase the ratio between the number of indicator voxels and network parameters, twice as many filters are employed in the three-dimensional case. For more details on both networks, c.f. Appendix~\ref{sec:appendixNN}. Nearest-Neighbor interpolation was used for upsampling, as it performed better than polynomial interpolation, as it preserves discontinuities needed to represent the voids. Due to their convergence accelerating behavior~\cite{jagtap_adaptive_2020}, adaptive activation functions are used. Specifically, PReLU~\cite{he_delving_2015} was used for all layers except the last layer, where an adaptive Sigmoid was used to ensure an output in the range $[0,1]$. The output was then rescaled to the range $[\epsilon, 1]$ to ensure numerical stability. 
Optimization is always conducted with Adam~\cite{kingma_adam_2017} to make the comparison more straightforward, although L-BFGS leads to fewer iterations. Gradient clipping~\cite{Goodfellow-et-al-2016, pascanu_difficulty_2013, zhang_why_2020} was found to be essential to ensure an early convergence. Without it, a plateau was observed, followed by a strongly delayed breakthrough. In addition, a learning rate scheduler with a polynomial decay $(b \cdot \text{epoch} + 1)^a$ was employed with the values $a=-0.5$, $b=0.2$. The initial network weights are set with Glorot initialization~\cite{kumar_weight_2017}.


%
%
%
%
%
%
%
%

\section{Numerical Results}\label{sec:results}

\subsection{Two-Dimensional Case}
In the subsequent sections, the individual components of PINNs will be investigated and compared to the classical adjoint optimization using the case presented in~\cref{sec:2dcase}. All simulations are performed on an Nvidia A100 SXM4 40GB GPU and the same convolutional neural network presented in~\cref{fig:neuralnetwork} and Appendix~\ref{sec:appendixNN}.

\subsubsection{Physics-Informed Neural Networks}
\begin{figure}
	\centering
	\begin{subfigure}[b]{0.49\textwidth}
		\includegraphics[width=\textwidth]{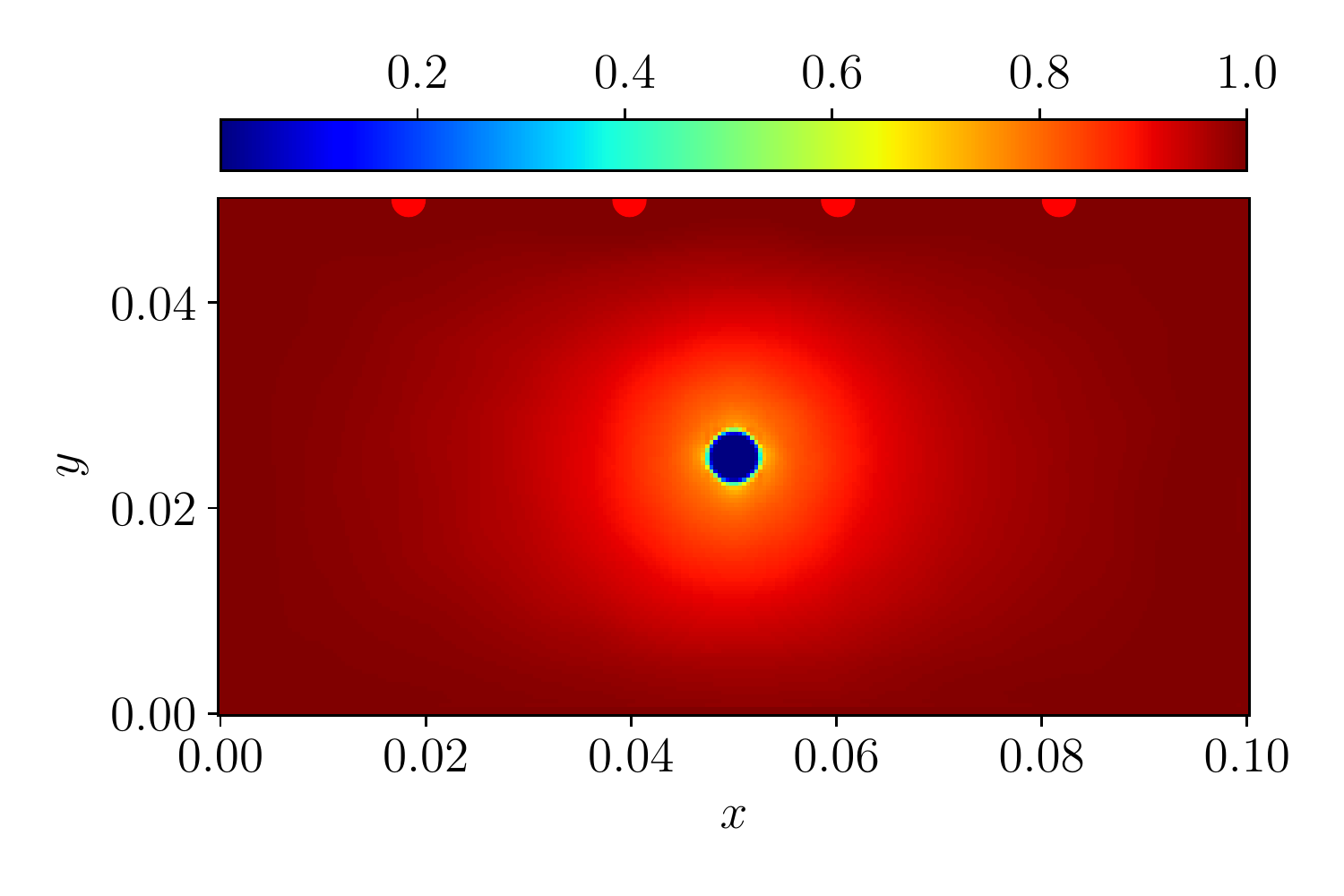}
		\caption{Prediction.}
	\end{subfigure}
	\hfill
	\begin{subfigure}[b]{0.49\textwidth}
		\includegraphics[width=\textwidth]{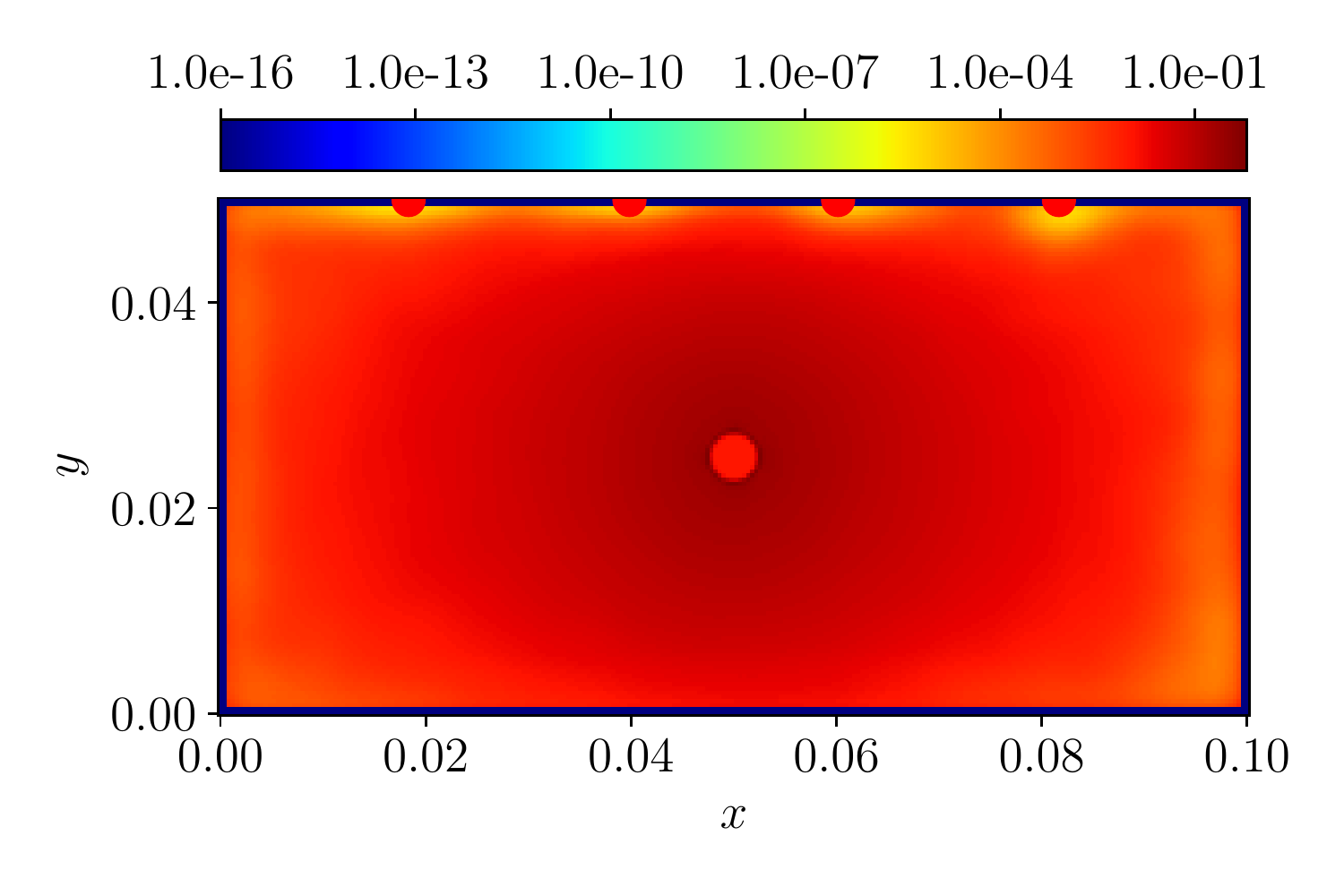}
		\caption{Absolute error.}
	\end{subfigure}
	\caption{Prediction and error of the material distribution after 100 epochs using a PINN with full domain knowledge. Average elapsed time per epoch is $0.51$ s.} \label{fig:pinninversion}
\end{figure}

The results presented in this section refer to case PINNs,~\cref{tab:methods}, whose methodology is described in~\cref{sec:pinns}. The various inversion tasks presented in \cite{rashtbehesht_physicsinformed_2022} ran for about $400\,000-600\,000$ epochs with computational times of about $0.2$ s per epoch\footnote{Based on the Python script provided by~\cite{rasht-behesht_physics-informed_2021} and tested on an Nvidia A100-SXM4-40GB GPU.}. The need for such a large number of epochs illustrates the complexity of the nested optimization task of simultaneously learning the forward and inverse solution for a single problem. Additional complications arise when multiple sources are applied, and a forward operator taking the different sources into account must be learned. In~\cite{rashtbehesht_physicsinformed_2022}, this problem was circumvented by using source stacking~\cite{fichtner_source_2011}, such that only one forward simulation had to be learned. \\
If only one field has to be learned, as in the case of inversion with full domain knowledge, as also presented in~\cite{shukla_physics-informed_2020}, the complexity of the optimization task decreases. To illustrate the decrease in complexity and the capabilities of PINNs in specialized settings, the problem presented in~\cref{sec:2dcase} is solved with a PINN using additional sensors. The sensors are placed throughout the entire domain, coinciding with the finite difference grid points. As the measurements are only given at discrete points, the residual of the differential equation~\eqref{eq:losspde} is evaluated numerically using finite differences for the derivative terms. Additionally, the indicator function is prescribed at the boundaries, i.e. set to one to ensure material at the boundaries, which was essential to retrieve an accurate inversion. The resulting inversion after 100 epochs is provided in~\cref{fig:pinninversion} with a clear identification of the circular void and only minor artifacts. It was essential to prescribe the density $\rho(\boldsymbol{x})$ at the boundary and use a relatively large learning rate $\alpha_{\min}=2\cdot 10^{-3}$. Furthermore the use of the learnable penalty weights using the attention mechanism was helpful with a learning rate of $\alpha_{\max}=2\cdot 10^{-2}$. \\
It is important to emphasize that the improvement from $400\,000-600\,000$ epochs to 100 epochs is a result of the simplification of the physical problem by using full domain knowledge, which should not be taken as a reference of what is possible with PINNs in the partial domain knowledge case presented in~\cref{sec:2dcase}. Instead, $400\,000-600\,000$ epochs are to be considered as the reference of the capabilities of PINNs for a similar inversion task.




\subsubsection{Physics-Informed Neural Networks with non-trainable Forward Operator} 
\begin{figure}
	\centering
	\begin{subfigure}[b]{0.49\textwidth}
		\includegraphics[width=\textwidth]{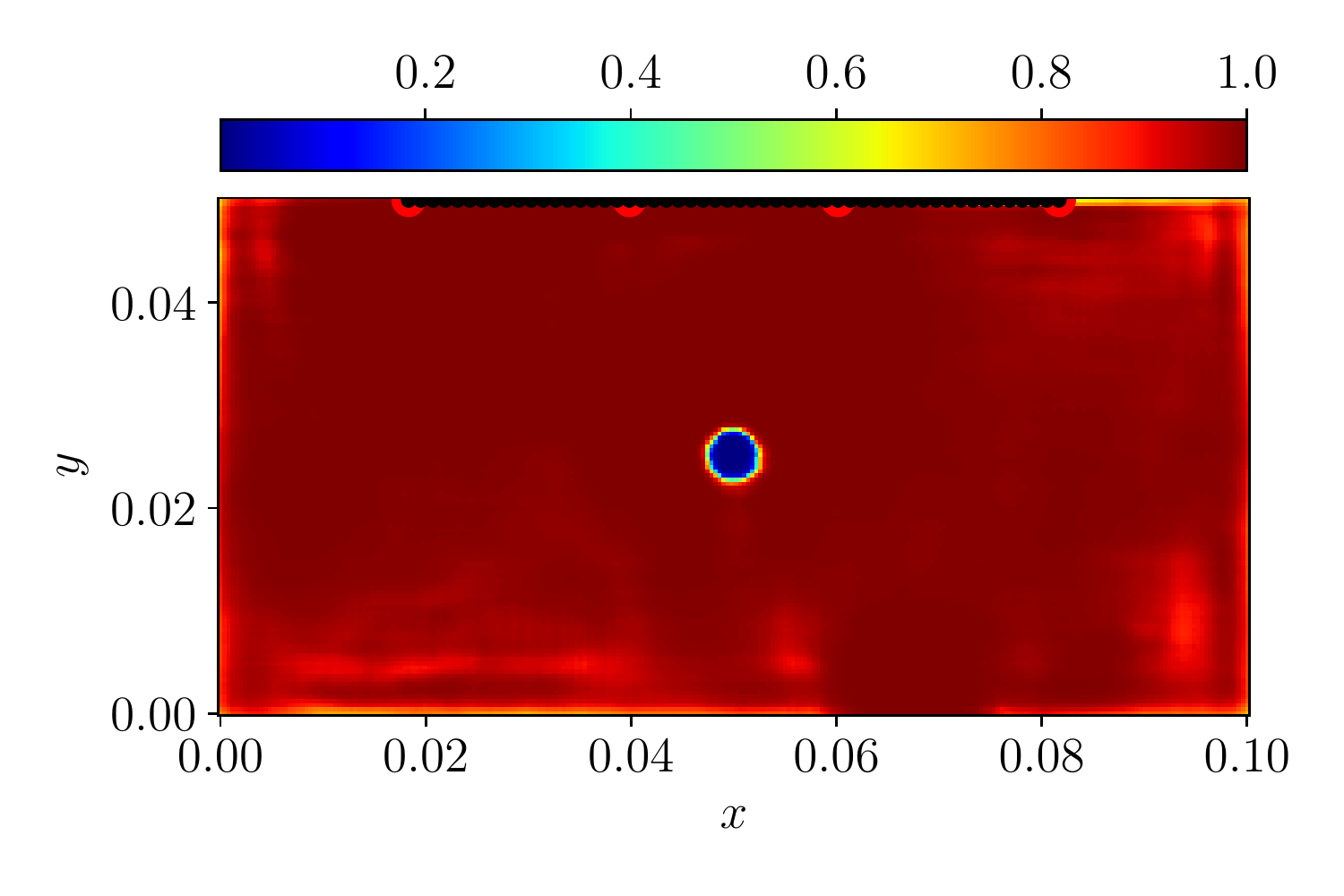}
		\caption{Prediction.}
	\end{subfigure}
	\hfill
	\begin{subfigure}[b]{0.49\textwidth}
		\includegraphics[width=\textwidth]{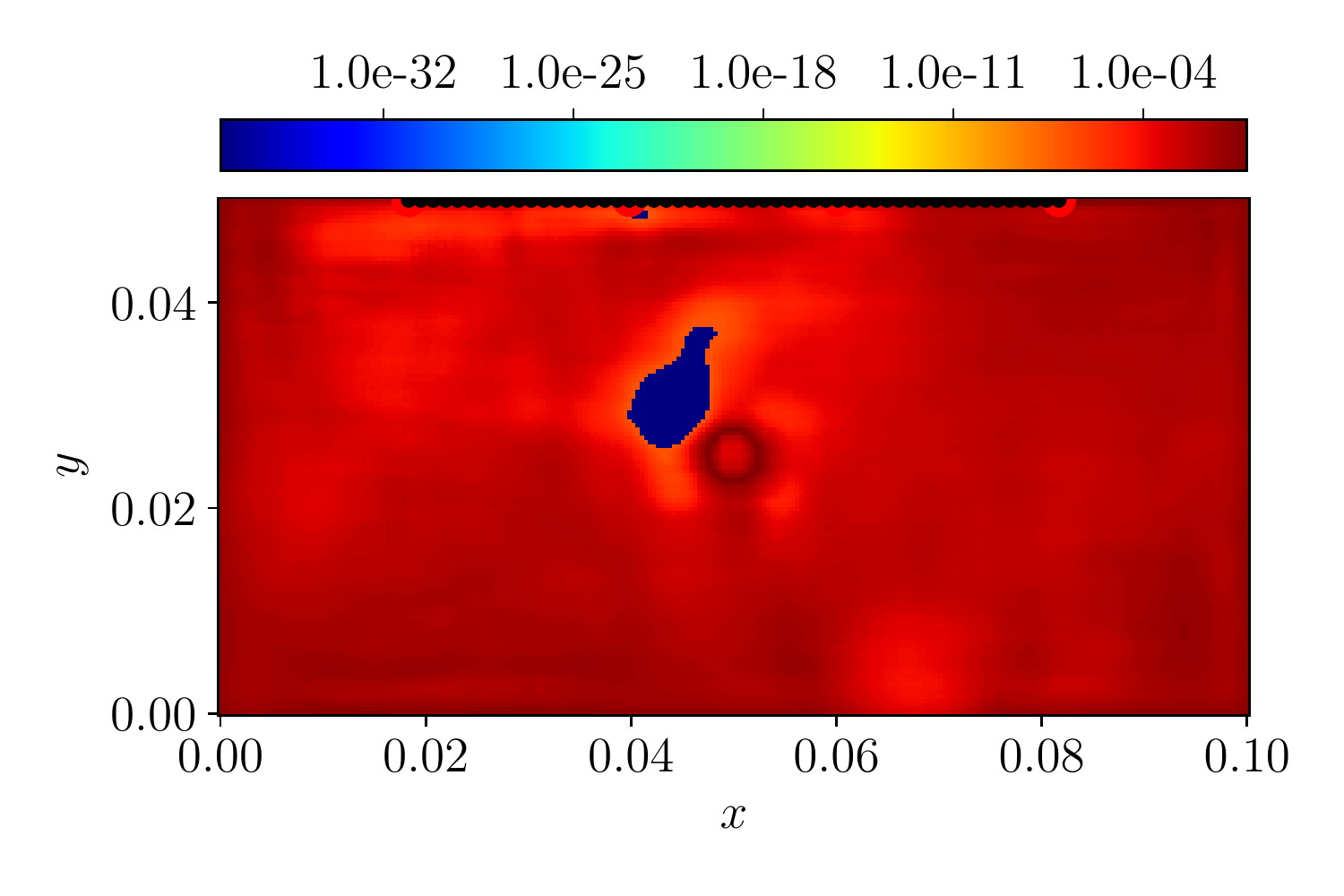}
		\caption{Absolute error.}
	\end{subfigure}
	\caption{Prediction and error of the material distribution after 50 epochs using automatic differentiation and a neural network Ansatz. Average elapsed time per epoch is $2.53$ s.} \label{fig:inversionForwardSolver}
\end{figure}

An alternative to reducing the complexity of learning both the forward and inverse solutions is using a non-trainable forward operator classified in~\cref{tab:methods} as PINNs with non-trainable forward operator and described in~\cref{sec:PINNSnonLearnableFO}. The neural network responsible for the forward solution is now replaced by a finite difference scheme. A reasonable inversion, shown in~\cref{fig:inversionForwardSolver}, is reached after only 50 epochs. Thus the convergence w.r.t. iterations is even faster than the PINN with full domain knowledge, with the detriment of a greater elapsed time per epoch at $2.53$ s. Overall, considering the $400\,000-600\,000$ epochs observed in~\cite{rashtbehesht_physicsinformed_2022} given only partial domain knowledge, the PINN is outperformed by using a non-trainable forward operator. For the optimization a learning rate of $\alpha_{\min}=2\cdot10^{-3}$ was employed with gradient clipping. The learnable penalty weights improved the results slightly, but were not found to be essential.

\subsubsection{Adjoint Method}
\begin{figure}
	\centering
	\begin{subfigure}[b]{0.49\textwidth}
		\includegraphics[width=\textwidth]{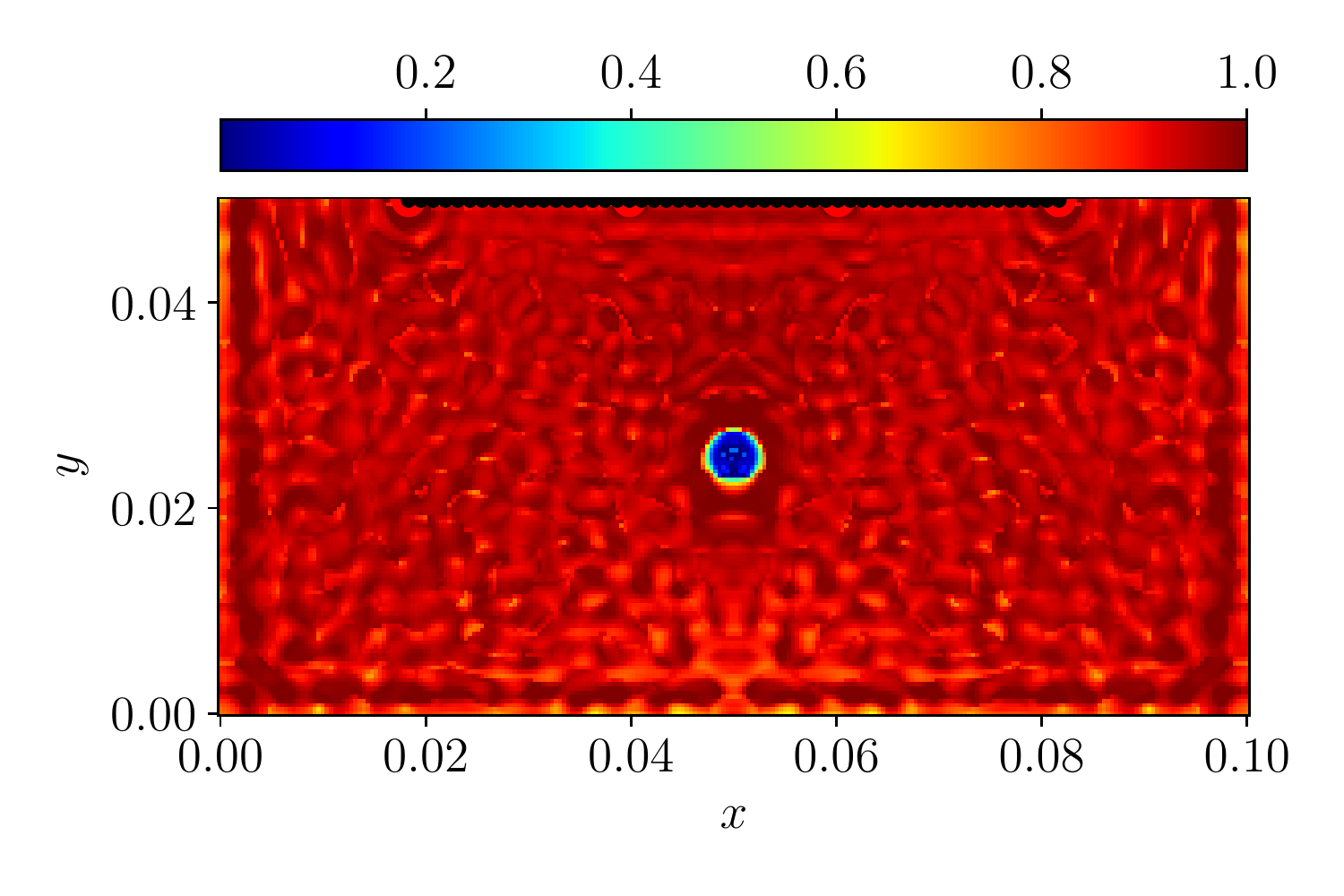}
		\caption{Prediction.}
	\end{subfigure}
	\hfill
	\begin{subfigure}[b]{0.49\textwidth}
		\includegraphics[width=\textwidth]{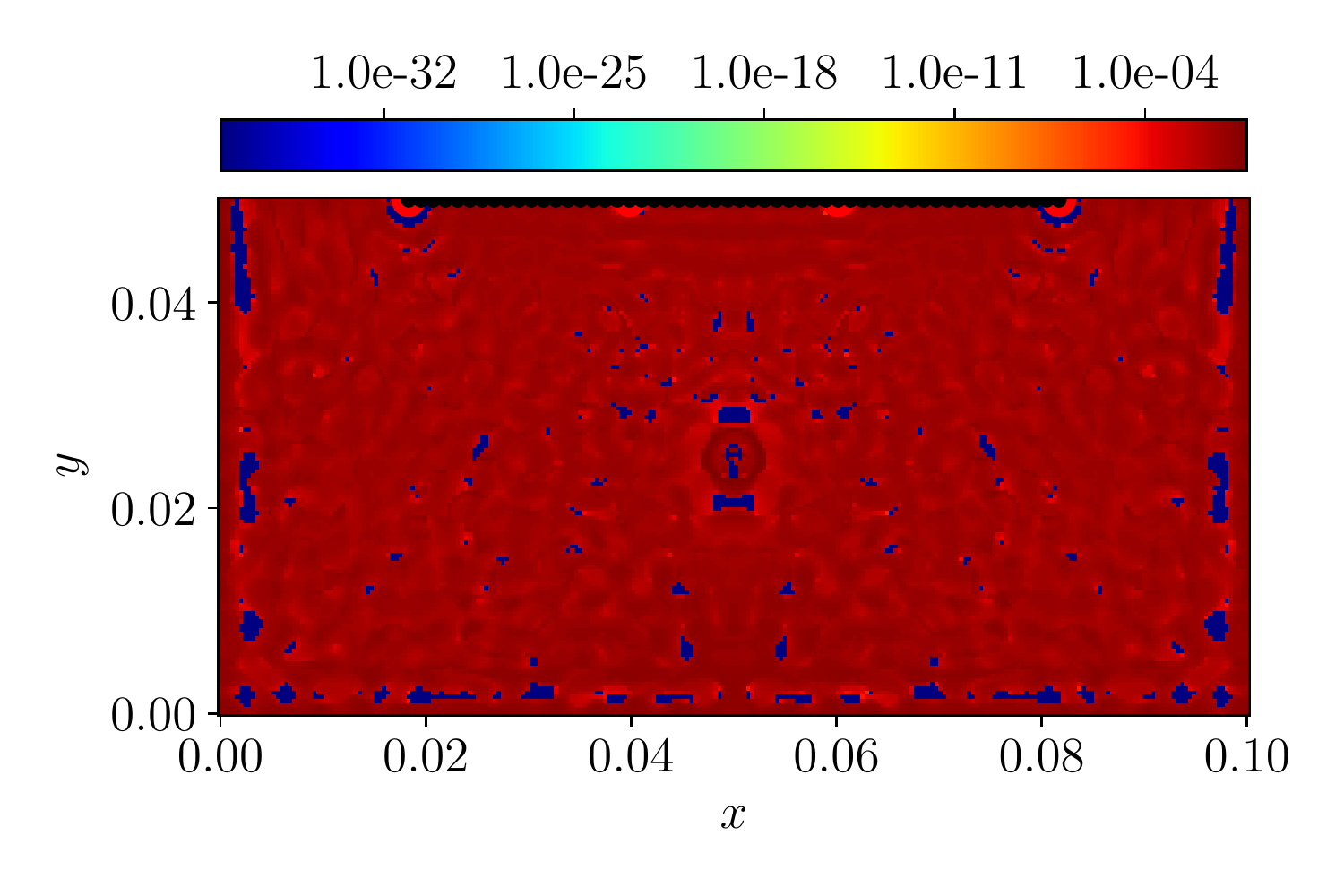}
		\caption{Absolute error.}
	\end{subfigure}
	\caption{Prediction and error of the material distribution after 50 epochs using automatic differentiation and a constant Ansatz. Average elapsed time per epoch is $2.46$ s.} \label{fig:constantAnsatz}
\end{figure}

Next, a performance evaluation of the adjoint method is provided in~\cref{tab:methods}. Section~\ref{sec:AdjointMethod} then contains more details and a special focus on the importance of the Ansatz and the gradient computation. A piece-wise constant Ansatz is considered instead of the neural network using automatic differentiation for the gradient computation. The void is clearly identified in the prediction after 50 epochs, shown in~\cref{fig:constantAnsatz}, and achieved with a slightly increased learning rate of $\alpha_{\min}=6\cdot 10^{-2}$. However, additional oscillatory artefacts are observed. This coincides with the observations in~\cite{buerchner_2022}, where linear shape functions are used as the Ansatz of the unknown material field.  However, interestingly changing the Ansatz of $\gamma$ to a neural network leads to a very different and better local minimum. Specifically, the mean squared error is reduced from $3.2\cdot 10^{-3}$ using the piecwise constant Ansatz for $\gamma$ as in~\cref{eq:PWConstAnsatz} to $9.5\cdot 10^{-4}$ using a neural network as in~\cref{eq:coefficient}. Furthermore, the smoothness of the solution is increased in regions where the true solution is constant, while the sharp gradients around the circular void are retained as depicted in~\cref{fig:inversionForwardSolver} and quantified in Appendix~\ref{sec:sharpness}. It is noteworthy that this increase in accuracy is achieved without a significant increase in computational effort. 

However, one major difference accompanying the change in Ansatz must be mentioned. The initial field of the constant Ansatz is chosen to be equal to one throughout the entire domain, which in the neural network case is impossible without pre-training. Pre-training to learn a constant field of one with a large regularization on the neural network was not found to be helpful. Thus the neural network starts with a random field, as later shown in the prediction histories in~\cref{fig:predictionHistories}. However, introducing a random field as initialization for the constant Ansatz is found to be detrimental, observable in an unsuccessful convergence toward a meaningful inversion. Thus, the strength of the neural network Ansatz must lie in its non-linear connectivity and strong over-parametrization, $526\,252$ parameters compared to the $252\times 124=31\,248$ parameters for the constant Ansatz. \\
It has to be noted that increasing the number of parameters in the constant Ansatz might be beneficial and yield a better comparison to the neural network Ansatz. However, when using finite differences, increasing the number of parameters, by e.g. increasing the number of voxels or the polynomial degree per finite difference grid point do not yield an improved indicator discretization. These methods are possible in the context of the finite element and cell method, where the added parameters can be taken into account via multiple Gauß points per element. Thus the benefit of neural networks as discretization becomes even clearer, as they provide a universal tool independent of the simulation method.

\begin{figure}
	\centering
	\begin{subfigure}[b]{0.49\textwidth}
		\includegraphics[width=\textwidth]{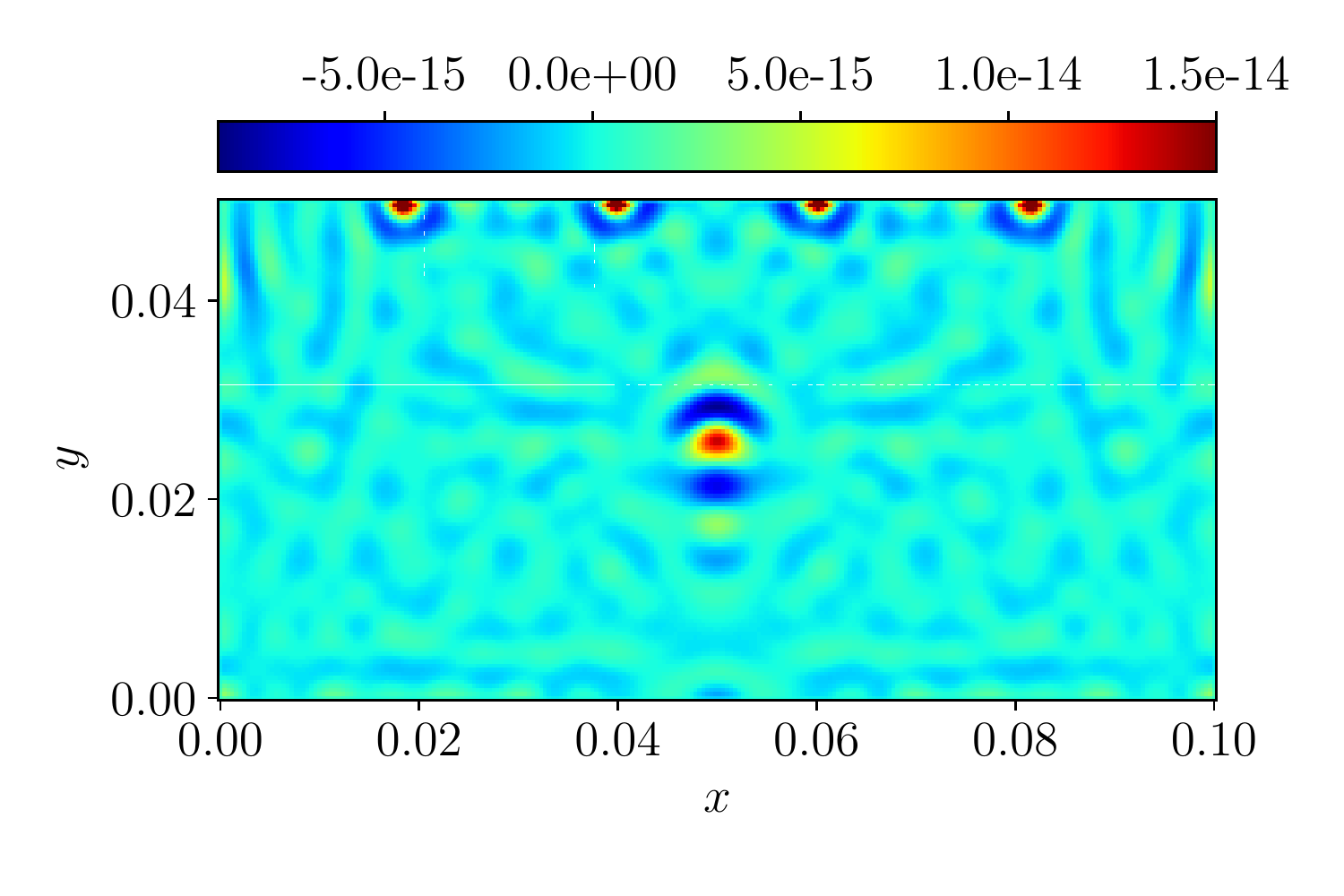}
		\caption{Gradient, automatic differentiation.}
	\end{subfigure}
	\hfill
	\begin{subfigure}[b]{0.49\textwidth}
		\includegraphics[width=\textwidth]{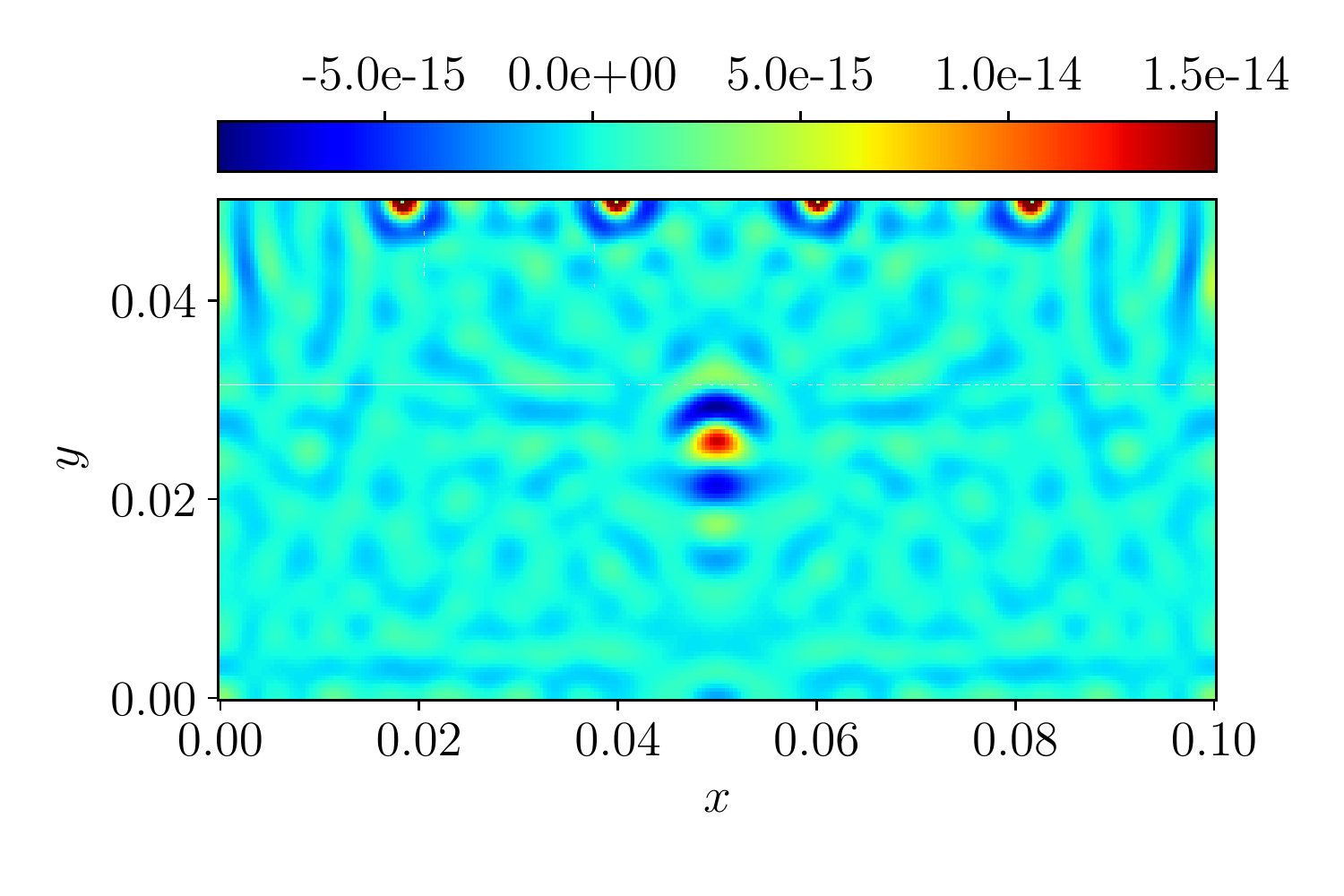}
		\caption{Gradient, adjoint method.}
	\end{subfigure} \\
	\begin{subfigure}[b]{0.49\textwidth}
		\includegraphics[width=\textwidth]{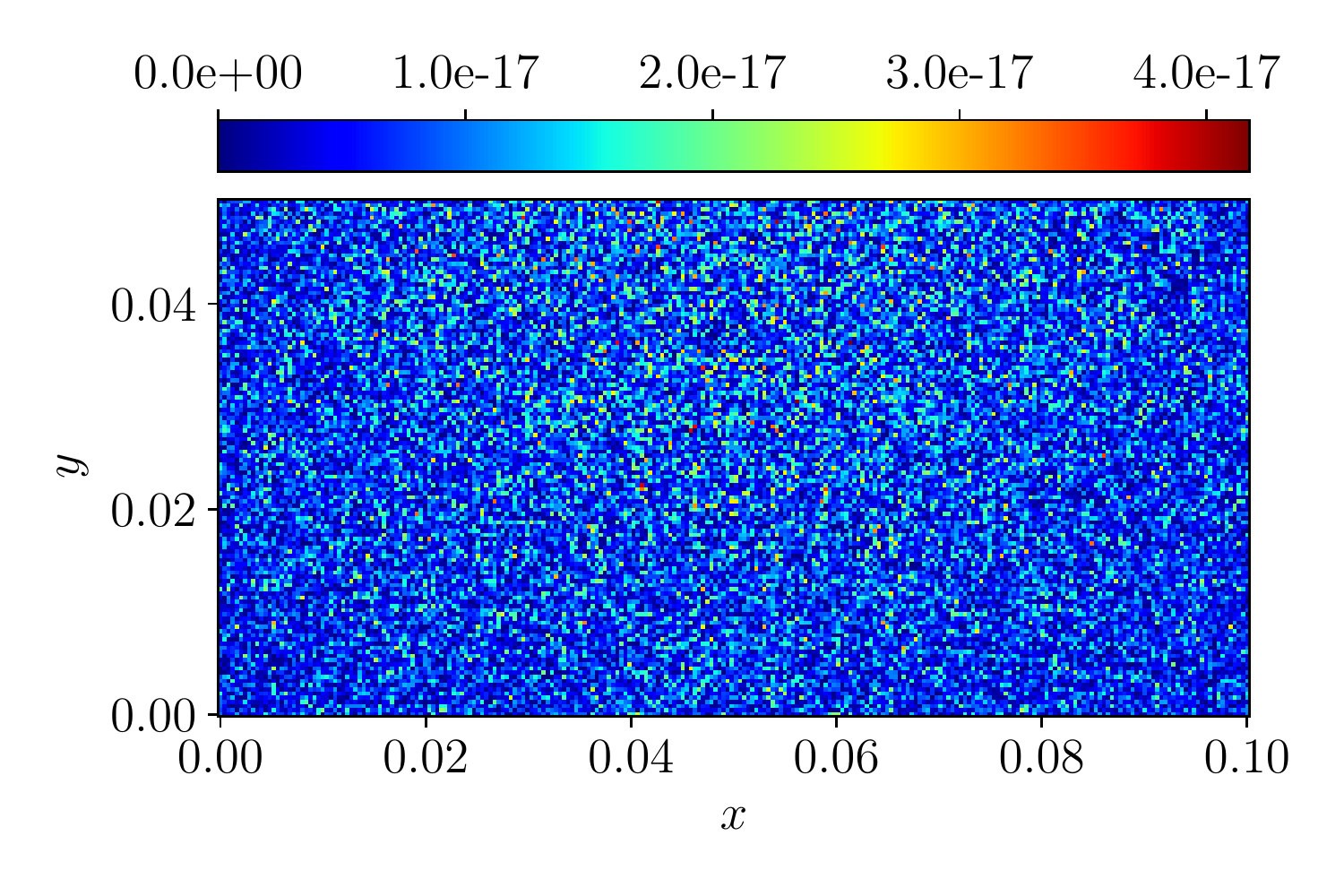}
		\caption{Absolute error, automatic differentiation.}
	\end{subfigure}
	\hfill
	\begin{subfigure}[b]{0.49\textwidth}
		\includegraphics[width=\textwidth]{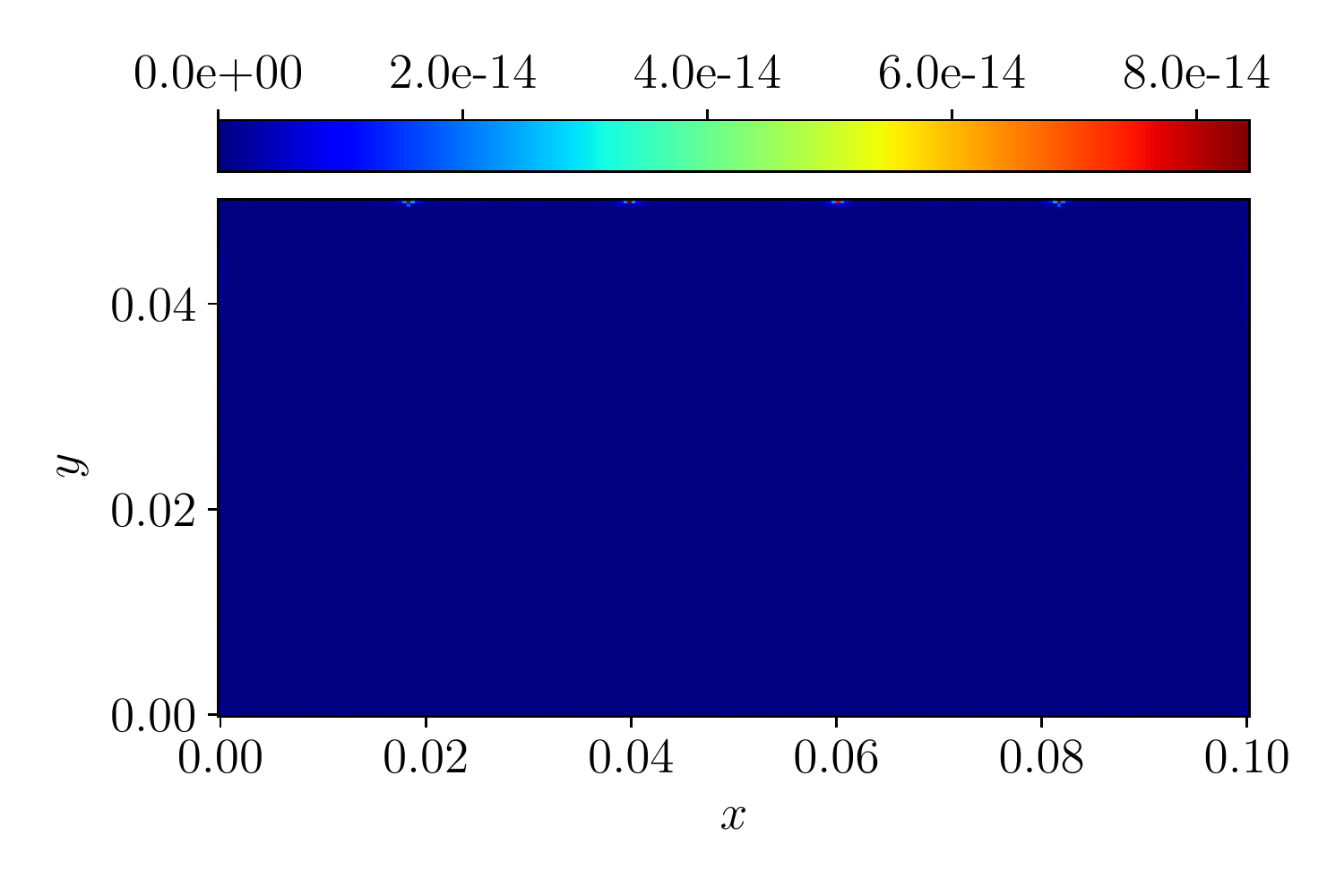}
		\caption{Absolute error, adjoint method.}
	\end{subfigure}
	\caption{Initial gradient using a constant Ansatz computed with automatic differentiation (a) and the adjoint method (b). The corresponding absolute errors in (c) and (d) are computed with a numerical reference solution obtained with finite differences.} \label{fig:initialGradientComparison}
\end{figure}

With this final modification, the change to the classical adjoint optimization as presented in~\cite{buerchner_2022} is almost complete, with the minor difference of using automatic differentiation for the gradient computation. However, as previously stated, moving from automatic differentiation to the continuous adjoint method for the gradient computation is not likely to have a significant impact. To assess this numerically, the initial gradient of the cost function w.r.t. the inverse field $\nabla_{\boldsymbol{c}} \mathcal{L}$ using a constant Ansatz is in both cases compared to a reference solution obtained with finite differences, by varying one parameter $c_i$ at a time. Only minor differences are observable at the boundaries, illustrated in~\cref{fig:initialGradientComparison}. The difference can be explained by the numerical evaluation of the Fréchet kernel,~\cref{eq:frechet}, which introduces the most significant numerical errors at the boundaries due to the central finite difference scheme. The absolute error w.r.t. the reference solution is, however, nonetheless negligible. Thus only minor differences in the learned prediction from~\cref{fig:constantAnsatz} are expected, as confirmed by the numerical comparison. The tiny differences accumulate throughout the iterations leading to slightly different local minima with a similar quality of inversion. Thus, no improvement is expected in the solution quality when using either automatic differentiation or the adjoint method for the gradient computation. However, using the adjoint method to compute the gradient $\nabla_{\hat{\gamma}} \mathcal{L}_{\mathcal{M}}$ leads to a significant drop in the computational effort from $2.46$ s to $1.05$ s per epoch for a simulation comparable to the one presented in~\cref{fig:constantAnsatz}.




\subsubsection{Hybrid Approach}
\begin{figure}
	\centering
	\begin{subfigure}[b]{0.49\textwidth}
		\includegraphics[width=\textwidth]{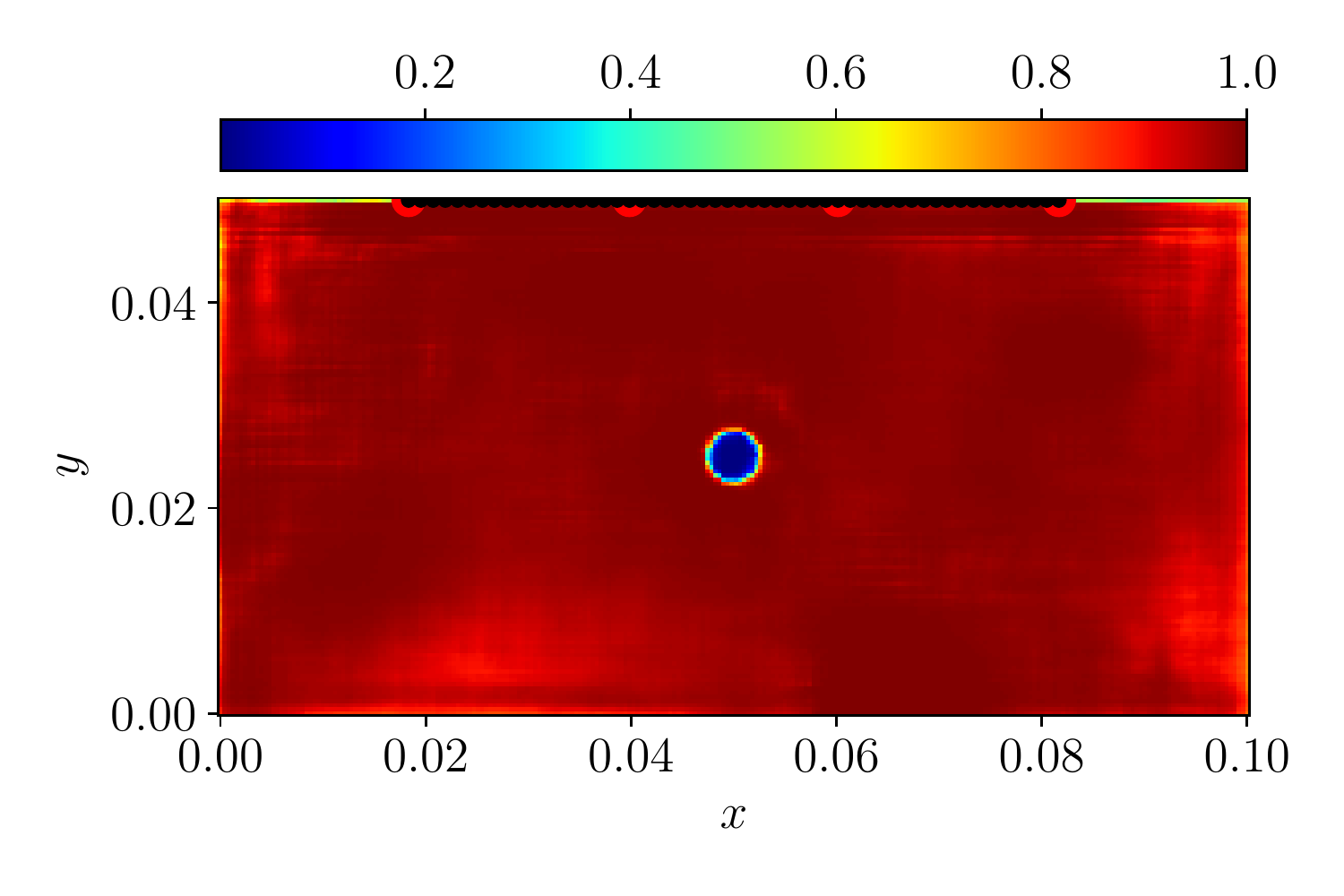}
		\caption{Prediction.}
	\end{subfigure}
	\hfill
	\begin{subfigure}[b]{0.49\textwidth}
		\includegraphics[width=\textwidth]{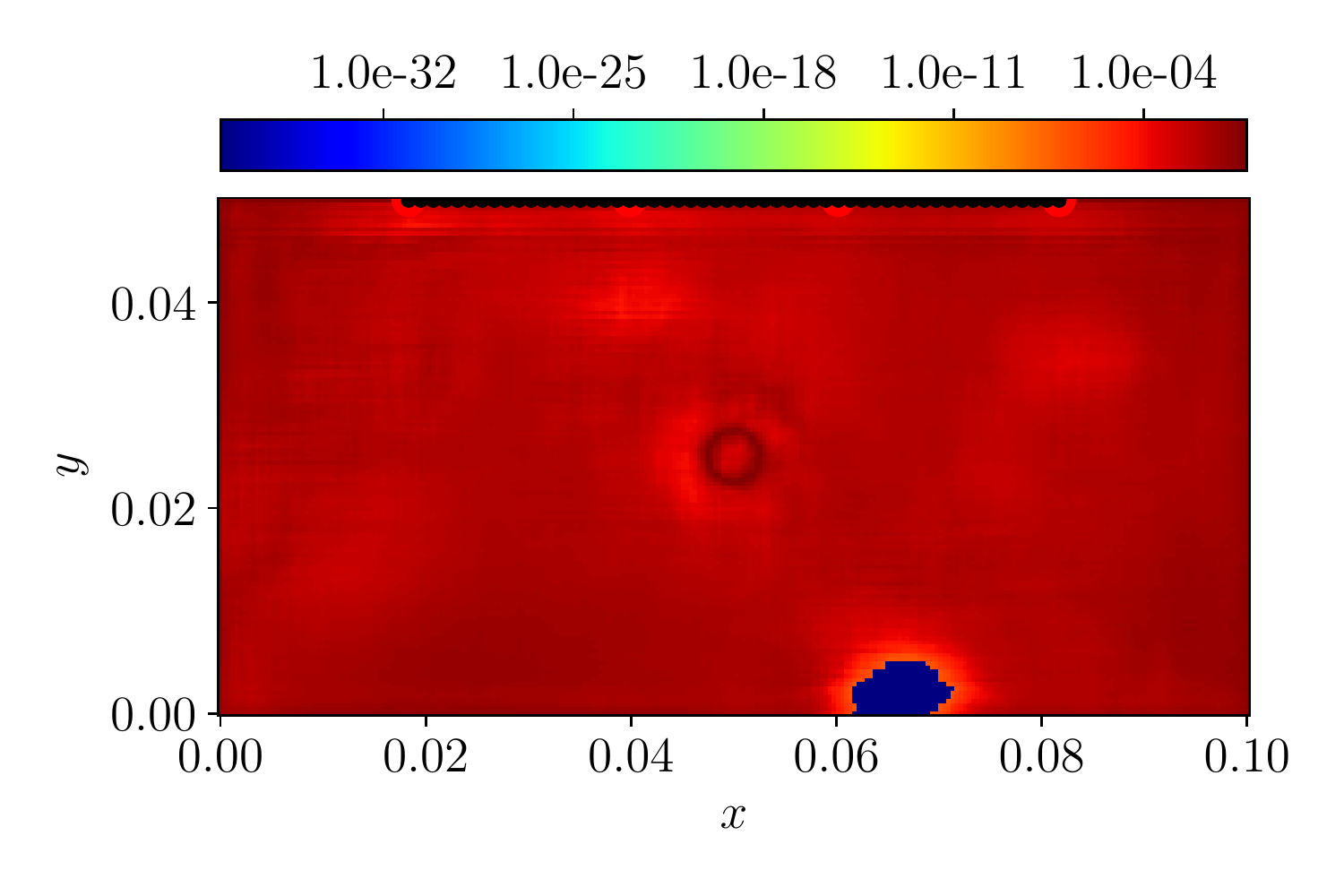}
		\caption{Absolute error.}
	\end{subfigure}
	\caption{Prediction and error of the material distribution after 50 epochs using the combination of automatic differentiation and adjoint method and a neural network Ansatz. Average elapsed time per epoch is $1.06$ s.} \label{fig:predictionHybrid}
\end{figure}

The proposed hybrid approach classified in~\cref{tab:methods} and described in detail in~\cref{sec:hybridApproach} combines the advantage of the neural network Ansatz with the lower computational effort of the adjoint method for the gradient computation $\nabla_{\hat{\gamma}} \mathcal{L}_{\mathcal{M}}$. The resulting prediction is presented in~\cref{fig:predictionHybrid} with only a minor decrease in prediction quality when compared to~\cref{fig:inversionForwardSolver}, while retaining computation times of $1.06$ s per epoch. A slight increase of the learning rate $\alpha_{\min}=4\cdot 10^{-3}$ was beneficial.

\subsubsection{Comparison}
\begin{figure}
	\centering
	\begin{subfigure}[b]{0.99\textwidth}
		\includegraphics[width=\textwidth]{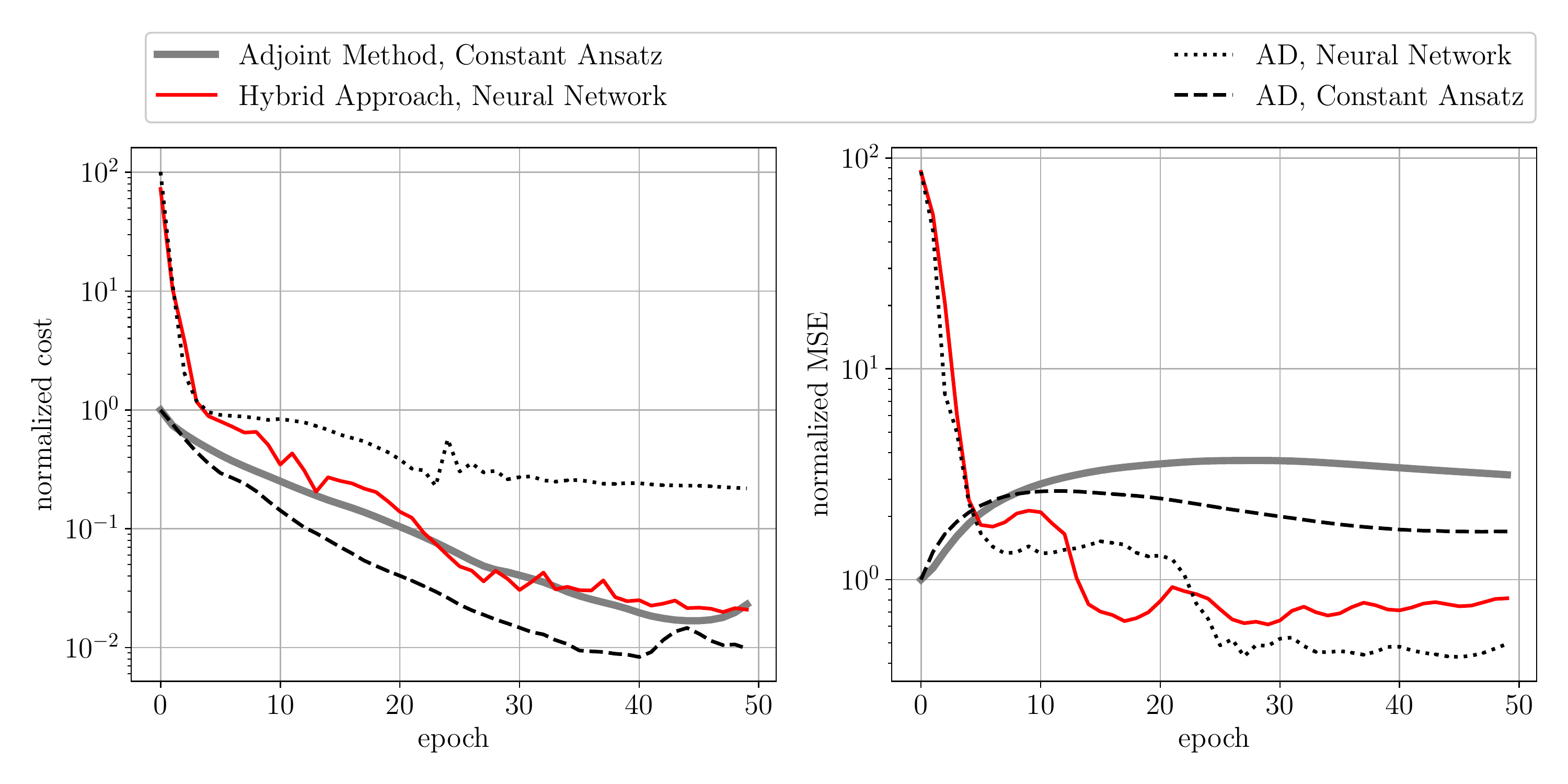}
	\end{subfigure}
	\caption{Comparison of the training histories for each presented approach. The cost function, which consists of the measurement loss $\mathcal{L}_{\mathcal{M}}$ is shown on the left, and the mean squared error (MSE) w.r.t. the true material distribution $\gamma(\boldsymbol{x})$ is on the right. The normalization of cost and mean squared error is with respect to the corresponding values with an indicator function equal to one in the entire domain.} \label{fig:learningcurves}
\end{figure}

To conclude the comparison between the investigated methodologies, the training histories are juxtaposed in~\cref{fig:learningcurves}. While the constant Ansatz enables the best minimization of the cost function, the mean squared error is lowest for the neural network Ansatz, confirming once again the advantage of the neural network Ansatz. Furthermore, the impact of using automatic differentiation or the adjoint method with finite differences becomes clearer with a consistently lower cost and mean squared error when using automatic differentiation. In view of the computational benefit, the slight loss in accuracy is acceptable. This difference only occurs due to the finite difference evaluations of the Fréchet Kernel,~\cref{eq:frechet}. In the case of finite elements, the Fréchet Kernel can be evaluated exactly by analytical differentiation of the shape functions. 

\begin{figure}
	\centering
	\begin{tikzpicture}
	\node at (0,0) (picA){
		\includegraphics[width=0.47\textwidth]{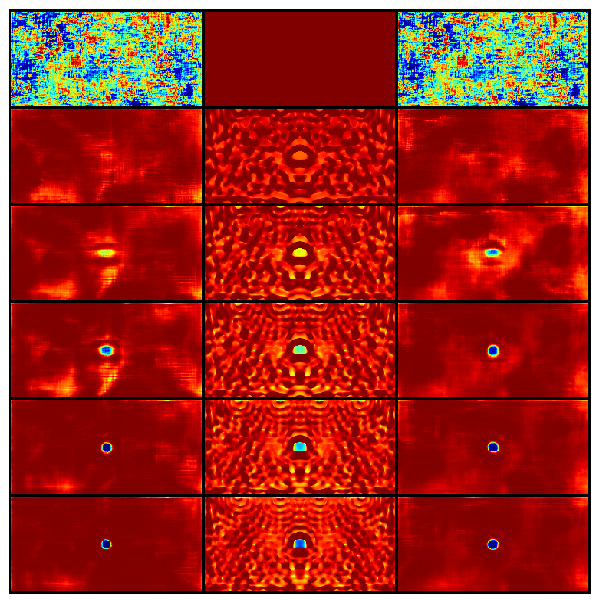}
	};
	\node at (-5,3.75-0.65) {epoch 0};
	\node at (-5,2.5-0.65) {epoch 5};
	\node at (-5,1.25-0.65) {epoch 10};
	\node at (-5,0-0.65) {epoch 15};
	\node at (-5,-1.25-0.65) {epoch 20};
	\node at (-5,-2.5-0.65) {epoch 25};
	
	\node at (-2.5,4.5) {\begin{tabular}{c} AD\\ neural network \end{tabular}};
	\node at (0,4.5) {\begin{tabular}{c} Adjoint\\ Method \end{tabular}};
	\node at (2.5,4.5) {\begin{tabular}{c} Hybrid\\ Approach \end{tabular}};
	
	\node at (0,-4.5) (picB) {
		\includesvg[width=0.45\textwidth]{code/output/colorbar.svg}	
	};
	
	\end{tikzpicture}
	\caption{Prediction histories starting at epoch 0 until epoch 25. From left to right, automatic differentiation with a neural network Ansatz, the adjoint method with a piece-wise constant Ansatz, and the hybrid approach with a neural network Ansatz are depicted.} \label{fig:predictionHistories}
\end{figure}

To gain an additional understanding of the training behaviour, consider the history of indicator predictions provided in~\cref{fig:predictionHistories}. Here, automatic differentiation with a neural network Ansatz, the adjoint method with a constant Ansatz and the hybrid approach are juxtaposed as snapshots showing the entire domain at specific instances of the iteration. In all three cases, the void is essentially recovered after 15 epochs. After the identification, only minor improvements in the remaining domain are made, where the effect is strongest for the hybrid approach and weakest for the adjoint method. Interestingly, all three approaches lead to three completely different local minima within a comparable amount of computational effort, despite initial noise in the case of the neural networks.

\begin{figure}
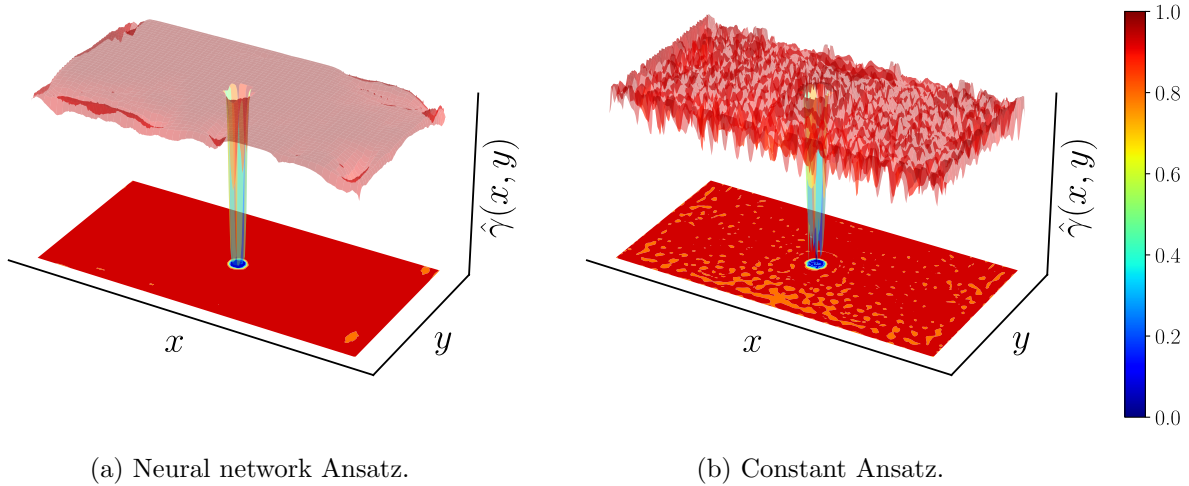

	\centering
	\begin{subfigure}[b]{0.44\textwidth}
		\hspace*{-3.4cm}
		\includesvg[width=1.8\textwidth]{code/output/IterativeSurfaceNeuralNetwork1.svg}
		\vspace*{-1.5cm}
		\caption{Neural network Ansatz.}
	\end{subfigure}
	\hfill
	\begin{subfigure}[b]{0.44\textwidth}
		\hspace*{-3.3cm}
		\includesvg[width=1.8\textwidth]{code/output/IterativeSurfaceConstant1.svg}
		\vspace*{-1.5cm}
		\caption{Constant Ansatz.}
	\end{subfigure}
	\hfill
	\begin{subfigure}[b]{0.09\textwidth}	
		\begin{tikzpicture}
	
		\node at (0,1) (picB) {
			\includesvg[width=0.51\textwidth]{code/output/colorbarRight.svg}	
		};	
	\end{tikzpicture}
	\caption*{ }
	\end{subfigure}
	\caption{Surface plot of the predicted material distribution after 50 epochs. The training was performed with a neural network and constant Ansatz using automatic differentiation.} \label{fig:surfaceplot}
\end{figure}

Thus two critical conclusions arise from this investigation. Firstly, learning the forward and inverse fields simultaneously for a single inversion problem yields a much more complex optimization task than using an non-trainable forward operator such as finite differences. The improvement of using a non-trainable forward operator is about four orders of magnitude when considering the number of iterations. Secondly, a neural network as Ansatz for the inverse fields outperforms classical approaches, such as a piece-wise constant Ansatz, with regard to the recovered field. To once again highlight the difference, a surface plot of the prediction of both the neural network Ansatz and the constant Ansatz after optimization using automatic differentiation is provided in~\cref{fig:surfaceplot}. The surfaceplot illustrates the difference in smoothness of the two solutions. Although the neural network solution is smoother, the sharpness retained (quantified in Appendix~\ref{sec:sharpness}). The plot indicates a regularization effect to be studied further using the three-dimensional case in the upcoming~\cref{sec:threeDim}.


\subsection{Three-Dimensional Case}\label{sec:threeDim}
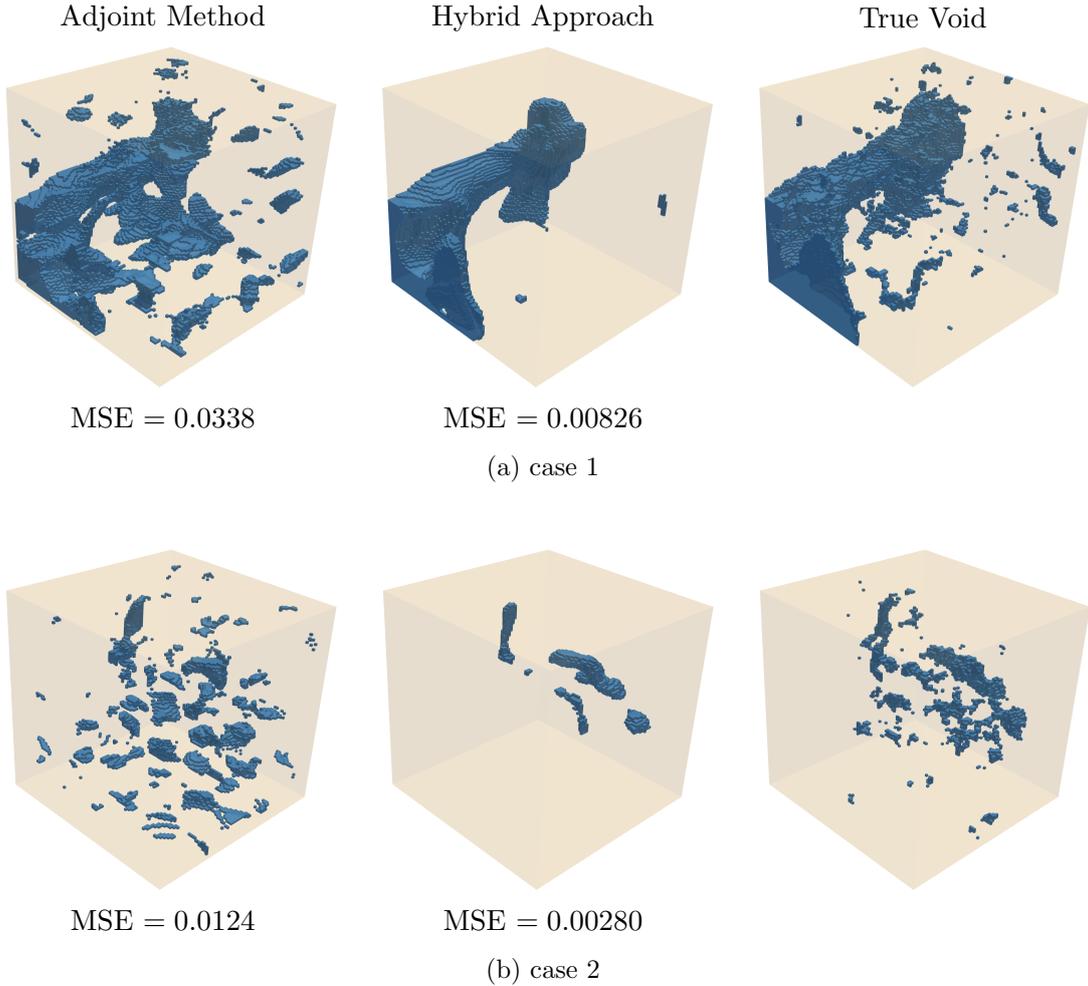
\begin{figure}
	\centering
	\begin{subfigure}[b]{1\textwidth}
		\centering
		\begin{tikzpicture}
		\node at (0,0) (picA){
			\includegraphics[width=0.9\textwidth]{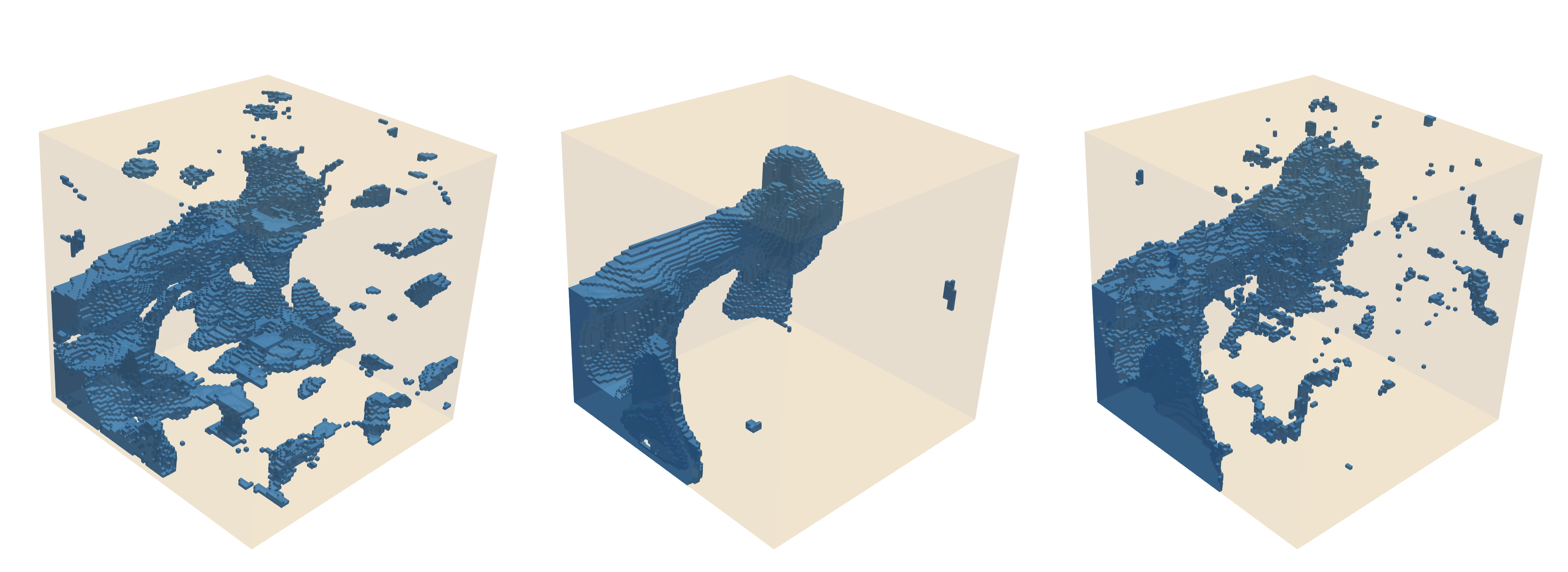}
		};
		\node at (-5,2.4) {Adjoint Method};
		\node at (0,2.4) {Hybrid Approach};
		\node at (5,2.4) {True Void};
		
		\node at (-5, -2.9) {MSE $=0.0338$};
		\node at (0,-2.9) {MSE $=0.00826$};
		\end{tikzpicture}
		\caption{case 1}
	\end{subfigure}
	\begin{subfigure}[b]{1\textwidth}
		\centering
		\begin{tikzpicture}
		\node at (0,0) (picA){
			\includegraphics[width=0.9\textwidth]{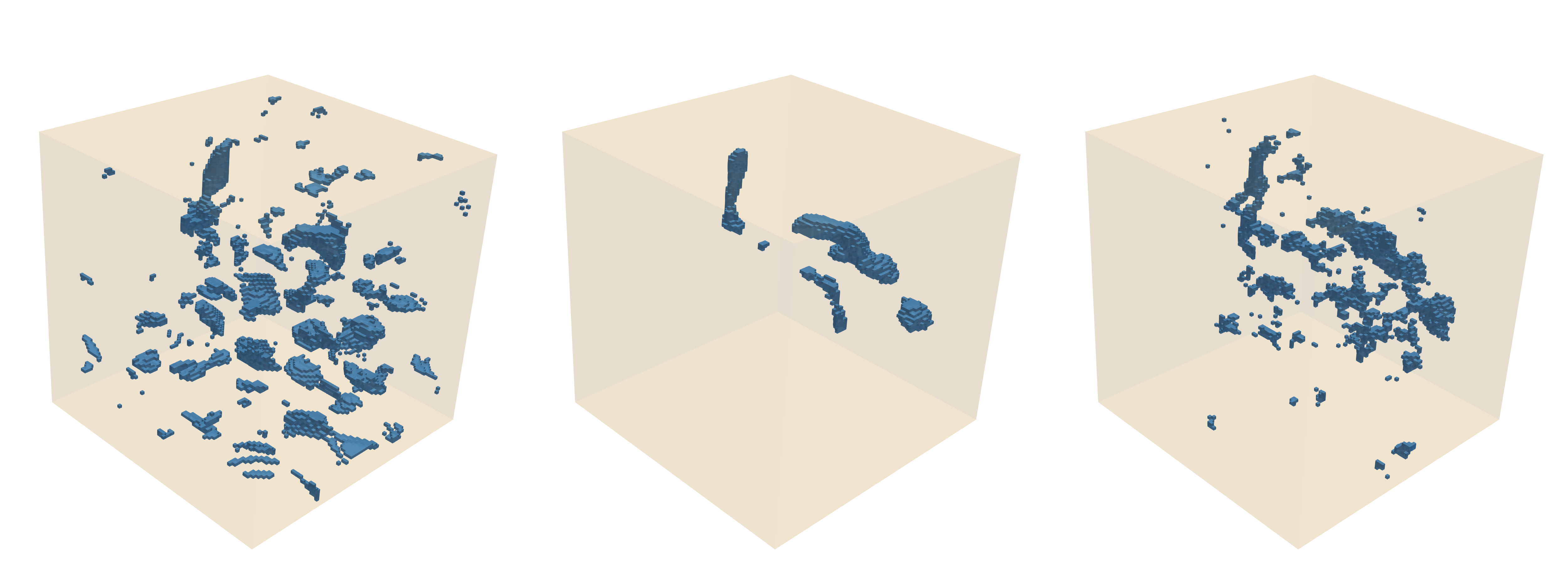}
		};
		\node at (-5,-2.9) {MSE $=0.0124$};
		\node at (0,-2.9) {MSE $=0.00280$};
		\end{tikzpicture}
		\caption{case 2}	
	\end{subfigure}
	\caption{Identified void by the adjoint method and the hybrid approach after 100 epochs in comparison to the true void. The illustrated void is extracted by thresholding indicator values lower than $0.5$. Average elapsed time per epoch is about $29.1$ s for both methods.} \label{fig:3dresult}
\end{figure}

\begin{figure}
	\centering
	\begin{subfigure}[b]{0.09\textwidth}		
		\begin{tikzpicture}
		\node [white] at (-3.7,-3) {a};
		\node at (-3.7,2.3) {\begin{tabular}{c}Adjoint\\Method\end{tabular}};
		\node at (-3.7,0+0.4) {\begin{tabular}{c}Hybrid\\Approach\end{tabular}};
		\node at (-3.7,-2.1+0.5) {\begin{tabular}{c}True\\Indicator\end{tabular}};		
		\end{tikzpicture}
	\end{subfigure}
	\hfill
	\begin{subfigure}[b]{0.35\textwidth}
		\begin{tikzpicture}
		\node at (0,0) (picA){
			\includegraphics[width=\textwidth]{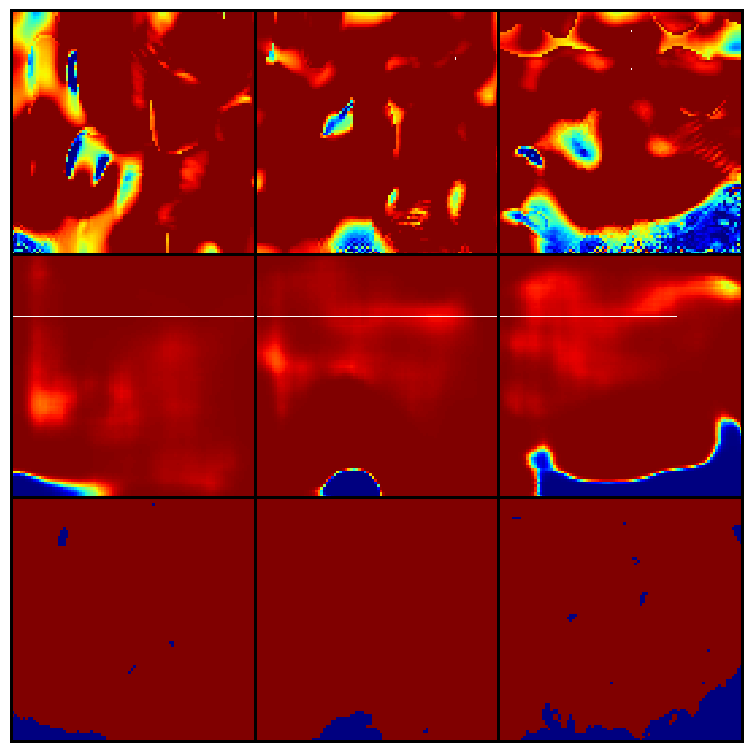}
		};
		\node at (-1.5,3.1) {$xy$};
		\node at (0,3.1) {$xz$};
		\node at (1.5,3.1) {$yz$};
		
		\end{tikzpicture}
		\caption{case 1}
	\end{subfigure}
	\hfill
	\begin{subfigure}[b]{0.35\textwidth}
		\begin{tikzpicture}
		\node at (0,0) (picA){
			\includegraphics[width=\textwidth]{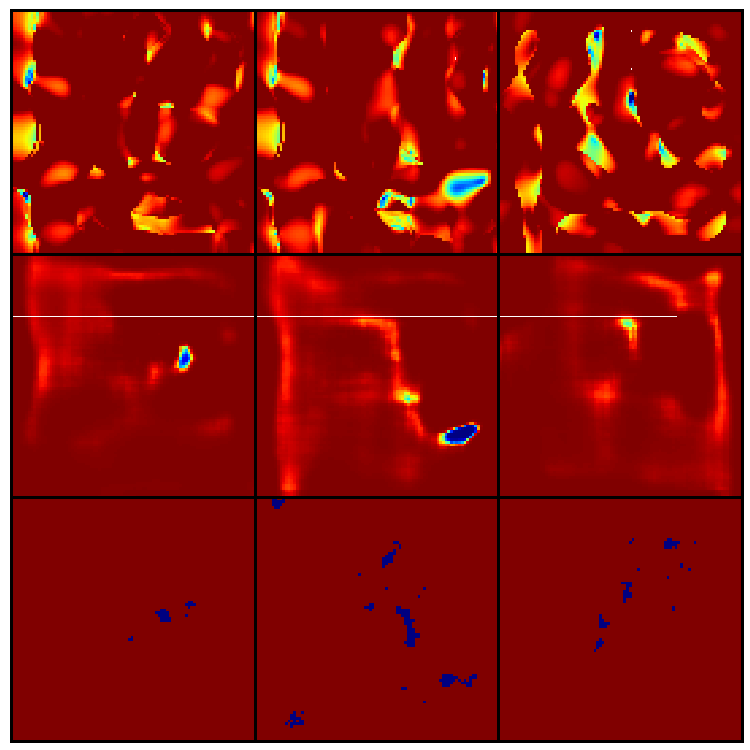}
		};	
		
		\node at (-1.5,3.1) {$xy$};
		\node at (0,3.1) {$xz$};
		\node at (1.5,3.1) {$yz$};
		
		\end{tikzpicture}
		\caption{case 2}
	\end{subfigure}
	\hfill
	\begin{subfigure}[b]{0.09\textwidth}	
		\begin{tikzpicture}
		
		\node at (0,1) (picB) {
			\includesvg[width=0.51\textwidth]{code/output/colorbarRight.svg}	
		};	
		\end{tikzpicture}
		\caption*{ }
	\end{subfigure}

	\caption{Center slices of the indicator function predictions using the adjoint method and the hybrid approach after 100 epochs. The sliced planes are annotated in the top. Average elapsed time per epoch is $29.1$ s for the hybrid approach and $29.1$ s for the adjoint method.} \label{fig:3dslices}
\end{figure}

Due to the quality of the inversion results using neural networks as Ansatz, a three-dimensional case is considered. As the additional spatial dimension introduces an additional computational burden, only the most efficient and simultaneously most promising methods are compared, i.e. the hybrid method and the adjoint method with a piece-wise constant Ansatz. The inversion results of the two approaches are illustrated in~\cref{fig:3dresult} in terms of identified voids for the two CT-scan slices from~\cref{fig:3dcase}. Once again, the neural networks enable a superior reconstruction, quantifiable by considering the mean squared error provided as annotations in~\cref{fig:3dresult}. Center slices are shown in~\cref{fig:3dslices} to provide an insight into the overall distribution and smoothness of the solution. From this, it can be observed that the reconstruction of the adjoint method suffers from the oscillatory artefacts to such an extent that the voids are no longer clearly identified. The hybrid approach exhibits much less artefacts and is able to recover the voids with superior accuracy. 

\begin{figure}
	\centering
	\includegraphics[width=0.8\textwidth]{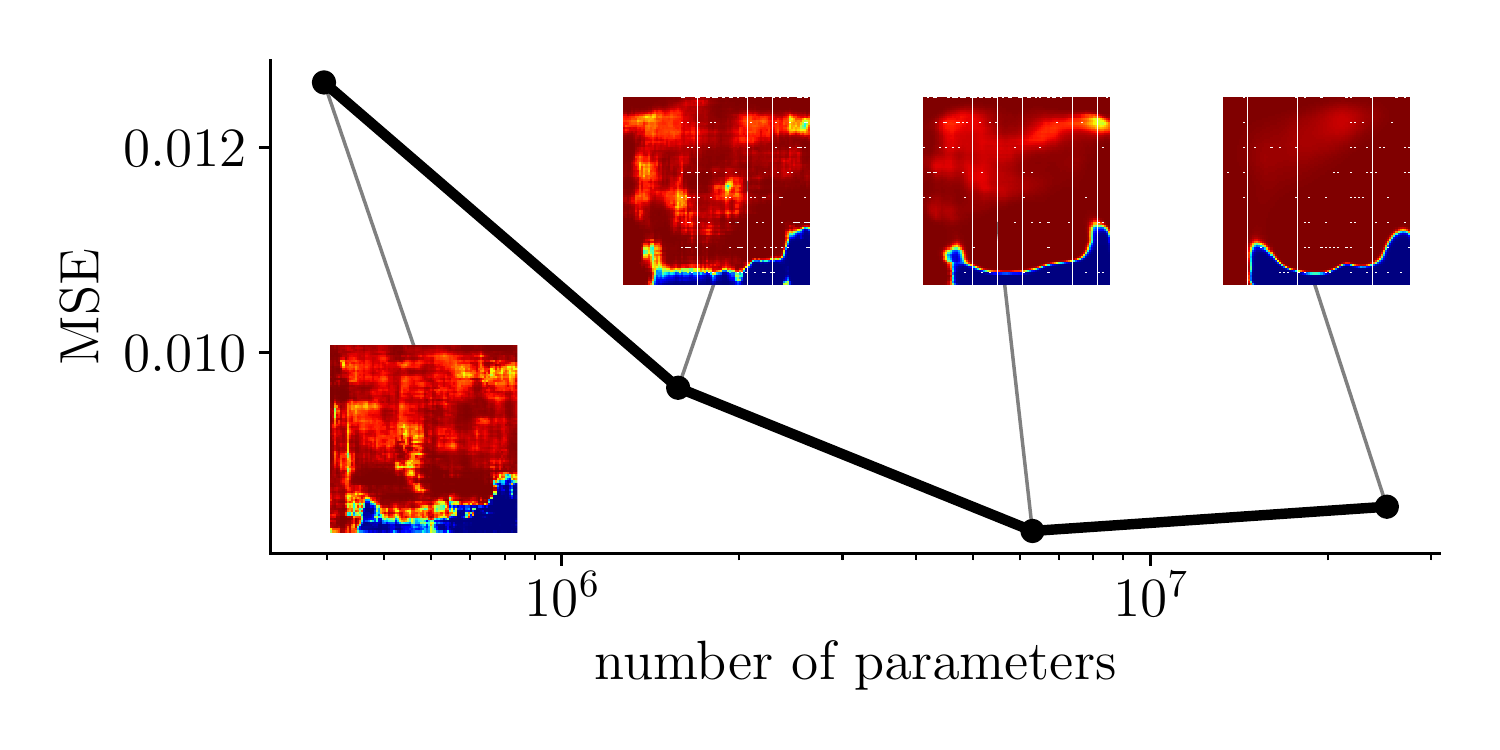}
	\caption{Mean squared error after 100 epochs for different number of parameters illustrated on case 1.}\label{fig:numparameters}
\end{figure}

We attribute this superiority to the different characteristics of the neural network Ansatz which not only leads to an over-parametrization of the indicator field but also has an in-built non-linearity. Specifically, $6\,306\,764$ parameters were used for the neural network compared to the $92^3=778\,688$ parameters in the classical discretization. Interestingly, the reconstruction is worse when using fewer neural network parameters. By using half and quarter as many filters, yielding only $1\,577\,452$ and $394\,748$ parameters, the quality of the reconstruction decreases. At the same time, using double the amount of filters, resulting in $25\,221\,004$ parameters, leads to a minor decrease in quality, due to the now dominating regularization effect of the over-parameterization. The mean squared error (MSE) w.r.t. the true material distribution presented for the three variations of the network in~\cref{fig:numparameters} quantifies this observation. Furthermore, a smoothing effect is noticeable, with an increase of parameters indicating a regularization-like behaviour. Thus, the over-parametrization enabled by the neural network is one, yet not the only key ingredient in enhancing the inversion results.

Several other aspects distinguish the neural network from the piece-wise constant Ansatz. The non-linear activations inside the neural network represent an important difference. Also, the neural network acts as a global function, while the constant Ansatz is defined on local patches. Furthermore, other network architectures, such as fully-connected or graph neural networks, may show different behaviour. Currently, it is unclear which aspects are critical aside from the over-parametrization. All aspects can be studied in a straightforward empirical manner and are subject to further research.

\section{Conclusion}\label{sec:conclusion}

%
%
%
%
%
%
%
%
%
%

The aim of this work was to compare PINNs with the classical adjoint optimization for full waveform inversion and thereby quantify the strengths of each approach by incremental modification. Firstly, estimating the forward solution via a nested optimization using two neural networks was shown to be disadvantageous. Instead, a classical forward solver, as in the adjoint optimization, is recommended. Secondly, a neural network as Ansatz of the inverse field is highly beneficial with a regularization-like effect without a loss of sharpness. Thirdly, no significant difference in accuracy exists in the sensitivity computation when comparing automatic differentiation and the continuous adjoint method. However, the continuous adjoint method is preferred in the gradient estimation of the cost w.r.t. the inverse field due to the lower computational effort. 

As a consequence of these insights, the hybrid method is proposed, utilizing a neural network to discretize the inverse field and computing the cost function gradient w.r.t. the inverse field with the adjoint method, which then is passed on to the backpropagation algorithm to update the neural network parameters. Comparable computational times are achieved, while the inverse reconstructions improved.

The performance of neural networks as the discretization of the inverse fields enables not only further research for other inverse problems but also a thorough investigation into the underlying cause of its success. Insight can be gained by considering the difference between global and local functions, different over-parametrizations, i.e. the number of parameters, the non-linearities in the network, different network architectures, and possibly a thorough investigations of the initialization of the wavefield and its effect on the inversion.

\section*{Acknowledgements}
The authors gratefully acknowledge the funding through the joint research project Geothermal-Alliance Bavaria (GAB) by the Bavarian State Ministry of Science and the Arts (StMWK) wich finance Leon Herrmann as well as the Georg Nemetschek Institut (GNI) under the project DeepMonitor. Furthermore, we gratefully acknowledge funds received by the Deutsche Forschungsgemeinschaft under Grant no. KO 4570/1-1 and RA 624/29-1 which support Tim Bürchner. Felix Dietrich would like to acknowledge the funds received by the Deutsche Forschungsgemeinschaft - project no. 468830823.

\section*{Declarations}
\textbf{Conflict of interest} No potential conflict of interest was reported by
the authors.

\section*{Data Availability}
We provide a PyTorch~\cite{paszke_pytorch_2019} implementation of all methods in \cite{herrmann_code}.

\newpage

\appendix
\renewcommand\thesection{\Alph{section}}

\section{Efficient Finite Difference Implementation in Deep Learning Frameworks}\label{sec:appendixFD}
The explicit finite difference scheme, shown for the two-dimensional case,
\footnotesize
\begin{align*}
u_{i,j}^{n+1}=-u_{i,j}^{n-1}+2u_{i,j}^n\\
+\frac{2}{\gamma_{i,j}}\left(\frac{c_0\Delta t}{\Delta x}\right)^2\left( \left(\frac{1}{\gamma_{i,j}}+\frac{1}{\gamma_{i+1,j}}\right)^{-1}(u_{i+1,j}^n-u_{i,j}^n)-\left(\frac{1}{\gamma_{i-1,j}}+\frac{1}{\gamma_{i,j}}\right)^{-1}(u^n_{i,j}-u_{i-1,j}^n)\right)\\
+\frac{2}{\gamma_{i,j}}\left(\frac{c_0\Delta t}{\Delta y}\right)^2\left( \left(\frac{1}{\gamma_{i,j}}+\frac{1}{\gamma_{i,j+1}}\right)^{-1}(u_{i,j+1}^n-u_{i,j}^n)-\left(\frac{1}{\gamma_{i,j-1}}+\frac{1}{\gamma_{i,j}}\right)^{-1}(u^n_{i,j}-u_{i,j-1}^n)\right)\\
+\frac{\Delta t^2}{\rho_0 \gamma_{i,j}}f_{i,j}^n
\end{align*}
\normalsize
was implemented in PyTorch \cite{paszke_pytorch_2019} using convolutional kernels $\boldsymbol{K}$, thus exploiting PyTorch's GPU capabilities.
\begin{align*}
\boldsymbol{u}^n=-\boldsymbol{u}^n+2\boldsymbol{u}^{n-1}\\
+\frac{2}{\boldsymbol{\gamma}}\Bigg( \frac{1}{\frac{1}{\boldsymbol{\gamma}}*\boldsymbol{K}_{\gamma_{x_0}}}(\boldsymbol{u}*\boldsymbol{K}_{u_{x_0}})-\frac{1}{\frac{1}{\boldsymbol{\gamma}}*\boldsymbol{K}_{\gamma_{x_1}}}(\boldsymbol{u}*\boldsymbol{K}_{u_{x_1}})\\
+\frac{1}{\frac{1}{\boldsymbol{\gamma}}*\boldsymbol{K}_{\gamma_{y_0}}}(\boldsymbol{u}*\boldsymbol{K}_{u_{y_0}})-\frac{1}{\frac{1}{\boldsymbol{\gamma}}*\boldsymbol{K}_{\gamma_{y_1}}}(\boldsymbol{u}*\boldsymbol{K}_{u_{y_1}})\\
+\frac{\Delta t^2}{2 \rho_0}\boldsymbol{f}\Bigg)
\end{align*}

The corresponding kernels are defined as
\begin{align*}
\boldsymbol{K}_{\gamma_{x_0}}=\begin{pmatrix}
0 & 0 & 0 \\
0 & 1 & 1 \\
0 & 0 & 0
\end{pmatrix},
\boldsymbol{K}_{u_{x_0}}=\begin{pmatrix}
0 & 0 & 0 \\
0 & -\left(\frac{c_0\Delta t}{\Delta x}\right)^2 & \left(\frac{c_0\Delta t}{\Delta x}\right)^2 \\
0 & 0 & 0
\end{pmatrix}, \\
\boldsymbol{K}_{\gamma_{x_1}}=\begin{pmatrix}
0 & 0 & 0 \\
1 & 1 & 0 \\
0 & 0 & 0
\end{pmatrix},
\boldsymbol{K}_{u_{x_1}}=\begin{pmatrix}
0 & 0 & 0 \\
-\left(\frac{c_0\Delta t}{\Delta x}\right)^2 & \left(\frac{c_0\Delta t}{\Delta x}\right)^2 & 0 \\
0 & 0 & 0
\end{pmatrix}, \\
\boldsymbol{K}_{\gamma_{y_0}}=\begin{pmatrix}
0 & 0 & 0 \\
0 & 1 & 0 \\
0 & 1 & 0
\end{pmatrix},
\boldsymbol{K}_{u_{y_0}}=\begin{pmatrix}
0 & 0 & 0 \\
0 & -\left(\frac{c_0\Delta t}{\Delta y}\right)^2 & 0 \\
0 & \left(\frac{c_0\Delta t}{\Delta y}\right)^2 & 0
\end{pmatrix}, \\
\boldsymbol{K}_{\gamma_{y_1}}=\begin{pmatrix}
0 & 1 & 0 \\
0 & 1 & 0 \\
0 & 0 & 0
\end{pmatrix},
\boldsymbol{K}_{u_{y_1}}=\begin{pmatrix}
0 & -\left(\frac{c_0\Delta t}{\Delta y}\right)^2 & 0 \\
0 & \left(\frac{c_0\Delta t}{\Delta y}\right)^2 & 0 \\
0 & 0 & 0
\end{pmatrix}
\end{align*}
The homogeneous Neumann boundary conditions were enforced via ghost cells (see e.g.~\cite{langtangen_finite_2017}) and implemented by using pooling layers, and manual adjustments according to the ghost cells. \\

\section{Network Architecture}\label{sec:appendixNN}
A quantitative summary of the presented network from~\cref{fig:neuralnetwork} is provided in~\cref{tab:2dNN} for the two-dimensional case and~\cref{tab:3dNN} for the three-dimensional case. The choice of network hyperparameters varies slightly for the different approaches and cases. The interested reader is referred to the code made available at \cite{herrmann_code}.


\begin{table}[h!]
	\caption{Two-dimensional convolutional neural network architecture for the prediction of the indicator $\gamma$. The total number of parameters is $526\,252$.} \label{tab:2dNN}
	\centering
	\begin{tabular}{lll}
		\hline
		\textbf{layer}                          & \textbf{shape after layer}        & \textbf{learnable parameters} \\ \hline
		input                          & $128\times 8\times 4$    & 0                    \\ \hline \hline
		upsample                       & $128\times 16\times 8$    & 0                    \\ \hline
		2D convolution \& PReLU         & $128\times 16\times 8$    & $147\,584 + 1$          \\ \hline
		2D convolution \& PReLU         & $128\times 16\times 8$    & $147\,584 + 1$          \\ \hline \hline
		upsample                       & $128\times 32\times 16$  & 0                    \\ \hline
		2D convolution \& PReLU        & $64\times 32\times 16$   & $73\,792 + 1$            \\ \hline
		2D convolution \& PReLU        & $64\times 32\times 16$   & $36\,928 + 1$               \\ \hline \hline
		upsample                       & $64\times 64\times 32$   & 0                    \\ \hline
		2D convolution \& PReLU        & $64\times 64\times 32$   & $36\,928 + 1$               \\ \hline
		2D convolution \& PReLU        & $64\times 64\times 32$   & $36\,928 + 1$               \\ \hline \hline
		upsample                       & $64\times 128\times 64$   & 0                    \\ \hline
		2D convolution \& PReLU        & $32\times 128\times 64$   & $18\,464 + 1$              \\ \hline
		2D convolution \& PReLU        & $32\times 128\times 64$   & $9\,248 + 1$               \\ \hline \hline
		upsample                       & $32\times 256\times 128$ & 0                    \\ \hline
		2D convolution \& PReLU        & $32\times 256\times 128$ & $9\,248 + 1$                \\ \hline
		2D convolution \& PReLU        & $32\times 256\times 128$ & $9\,248 + 1$               \\ \hline \hline
		2D convolution without padding & $1\times 254\times 126$ & $289 + 1$ \\ 
		\& adaptive Sigmoid & &                   \\ \hline
	\end{tabular}
\end{table}

\begin{table}[h!]
	\caption{Three-dimensional convolutional neural network architecture for the prediction of the indicator $\gamma$. The number of filters per layer are double compared to the two-dimensional convolutional neural network. The total number of parameters is $6\,306\,764$.} \label{tab:3dNN}
	\centering
	\begin{tabular}{lll}
		\hline
		\textbf{layer}                          & \textbf{shape after layer}        & \textbf{learnable parameters} \\ \hline
		input                          & $256\times 3\times 3 \times 3$    & 0                    \\ \hline \hline
		upsample                       & $256\times 6\times 6 \times 6$    & 0                    \\ \hline
		3D convolution \& PReLU         & $256\times 6\times 6 \times 6$    & $1\,769\,728 + 1$          \\ \hline
		3D convolution \& PReLU         & $256\times 6\times 6 \times 6$    & $1\,769\,728 + 1$          \\ \hline \hline
		upsample                       & $256\times 12\times 12 \times 12$  & 0                    \\ \hline
		3D convolution \& PReLU        & $128\times 12\times 12 \times 12$   & $884\,864 + 1$            \\ \hline
		3D convolution \& PReLU        & $128\times 12\times 12 \times 12$   & $442\,496 + 1$               \\ \hline \hline
		upsample                       & $128\times 24\times 24 \times 24$   & 0                    \\ \hline
		3D convolution \& PReLU        & $128\times 24\times 24 \times 24$   & $442\,496 + 1$               \\ \hline
		3D convolution \& PReLU        & $128\times 24\times 24 \times 24$   & $442\,496 + 1$               \\ \hline \hline
		upsample                       & $128\times 48\times 48 \times 48$   & 0                    \\ \hline
		3D convolution \& PReLU        & $64\times 48\times 48 \times 48$   & $221\,248 + 1$              \\ \hline
		3D convolution \& PReLU        & $64\times 48\times 48 \times 48$   & $110\,656 + 1$               \\ \hline \hline
		upsample                       & $64\times 96\times 96 \times 96$ & 0                    \\ \hline
		3D convolution \& PReLU        & $64\times 96\times 96 \times 96$ & $110\,656 + 1$                \\ \hline
		3D convolution \& PReLU        & $64\times 96\times 96 \times 96$ & $110\,656 + 1$               \\ \hline \hline
		3D convolution without padding & $1\times 94\times 94 \times 94$ & $1\,729 + 1$ \\ 
		\& adaptive Sigmoid & &                   \\ \hline
	\end{tabular}
\end{table}

\section{Sharpness Quantification}\label{sec:sharpness}

\begin{figure}[htb] 
	\centering
	\begin{subfigure}[b]{0.45\textwidth}
		\includegraphics[width=\textwidth]{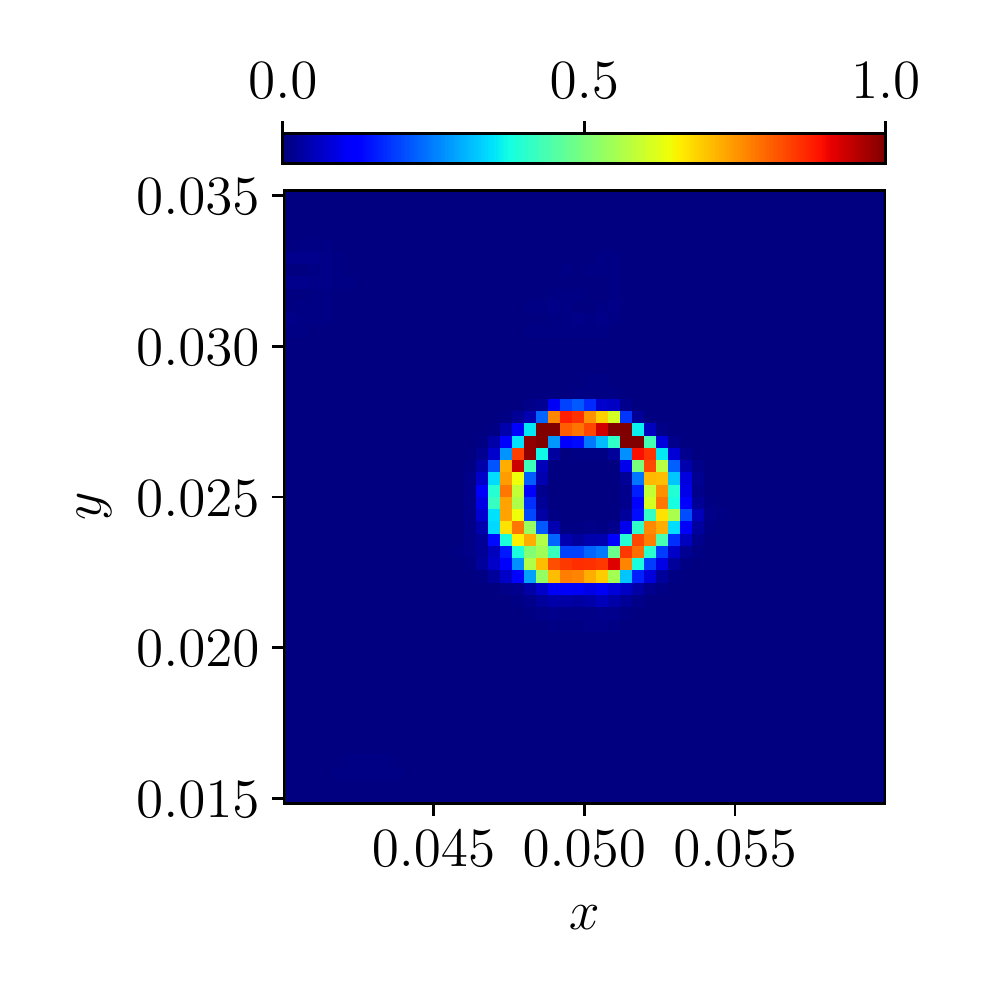}
		\caption{Neural network Ansatz from~\cref{fig:inversionForwardSolver}.}
	\end{subfigure}
	\begin{subfigure}[b]{0.45\textwidth}
		\includegraphics[width=\textwidth]{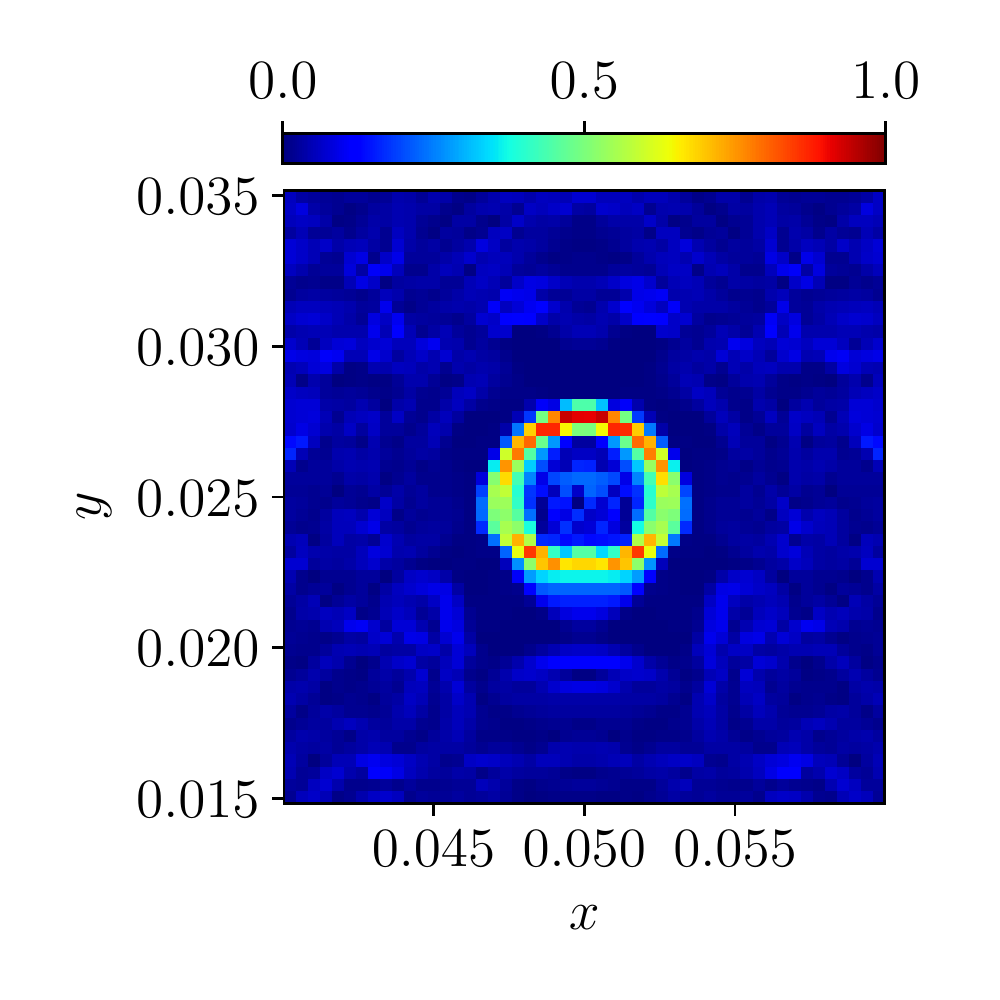}
		\caption{Constant Ansatz from~\cref{fig:constantAnsatz}.} 
	\end{subfigure}
	\caption{The $L_2$-norm of the spatial gradients of the recovered material distribution $||\nabla \hat{\gamma}(\boldsymbol{x})||_2$.} \label{fig:sharpness}
\end{figure}

The sharpness of the identified material distributions can be quantified by the $L_2$-norm of the spatial gradients of the recovered material, inspired by Sobel filters commonly used in image processing \cite{gonzalez_digital_2018}:
\begin{equation}
||\nabla \hat{\gamma}(\boldsymbol{x})||_2=\sqrt{\nabla_x \hat{\gamma}(\boldsymbol{x}) + \nabla_y \hat{\gamma}(\boldsymbol{x})} \label{eq:norm}
\end{equation}
A greater norm indicates a sharper edge. To highlight the improved sharpness recovery by the neural network Ansatz compared to the constant Ansatz, the norm from~\cref{eq:norm} around the circular void is illustrated in~\cref{fig:sharpness} for the two-dimensional case. At the cirular boundary, the norm is mostly larger in the neural network case showcasing a sharper inversion, which can be quantified by considering the mean of the non-zero values, which are extracted by selecting the points with a norm greater than a positive threshold, with the purpose of excluding the noise present in the constant Ansatz solution. For any threshold, the mean is larger for the neural network case, confirming an overall sharper recovered material distribution. 

\newpage

\bibliographystyle{ieeetr}



\setlength{\bibsep}{3pt}
\setlength{\bibhang}{0.75cm}{\fontsize{9}{9}\selectfont\bibliography{pinncomparison}}

\begin{thebibliography}{10}

\bibitem{lagaris_artificial_1998}
I.~Lagaris, A.~Likas, and D.~Fotiadis, ``Artificial neural networks for solving
  ordinary and partial differential equations,'' {\em IEEE Transactions on
  Neural Networks}, vol.~9, pp.~987--1000, Sept. 1998.

\bibitem{psichogios_hybrid_1992}
D.~C. Psichogios and L.~H. Ungar, ``A hybrid neural network-first principles
  approach to process modeling,'' {\em AIChE Journal}, vol.~38, pp.~1499--1511,
  Oct. 1992.

\bibitem{raissi_physics-informed_2019}
M.~Raissi, P.~Perdikaris, and G.~Karniadakis, ``Physics-informed neural
  networks: {A} deep learning framework for solving forward and inverse
  problems involving nonlinear partial differential equations,'' {\em Journal
  of Computational Physics}, vol.~378, pp.~686--707, Feb. 2019.

\bibitem{cuomo_scientific_2022}
S.~Cuomo, V.~S. Di~Cola, F.~Giampaolo, G.~Rozza, M.~Raissi, and F.~Piccialli,
  ``Scientific {Machine} {Learning} {Through} {Physics}–{Informed} {Neural}
  {Networks}: {Where} we are and {What}’s {Next},'' {\em Journal of
  Scientific Computing}, vol.~92, p.~88, July 2022.

\bibitem{karniadakis_physics-informed_2021}
G.~E. Karniadakis, I.~G. Kevrekidis, L.~Lu, P.~Perdikaris, S.~Wang, and
  L.~Yang, ``Physics-informed machine learning,'' {\em Nature Reviews Physics},
  vol.~3, pp.~422--440, June 2021.

\bibitem{Kollmannsberger2021a}
S.~Kollmannsberger, D.~D'Angella, M.~Jokeit, and L.~Herrmann, {\em
  Physics-{{Informed Neural Networks}}}, vol.~977, pp.~55--84.
\newblock {Cham}: {Springer International Publishing}, 2021.

\bibitem{paszke_pytorch_2019}
A.~Paszke, S.~Gross, F.~Massa, A.~Lerer, J.~Bradbury, G.~Chanan, T.~Killeen,
  Z.~Lin, N.~Gimelshein, L.~Antiga, A.~Desmaison, A.~Köpf, E.~Yang, Z.~DeVito,
  M.~Raison, A.~Tejani, S.~Chilamkurthy, B.~Steiner, L.~Fang, J.~Bai, and
  S.~Chintala, ``{PyTorch}: {An} {Imperative} {Style}, {High}-{Performance}
  {Deep} {Learning} {Library},'' Dec. 2019.
\newblock arXiv:1912.01703 [cs, stat].

\bibitem{tensorflow2015-whitepaper}
M.~Abadi, A.~Agarwal, P.~Barham, E.~Brevdo, Z.~Chen, C.~Citro, G.~S. Corrado,
  A.~Davis, J.~Dean, M.~Devin, S.~Ghemawat, I.~Goodfellow, A.~Harp, G.~Irving,
  M.~Isard, Y.~Jia, R.~Jozefowicz, L.~Kaiser, M.~Kudlur, J.~Levenberg,
  D.~Man\'{e}, R.~Monga, S.~Moore, D.~Murray, C.~Olah, M.~Schuster, J.~Shlens,
  B.~Steiner, I.~Sutskever, K.~Talwar, P.~Tucker, V.~Vanhoucke, V.~Vasudevan,
  F.~Vi\'{e}gas, O.~Vinyals, P.~Warden, M.~Wattenberg, M.~Wicke, Y.~Yu, and
  X.~Zheng, ``{TensorFlow}: Large-scale machine learning on heterogeneous
  systems,'' 2015.
\newblock Software available from tensorflow.org.

\bibitem{thuerey_deep_2020}
N.~Thuerey, K.~Weißenow, L.~Prantl, and X.~Hu, ``Deep {Learning} {Methods} for
  {Reynolds}-{Averaged} {Navier}–{Stokes} {Simulations} of {Airfoil}
  {Flows},'' {\em AIAA Journal}, vol.~58, pp.~25--36, Jan. 2020.

\bibitem{zhu_physics-constrained_2019}
Y.~Zhu, N.~Zabaras, P.-S. Koutsourelakis, and P.~Perdikaris,
  ``Physics-constrained deep learning for high-dimensional surrogate modeling
  and uncertainty quantification without labeled data,'' {\em Journal of
  Computational Physics}, vol.~394, pp.~56--81, Oct. 2019.

\bibitem{lino_simulating_2021}
M.~Lino, C.~Cantwell, A.~A. Bharath, and S.~Fotiadis, ``Simulating {Continuum}
  {Mechanics} with {Multi}-{Scale} {Graph} {Neural} {Networks},'' June 2021.
\newblock arXiv:2106.04900 [physics].

\bibitem{sanchez-gonzalez_learning_2020}
A.~Sanchez-Gonzalez, J.~Godwin, T.~Pfaff, R.~Ying, J.~Leskovec, and P.~W.
  Battaglia, ``Learning to {Simulate} {Complex} {Physics} with {Graph}
  {Networks},'' Sept. 2020.
\newblock arXiv:2002.09405 [physics, stat].

\bibitem{pfaff_learning_2021}
T.~Pfaff, M.~Fortunato, A.~Sanchez-Gonzalez, and P.~W. Battaglia, ``Learning
  {Mesh}-{Based} {Simulation} with {Graph} {Networks},'' June 2021.
\newblock arXiv:2010.03409 [cs].

\bibitem{bhatnagar_prediction_2019}
S.~Bhatnagar, Y.~Afshar, S.~Pan, K.~Duraisamy, and S.~Kaushik, ``Prediction of
  aerodynamic flow fields using convolutional neural networks,'' {\em
  Computational Mechanics}, vol.~64, pp.~525--545, Aug. 2019.

\bibitem{lu_learning_2021}
L.~Lu, P.~Jin, G.~Pang, Z.~Zhang, and G.~E. Karniadakis, ``Learning nonlinear
  operators via {DeepONet} based on the universal approximation theorem of
  operators,'' {\em Nature Machine Intelligence}, vol.~3, pp.~218--229, Mar.
  2021.

\bibitem{markidis_old_2021}
S.~Markidis, ``The {Old} and the {New}: {Can} {Physics}-{Informed}
  {Deep}-{Learning} {Replace} {Traditional} {Linear} {Solvers}?,'' {\em
  Frontiers in Big Data}, vol.~4, 2021.

\bibitem{leiteritz_how_2021}
R.~Leiteritz and D.~Pflüger, ``How to {Avoid} {Trivial} {Solutions} in
  {Physics}-{Informed} {Neural} {Networks},'' Dec. 2021.
\newblock arXiv:2112.05620 [cs, stat].

\bibitem{hughes_finite_2000}
T.~J.~R. Hughes, {\em The finite element method: linear static and dynamic
  finite element analysis}.
\newblock Mineola, NY: Dover Publications, 2000.

\bibitem{langtangen_finite_2017}
H.~P. Langtangen, {\em Finite difference computing with {PDEs}: a modern
  software approach}.
\newblock No.~vol. 16 in Texts in computational science and engineering, Cham,
  Switzerland: Springer Open, 2017.

\bibitem{han_solving_2018}
J.~Han, A.~Jentzen, and W.~E, ``Solving high-dimensional partial differential
  equations using deep learning,'' {\em Proceedings of the National Academy of
  Sciences}, vol.~115, pp.~8505--8510, Aug. 2018.

\bibitem{sirignano_dgm_2018}
J.~Sirignano and K.~Spiliopoulos, ``{DGM}: {A} deep learning algorithm for
  solving partial differential equations,'' {\em Journal of Computational
  Physics}, vol.~375, pp.~1339--1364, Dec. 2018.

\bibitem{goswami_physics-informed_2022}
S.~Goswami, A.~Bora, Y.~Yu, and G.~E. Karniadakis, ``Physics-{Informed} {Deep}
  {Neural} {Operator} {Networks},'' July 2022.
\newblock arXiv:2207.05748 [cs, math].

\bibitem{oldenburg_geometry_2022}
J.~Oldenburg, F.~Borowski, A.~Öner, K.-P. Schmitz, and M.~Stiehm, ``Geometry
  aware physics informed neural network surrogate for solving
  {Navier}–{Stokes} equation ({GAPINN}),'' {\em Advanced Modeling and
  Simulation in Engineering Sciences}, vol.~9, p.~8, June 2022.

\bibitem{wong_improved_2021}
J.~C. Wong, C.~Ooi, P.-H. Chiu, and M.~H. Dao, ``Improved {Surrogate}
  {Modeling} of {Fluid} {Dynamics} with {Physics}-{Informed} {Neural}
  {Networks},'' May 2021.
\newblock arXiv:2105.01838 [physics].

\bibitem{shukla_physics-informed_2020}
K.~Shukla, P.~C. Di~Leoni, J.~Blackshire, D.~Sparkman, and G.~E. Karniadakis,
  ``Physics-{Informed} {Neural} {Network} for {Ultrasound} {Nondestructive}
  {Quantification} of {Surface} {Breaking} {Cracks},'' {\em Journal of
  Nondestructive Evaluation}, vol.~39, p.~61, Sept. 2020.

\bibitem{rashtbehesht_physicsinformed_2022}
M.~Rasht‐Behesht, C.~Huber, K.~Shukla, and G.~E. Karniadakis,
  ``Physics‐{Informed} {Neural} {Networks} ({PINNs}) for {Wave} {Propagation}
  and {Full} {Waveform} {Inversions},'' {\em Journal of Geophysical Research:
  Solid Earth}, vol.~127, May 2022.

\bibitem{cai_flow_2021}
S.~Cai, Z.~Wang, F.~Fuest, Y.~J. Jeon, C.~Gray, and G.~E. Karniadakis, ``Flow
  over an espresso cup: inferring 3-{D} velocity and pressure fields from
  tomographic background oriented {Schlieren} via physics-informed neural
  networks,'' {\em Journal of Fluid Mechanics}, vol.~915, May 2021.

\bibitem{wang_deep_2021}
S.~Wang and P.~Perdikaris, ``Deep learning of free boundary and {Stefan}
  problems,'' {\em Journal of Computational Physics}, vol.~428, p.~109914, Mar.
  2021.

\bibitem{mao_physics-informed_2020}
Z.~Mao, A.~D. Jagtap, and G.~E. Karniadakis, ``Physics-informed neural networks
  for high-speed flows,'' {\em Computer Methods in Applied Mechanics and
  Engineering}, vol.~360, p.~112789, Mar. 2020.

\bibitem{jagtap_conservative_2020}
A.~D. Jagtap, E.~Kharazmi, and G.~E. Karniadakis, ``Conservative
  physics-informed neural networks on discrete domains for conservation laws:
  {Applications} to forward and inverse problems,'' {\em Computer Methods in
  Applied Mechanics and Engineering}, vol.~365, p.~113028, June 2020.

\bibitem{jagtap_physics-informed_2022}
A.~D. Jagtap, Z.~Mao, N.~Adams, and G.~E. Karniadakis, ``Physics-informed
  neural networks for inverse problems in supersonic flows,'' {\em Journal of
  Computational Physics}, vol.~466, p.~111402, Oct. 2022.
\newblock arXiv:2202.11821 [cs, math].

\bibitem{chen_physics-informed_2020}
Y.~Chen, L.~Lu, G.~E. Karniadakis, and L.~Dal~Negro, ``Physics-informed neural
  networks for inverse problems in nano-optics and metamaterials,'' {\em Optics
  Express}, vol.~28, p.~11618, Apr. 2020.

\bibitem{givoli_tutorial_2021}
D.~Givoli, ``A tutorial on the adjoint method for inverse problems,'' {\em
  Computer Methods in Applied Mechanics and Engineering}, vol.~380, p.~113810,
  July 2021.

\bibitem{plessix_review_2006}
R.-E. Plessix, ``A review of the adjoint-state method for computing the
  gradient of a functional with geophysical applications,'' {\em Geophysical
  Journal International}, vol.~167, pp.~495--503, Nov. 2006.

\bibitem{fichtner_full_2011}
A.~Fichtner and F.~Bleibinhaus, {\em Full seismic waveform modelling and
  inversion}.
\newblock Advances in geophysical and environmental mechanics and mathematics,
  Berlin Heidelberg: Springer, 2011.

\bibitem{sayag_shape_2022}
A.~Sayag and D.~Givoli, ``Shape identification of scatterers {Using} a
  time-dependent adjoint method,'' {\em Computer Methods in Applied Mechanics
  and Engineering}, vol.~394, p.~114923, May 2022.

\bibitem{buerchner_2022}
T.~Bürchner, P.~Kopp, S.~Kollmannsberger, and E.~Rank, ``Immersed boundary
  parametrizations for full waveform inversion,'' {\em Computer Methods in
  Applied Mechanics and Engineering}, vol.~406, p.~115893, Mar. 2023.

\bibitem{baydin_automatic_2018}
A.~G. Baydin, B.~A. Pearlmutter, A.~A. Radul, and J.~M. Siskind, ``Automatic
  differentiation in machine learning: a survey,'' Feb. 2018.
\newblock arXiv:1502.05767 [cs, stat].

\bibitem{geneva_modeling_2020}
N.~Geneva and N.~Zabaras, ``Modeling the {Dynamics} of {PDE} {Systems} with
  {Physics}-{Constrained} {Deep} {Auto}-{Regressive} {Networks},'' {\em Journal
  of Computational Physics}, vol.~403, p.~109056, Feb. 2020.
\newblock arXiv: 1906.05747.

\bibitem{wang_residual_2017}
F.~Wang, M.~Jiang, C.~Qian, S.~Yang, C.~Li, H.~Zhang, X.~Wang, and X.~Tang,
  ``Residual {Attention} {Network} for {Image} {Classification},'' in {\em 2017
  {IEEE} {Conference} on {Computer} {Vision} and {Pattern} {Recognition}
  ({CVPR})}, pp.~6450--6458, July 2017.
\newblock ISSN: 1063-6919.

\bibitem{zhang_occluded_2018}
S.~Zhang, J.~Yang, and B.~Schiele, ``Occluded {Pedestrian} {Detection}
  {Through} {Guided} {Attention} in {CNNs},'' in {\em 2018 {IEEE}/{CVF}
  {Conference} on {Computer} {Vision} and {Pattern} {Recognition}},
  pp.~6995--7003, June 2018.
\newblock ISSN: 2575-7075.

\bibitem{nandwani_primal-dual_2019}
Y.~Nandwani, A.~Pathak, {Mausam}, and P.~Singla, ``A primal-dual formulation
  for deep learning with constraints,'' in {\em Proceedings of the 33rd
  {International} {Conference} on {Neural} {Information} {Processing}
  {Systems}}, no.~1091, pp.~12179--12190, Red Hook, NY, USA: Curran Associates
  Inc., Dec. 2019.

\bibitem{mcclenny_self-adaptive_2022}
L.~McClenny and U.~Braga-Neto, ``Self-{Adaptive} {Physics}-{Informed} {Neural}
  {Networks} using a {Soft} {Attention} {Mechanism},'' Apr. 2022.
\newblock arXiv:2009.04544 [cs, stat].

\bibitem{kingma_adam_2017}
D.~P. Kingma and J.~Ba, ``Adam: {A} {Method} for {Stochastic} {Optimization},''
  Jan. 2017.
\newblock arXiv:1412.6980 [cs].

\bibitem{liu_limited_1989}
D.~C. Liu and J.~Nocedal, ``On the limited memory {BFGS} method for large scale
  optimization,'' {\em Mathematical Programming}, vol.~45, pp.~503--528, Aug.
  1989.

\bibitem{samaniego_energy_2020}
E.~Samaniego, C.~Anitescu, S.~Goswami, V.~Nguyen-Thanh, H.~Guo, K.~Hamdia,
  X.~Zhuang, and T.~Rabczuk, ``An energy approach to the solution of partial
  differential equations in computational mechanics via machine learning:
  {Concepts}, implementation and applications,'' {\em Computer Methods in
  Applied Mechanics and Engineering}, vol.~362, p.~112790, Apr. 2020.

\bibitem{moseley_solving_2020}
B.~Moseley, A.~Markham, and T.~Nissen-Meyer, ``Solving the wave equation with
  physics-informed deep learning,'' June 2020.
\newblock arXiv:2006.11894 [physics].

\bibitem{song_solving_2021}
C.~Song, T.~Alkhalifah, and U.~B. Waheed, ``Solving the frequency-domain
  acoustic {VTI} wave equation using physics-informed neural networks,'' {\em
  Geophysical Journal International}, vol.~225, pp.~846--859, Feb. 2021.

\bibitem{song_solving_2021-1}
C.~Song, T.~Alkhalifah, and U.~b. Waheed, ``Solving the acoustic {VTI} wave
  equation using physics-informed neural networks,'' {\em Geophysical Journal
  International}, vol.~225, pp.~846--859, Feb. 2021.
\newblock arXiv:2008.01865 [physics].

\bibitem{karimpouli_physics_2020}
S.~Karimpouli and P.~Tahmasebi, ``Physics informed machine learning: {Seismic}
  wave equation,'' {\em Geoscience Frontiers}, vol.~11, pp.~1993--2001, Nov.
  2020.

\bibitem{rasht-behesht_physics-informed_2021}
M.~Rasht-Behesht, C.~Huber, K.~Shukla, and G.~E. Karniadakis,
  ``Physics-{Informed} {Deep} {Learning} for {Wave} {Propagation} and {Full}
  {Waveform} {Inversions},'' 2021.

\bibitem{michea_accelerating_2010}
D.~Michéa and D.~Komatitsch, ``Accelerating a three-dimensional
  finite-difference wave propagation code using {GPU} graphics cards:
  {Accelerating} a wave propagation code using {GPUs},'' {\em Geophysical
  Journal International}, pp.~no--no, May 2010.

\bibitem{norgaard_applications_2017}
S.~A. Nørgaard, M.~Sagebaum, N.~R. Gauger, and B.~S. Lazarov, ``Applications
  of automatic differentiation in topology optimization,'' {\em Structural and
  Multidisciplinary Optimization}, vol.~56, pp.~1135--1146, Nov. 2017.

\bibitem{dilgen_topology_2018}
C.~B. Dilgen, S.~B. Dilgen, D.~R. Fuhrman, O.~Sigmund, and B.~S. Lazarov,
  ``Topology optimization of turbulent flows,'' {\em Computer Methods in
  Applied Mechanics and Engineering}, vol.~331, pp.~363--393, Apr. 2018.

\bibitem{chen_neural_2019}
R.~T.~Q. Chen, Y.~Rubanova, J.~Bettencourt, and D.~Duvenaud, ``Neural
  {Ordinary} {Differential} {Equations},'' Dec. 2019.
\newblock arXiv:1806.07366 [cs, stat].

\bibitem{williams_gradient-based_1995}
R.~J. Williams and D.~Zipser, ``Gradient-based learning algorithms for
  recurrent networks and their computational complexity,'' in {\em
  Backpropagation: theory, architectures, and applications}, pp.~433--486, USA:
  L. Erlbaum Associates Inc., Jan. 1995.

\bibitem{sutskever_training_2013}
I.~Sutskever, {\em Training recurrent neural networks}.
\newblock phd, University of Toronto, CAN, 2013.
\newblock AAINS22066 ISBN-13: 9780499220660.

\bibitem{chandrasekhar_tounn_2021}
A.~Chandrasekhar and K.~Suresh, ``{TOuNN}: {Topology} {Optimization} using
  {Neural} {Networks},'' {\em Structural and Multidisciplinary Optimization},
  vol.~63, pp.~1135--1149, Mar. 2021.

\bibitem{wirgin_inverse_2004}
A.~Wirgin, ``The inverse crime,'' Jan. 2004.
\newblock arXiv:math-ph/0401050.

\bibitem{parvizian_finite_2007}
J.~Parvizian, A.~Düster, and E.~Rank, ``Finite cell method: h- and p-extension
  for embedded domain problems in solid mechanics,'' {\em Computational
  Mechanics}, vol.~41, pp.~121--133, Sept. 2007.

\bibitem{hicks_arbitrary_2002}
G.~J. Hicks, ``Arbitrary source and receiver positioning in finite‐difference
  schemes using {Kaiser} windowed sinc functions,'' {\em GEOPHYSICS}, vol.~67,
  pp.~156--165, Jan. 2002.

\bibitem{hug_three-field_2022}
L.~Hug, M.~Potten, G.~Stockinger, K.~Thuro, and S.~Kollmannsberger, ``A
  three-field phase-field model for mixed-mode fracture in rock based on
  experimental determination of the mode {II} fracture toughness,'' {\em
  Engineering with Computers}, July 2022.

\bibitem{karras_progressive_2018}
T.~Karras, T.~Aila, S.~Laine, and J.~Lehtinen, ``Progressive {Growing} of
  {GANs} for {Improved} {Quality}, {Stability}, and {Variation},'' Feb. 2018.
\newblock arXiv:1710.10196 [cs, stat].

\bibitem{jagtap_adaptive_2020}
A.~D. Jagtap, K.~Kawaguchi, and G.~E. Karniadakis, ``Adaptive activation
  functions accelerate convergence in deep and physics-informed neural
  networks,'' {\em Journal of Computational Physics}, vol.~404, p.~109136, Mar.
  2020.

\bibitem{he_delving_2015}
K.~He, X.~Zhang, S.~Ren, and J.~Sun, ``Delving {Deep} into {Rectifiers}:
  {Surpassing} {Human}-{Level} {Performance} on {ImageNet} {Classification},''
  Feb. 2015.
\newblock arXiv:1502.01852 [cs] version: 1.

\bibitem{Goodfellow-et-al-2016}
I.~Goodfellow, Y.~Bengio, and A.~Courville, {\em Deep Learning}.
\newblock MIT Press, 2016.
\newblock \url{http://www.deeplearningbook.org}.

\bibitem{pascanu_difficulty_2013}
R.~Pascanu, T.~Mikolov, and Y.~Bengio, ``On the difficulty of training
  {Recurrent} {Neural} {Networks},'' Feb. 2013.
\newblock arXiv:1211.5063 [cs].

\bibitem{zhang_why_2020}
J.~Zhang, T.~He, S.~Sra, and A.~Jadbabaie, ``Why gradient clipping accelerates
  training: {A} theoretical justification for adaptivity,'' Feb. 2020.
\newblock arXiv:1905.11881 [cs, math].

\bibitem{kumar_weight_2017}
S.~K. Kumar, ``On weight initialization in deep neural networks,'' May 2017.
\newblock arXiv:1704.08863 [cs].

\bibitem{fichtner_source_2011}
A.~Fichtner, ``Source {Stacking} {Data} {Reduction} for {Full} {Waveform}
  {Tomography} at the {Global} {Scale},'' in {\em Full {Seismic} {Waveform}
  {Modelling} and {Inversion}} (A.~Fichtner, ed.), Advances in {Geophysical}
  and {Environmental} {Mechanics} and {Mathematics}, pp.~281--299, Berlin,
  Heidelberg: Springer, 2011.

\bibitem{herrmann_code}
L.~Herrmann, T.~Bürchner, F.~Dietrich, and S.~Kollmannsberger, ``On the use of
  neural networks for full waveform inversion [software],'' {\em Mendeley
  Data}, 2023.

\bibitem{gonzalez_digital_2018}
R.~C. Gonzalez and R.~E. Woods, {\em Digital image processing}.
\newblock Pearson, fourth edition, global edition~ed.

\end{thebibliography}

\end{document}